\documentclass{article}
\usepackage{amsmath,amsthm,amsfonts,amscd,amssymb,eucal,latexsym}

\def\rest{\upharpoonright}

\def\Bbb{\mathbb}

\def\ov{\overline}
\def\cA{{\cal A}}
\def\cB{{\cal B}}
\def\cC{{\cal C}}
\def\cD{{\cal D}}

\def\cF{{\cal F}}

\def\cH{{\cal H}}

\def\cJ{{\cal J}}

\def\cS{{\cal S}}
\def\cT{{\cal T}}

\begin{document}

\title{A duality theorem for ergodic actions of compact quantum groups on
$C^*$--algebras}
\author{Claudia Pinzari\\
Dipartimento di Matematica, Universit\`a di Roma La Sapienza\\
00185--Roma, Italy\\ \\
John E. Roberts\\
Dipartimento di Matematica, Universit\`a di Roma Tor Vergata\\
00133--Roma, Italy}
\date{}
\maketitle

\begin{abstract} 
The spectral functor of an ergodic action of a compact
quantum group $G$ on a unital $C^*$--algebra is quasitensor, in
the sense that the
tensor product of two spectral subspaces
is isometrically contained in 
the spectral subspace of the tensor product representation, and the
inclusion maps satisfy natural properties. We show that any
quasitensor $^*$--functor from $\text{Rep}(G)$
to the category
of Hilbert spaces is the spectral functor of an
ergodic action of
$G$ on a unital $C^*$--algebra. 

As an application, we associate an
ergodic $G$--action on a unital $C^*$--algebra to an
inclusion of 
$\text{Rep}(G)$ into an abstract tensor $C^*$--category ${\cal T}$.

If
the inclusion arises from a quantum subgroup $K$ of $G$, the
associated 
$G$--system is just the quotient space $K\backslash G$.
If $G$ is a group and ${\cal T}$ has permutation symmetry, the
associated 
$G$--system
is commutative, and therefore  isomorphic to the classical quotient
space
by a subgroup of $G$.

If a tensor $C^*$--category  has a Hecke symmetry making an
object $\rho$ of 
dimension $d$ and $\mu$--determinant $1$ then there is an ergodic action
of
$S_\mu U(d)$ on a unital $C^*$--algebra having the $(\iota,\rho^r)$ as its
spectral subspaces. The special case of $S_\mu U(2)$ is discussed.

  \end{abstract}

\begin{section} {Introduction}
A theorem in \cite{DRinventiones}  asserts that any abstract tensor
$C^*$--category 
with conjugates and permutation symmetry is the representation category 
of a unique compact group, thus generalizing the
classical
Tannaka--Krein duality theorem, where one starts from a subcategory
of the category of Hilbert spaces (see, e.g., \cite{K}).

The content of this paper fits into the program of generalizing
the abstract duality theorem of \cite{DRinventiones} to tensor
$C^*$--categories
without permutation symmetry. Our interest in this problem 
is motivated by the fact that tensor $C^*$--categories 
 with conjugates, but also with a unitary symmetry of the braid group,
arise from low dimensional QFT \cite{Haag}.

In \cite{Wcmp}, \cite{Wtk} and \cite{Wlh} Woronowicz introduced  
 compact quantum groups, and generalized the classical
representation
theory of compact groups. He proved a Tannaka--Krein
duality theorem, asserting that the representation categories of compact
quantum groups are precisely the tensor $^*$--subcategories of categories 
of Hilbert spaces where every object has a conjugate. This theorem 
allowed him to construct the quantum deformations $S_\mu U(d)$ of the 
classical $SU(d)$ groups by  real a parameter $\mu$.

Therefore 
if a tensor $C^*$--category ${\cal T}$ with conjugates can be embedded
into a  category of Hilbert spaces, then ${\cal T}$ is necessarily the
representation
category of a compact quantum group.

In \cite{P1} the first named author characterized the representation
category of  $S_\mu U(d)$ among braided tensor
$C^*$--categories ${\cal T}$ with conjugates. 
It follows that if an object $\rho$ has a
 symmetry of the Hecke algebra type $H_\infty(\mu^2)$ making
$\rho$ of dimension $d$ and with $\mu$--determinant one, then there is a faithful tensor $^*$--functor
$\text{Rep}(S_\mu U(d))\to{\cal T}$.

The notion of  subgroup of a compact quantum group $G$ was
given by Podles in
\cite{Podles}, who computed all the subgroups of the  
quantum $SU(2)$ and $SO(3)$  groups. In the same paper the author introduced
 quantum quotient spaces, and he showed that these spaces have an
action of $G$ which splits into the direct sum of irreducibles, with
multiplicity bounded above by the Hilbert space dimension.

Later Wang proved in \cite{Wang} that compact quantum group
actions on quotient spaces are 
ergodic, as in the classical case, and he found an example of an ergodic
action on a commutative $C^*$--algebra which is not a quotient action.

It turns out that a compact quantum subgroup $K$ of $G$ gives rise to
an inclusion of tensor $C^*$--categories
$\text{Rep}(G)\to\text{Rep}(K)$, and that, by Tannaka--Krein duality, 
every tensor $^*$--inclusion of
$\text{Rep}(G)$ into a subcategory  of the  category of Hilbert
spaces,
taking a representation $u$ to its Hilbert space and acting trivially on
the arrows, is of this form.

In the case where a tensor $^*$--inclusion 
$\rho:\text{Rep}(G)\to{\cal T}$ into an abstract tensor $C^*$--category is
given, it is natural to look for a tensor $^*$--functor of ${\cal T}$
into the category of Hilbert spaces, acting as the embedding functor
on $\text{Rep}(G)$.  This  amounts to asking
 whether $\rho$
arises as an inclusion
associated with a quantum subgroup of $G$.

The aim of this paper is twofold. Assume that we have a tensor
$^*$--inclusion $\rho: \text{Rep}(G)\to{\cal T}$. The first result is the
construction 
of a unital $C^*$--algebra ${\cal B}$ associated with the inclusion,
and an ergodic action of 
$G$ on ${\cal B}$ whose spectral subspaces are the spaces $(\iota,\rho_u)$,
where $u$ varies in the set of unitary irreducible representations of $G$. 
The relevance of this construction to
 abstract duality for compact quantum groups will be discussed elsewhere
\cite{DPR}. 
We exhibit
 two interesting particular cases of this construction.

If  ${\cal T}=\text{Rep}(K)$, with $K$ 
a compact quantum subgroup of $G$, the ergodic system
thus
obtained is just the quotient space $K\backslash G$ (Theorem 11.1). If
instead $G$ is a group
and ${\cal T}$ has a permutation symmetry, then ${\cal B}$ turns out to be 
commutative, and can therefore be identified with the $C^*$--algebra 
of continuous functions over a quotient of $G$ by a point stabilizer subgroup
(Theorem 10.2). 

We apply our construction to the case where an abstract tensor
$C^*$--category ${\cal T}$ has an object of dimension $d$ and
$\mu$--determinant one  and we find ergodic $C^*$--actions 
of $S_\mu U(d)$ (Theorem 10.3). We also discuss the particular case where
$d=2$ (Cor.\
10.5).

Our second aim is to characterize, among all ergodic actions,
 those which are isomorphic to the quotient spaces by some
quantum subgroup. 
By a well known theorem by H\o egh--Krohn, Landstad and St\o rmer
\cite{HLS}, an irreducible representation of a compact group $G$
appears in the spectrum of an ergodic action of $G$
on a unital
$C^*$--algebra with multiplicity bounded above by its dimension.

As a generalization to compact quantum groups, 
Boca shows  that the 
multiplicity of a unitary irreducible representation of the quantum group 
in the action
is, instead,  bounded above by its quantum dimension \cite{Boca}.

In \cite{BRV}, Bichon, De Rijdt and Vaes construct examples of ergodic
actions of $S_\mu U(2)$ in which the multiplicity of the fundamental
representation can be any integer  $n$ with $2\leq n\leq
\mu+\frac{1}{\mu}$. Therefore
these
actions are not on quotient spaces of $S_\mu U(2)$. 
They also give a simpler proof of Boca's result by introducting
a new invariant, the {\it quantum multiplicity} of a spectral
representation, which they show to be
bounded below by the multiplicity
and above by the quantum dimension.
The main tool for constructing their 
examples is a duality theorem, proved in that paper,
  between ergodic quantum actions for which the quantum multiplicity
equals the quantum dimension,
and certain maps  associating to any 
irreducible representation of $G$ a finite dimensional Hilbert space, and 
to any  intertwining operator between tensor products of irreducible
representations, 
a linear map between the tensor products of the corresponding associated
Hilbert spaces, respecting  composition,
tensor products and the $^*$--involution.

In our generalization to the case where the quantum multiplicity
is not maximal, our main
tool is the 
{\it
spectral functor} associated with a generic  ergodic action (Sect.\ 7).
This
is the covariant $^*$--functor that
 associates to {\it any} unitary, finite dimensional
representation
$u$ of $G$, 
the dual space $\overline{L}_u$ of the space all multiplets in ${\cal B}$
transforming like
$u$ under the
action $\eta$. The space $\overline{L}_u$, by ergodicity, is known to
become
a
Hilbert
space in a
natural way.

We stress that the functor $\overline{L}$ satisfies two crucial
properties: the first
one is  that 
$\overline{L}_{u\otimes v}$ naturally contains a copy of
$\overline{L}_u\otimes\overline{L}_v$, in such a way that the copy of 
$\overline{L}_u\otimes \overline{L}_v\otimes \overline{L}_z$ contained in
both 
$\overline{L}_{u\otimes v}\otimes\overline{L}_z$ and
$\overline{L}_u\otimes\overline{L}_{v\otimes z}$ is the same. The second
property is that
the projection from $\overline{L}_{u\otimes v\otimes z}$ to
$\overline{L}_{u\otimes v}\otimes \overline{L}_z$ actually takes 
$\overline{L}_{u}\otimes\overline{L}_{v\otimes z}$ onto 
$\overline{L}_u\otimes\overline{L}_v\otimes\overline{L}_z$ (Theorem 7.3).

We call any functor ${\cal F}$ from a generic  tensor $C^*$--category ${\cal T}$
to the
category of Hilbert spaces satisfying the above properties,
{\it quasitensor} (Sect.\ 3). We show that
quasitensor functors, like the tensor ones, have the property that if 
 $\overline{\rho}$ is a conjugate 
of $\rho$ then the Hilbert space ${\cal F}(\overline{\rho})$ must be a
conjugate of ${\cal F}(\rho)$, although this conjugate must be found
in the image category of ${\cal F}$ enriched with the projection maps
from the spaces 
${\cal F}(\rho\otimes\sigma)$ onto ${\cal F}(\rho)\otimes{\cal F}(\sigma)$
(Theorem 3.7).

This result easily shows that
 ${\cal F}(\rho)$ is automatically
finite dimensional and endowed with an {\it intrinsic dimension}, in the
sense of \cite{LongoRoberts}, bounded below 
by the Hilbert space dimension of ${\cal F}(\rho)$ and above by the
intrinsic dimension of
$\rho$ (Cor.\ 3.8).

This result applied to the spectral functor $\ov{L}$
 allows us to identify the  quantum
multiplicity of a spectral irreducible representation of an ergodic
action of  \cite{BRV}, with the intrinsic dimension of
$\ov{L}_u$, and to
recover the multiplicity bound  theorems of \cite{Boca} and
\cite{BRV}. The maximal quantum multiplicity  case corresponds to the case
where
for any irreducible $u$ with conjugate $\overline{u}$,
$\ov{L}_{\overline{u}}$ is
already a conjugate of 
$\ov{L}_{u}$ in the image of $\ov{L}$ (see Cor.\ 7.5 for a precise
statement).

The spectral functor determines uniquely the $^*$--algebra
structure of the dense $^*$--subalgebra of spectral elements.
Our main result is
a duality theorem for ergodic  $C^*$--actions
of compact quantum groups,
 showing that  
 any quasitensor $^*$--functor ${\cal
F}$ from the
representation category of a compact quantum group $G$ to the category
of Hilbert spaces is the spectral functor
of  an ergodic $G$--action 
over a unital $C^*$--algebra
(Theorem 9.1).

 Isomorphisms between two constructed ergodic systems 
correspond bijectively to unitary natural transformations between the corresponding
functors splitting as tensor products over a tensor product subspace
(Prop.\ 9.4).

Our main application of the Duality Theorem 9.1 is to inclusions of
tensor
$C^*$--categories.
In fact, quasitensor $^*$--functors defined on the representation
category of a compact quantum group $G$ arise very naturally 
from tensor $^*$--functors $\rho:\text{Rep}(G)\to {\cal T}$:
just take the map associating with a unitary $G$--representation $u$
the Hilbert space $(\iota,\rho_u)$, and with an intertwiner $T\in(u,v)$
between two representations, 
 the map acting on $(\iota,\rho_u)$ as left composition with $\rho(T)$
(Example 3.5).
Therefore
for any tensor $^*$--functor $\rho:\text{Rep}(G)\to{\cal T}$,
our duality theorem 
provides us with  an ergodic $G$--action over a
unital $C^*$--algebra, having the spaces $(\iota,\rho_u)$ as its
 spectral subspaces (Theorem 10.1).

Among quasitensor functors, those arising from quantum quotient spaces 
share the  property of being subobjects of the embedding functor
$H$ associating to any representation  $u$ its Hilbert space $H_u$.
In fact, for these functors, one can find a natural unitary transformation 
identifying the spectral subspace $\ov{L}_u$ with the subspace $K_u$
of $H_u$ of all vectors invariant under the restriction of $u$ to the
subgroup (Theorem 7.7).  For
maximal compact quantum groups, we
characterize algebraically 
the spaces of invariant vectors $K_u$ for a unique maximal subgroup $K$
(Theorem 5.5).

From this we derive our second main 
result.  In order that a maximal ergodic action $({\cal
B},\eta)$ be isomorphic to a compact quantum quotient space, it is
necessary and sufficient that there
exists, for any representation $u$ of $G$, an isometry from the
spectral subspace $\ov{L}_u$ onto some subspace $K_u$ of the 
representation Hilbert space, satisfying
certain coherence properties with
the tensor products (Theorem 11.3).
\end{section}

\begin{section}{Preliminaries}

\noindent{\it 2.1 Compact quantum groups and their representations}
\bigskip

In this paper  $G=({\cal A},\Delta)$ will always denote  a compact quantum
group in the sense of
\cite{Wlh}:
a unital $C^*$--algebra ${\cal A}$   with a 
unital $^*$--homomorphism (the coproduct)
$\Delta:{\cal A}\to{\cal A}\otimes{\cal A}$ such that
\begin{description}
\item a) $\Delta\otimes\iota\circ\Delta=\iota\otimes\Delta\circ\Delta$,
with $\iota$ the identity map on ${\cal A}$,
\item b) the sets
$\{b\otimes I\Delta(c), b,c\in{\cal A}\}$
$\{I\otimes b\Delta(c), b,c\in{\cal A}\}$
both span dense subspaces of ${\cal A}\otimes{\cal A}$.
\end{description}

Let $H$ be a finite dimensional Hilbert space, and form the free right
${\cal A}$--module $H\otimes{\cal A}$. 
The natural ${\cal A}$--valued inner product makes it into a
right Hilbert ${\cal A}$--module. 

A  unitary representation of $G$
with Hilbert space  $H_u$ can be defined as a  
linear map $u:H_u\to H_u\otimes{\cal
A}$,
with $H_u$ a finite dimensional Hilbert space, such that
$$(u(\psi), u(\phi))=(\psi,\phi) I,\quad \psi,\phi\in H_u,$$ 
$$u\otimes\iota\circ u=\iota\otimes\Delta\circ u,$$
$$u(H_u)I\otimes{\cal A}\text{ is total in } H_u\otimes{\cal A}.$$
If $u$ is a unitary representation, 
the matrix coefficients of $u$ 
are the elements of ${\cal A}$: 
$$u_{\phi,\psi}:=\ell_{\phi}^*\circ u(\psi),\quad \psi,\phi\in H_u,$$
with $\ell_\phi: {\cal A}\to H\otimes{\cal A}$  the operator
of tensoring on the left by $\phi$.
Let $u:H_u\to H_u\otimes{\cal A}$ be any linear map, and let
$(\phi_i)$ be an
orthonormal basis of $H_u$. Consider the matrix $(u_{ij})$ with entries
in ${\cal A}$, where $u_{ij}:=\ell_{\phi_i}^*u(\phi_j)$. Then $u$ is a
unitary representation if and only if the matrix $(u_{ij})$ is unitary
and $\Delta(u_{ij})=\sum_k u_{ik}\otimes u_{kj}$.

The linear span ${\cal A}_\infty$ of all the matrix coefficients is known 
to be a unital dense $^*$--subalgebra of ${\cal A}$, and a Hopf
$^*$--algebra 
$G_\infty=({\cal A}_\infty,\Delta)$
with the restricted coproduct \cite{Wlh}.

The category $\text{Rep}(G)$ of unitary  representations of $G$ 
with arrows
the spaces
$(u, v)$ of linear maps $T: H_u\to H_v$ intertwining $u$ and $v$:
$T\otimes I\circ u=v\circ T$
is a
tensor $C^*$--category with conjugates \cite{Wtk}.
Recall that the tensor product $u\otimes v$  of two representations $u$
and $v$ and the conjugate representation $\overline{u}$ 
are those representations 
with Hilbert spaces
$H_u\otimes H_v$ and $H_{\overline{u}}$  and
coefficients, 
$$(u\otimes
v)_{\phi\otimes\phi',\psi\otimes\psi'}:=u_{\phi,\psi}v_{\phi',\psi'},\quad
\phi,\psi\in H_u, \phi',\psi'\in H_v,$$
$${\overline u}_{j\phi, {j^*}^{-1}\psi}:=
({u}_{\phi,\psi})^*,\quad \phi,\psi\in H_u,$$
respectively,
where $j: H_u\to H_{\overline{u}}$ is an antilinear invertible 
intertwiner.
\bigskip

\noindent{\it 2.2 Spectra of Hopf $^*$--algebra actions}
\bigskip

Let ${\cal C}$ be a unital $^*$--algebra and $G=({\cal A},\Delta)$ a
compact quantum group.
Consider the dense Hopf $^*$--subalgebra 
$G_\infty=({\cal A}_\infty,\Delta)$ of 
$G$ and a unital  action  
$\eta:{\cal C}\to{\cal C}\odot {\cal A}_\infty$ of it 
on ${\cal C}$ (where $\odot$ denotes the algebraic tensor product): $\eta$ is a unital
$^*$--homomorphism such that $\eta\otimes\iota\circ\eta=\iota\otimes\Delta\circ\eta$.

We define the ${\it spectrum}$ of $\eta$, denoted $sp(\eta)$, to be the
set 
of all unitary $G$--representations $u: H_u\to H_u\otimes{\cal A}$
for 
which there is a faithful linear map $T: H_u\to{\cal C}$ intertwining $u$
with $\eta$:
$$\eta\circ T=T\otimes \iota \circ u.$$

In other words, representing $u$ as a matrix $u=(u_{ij})$ with respect to
some orthonormal basis of $H_u$, $u\in sp(\eta)$ if and only if there
exists
a multiplet $(c_1,\dots, c_d)$, with $d$ the dimension of $u$, constituted
by linearly independent elements of ${\cal C}$, such that
$$\eta(c_i)=\sum_j c_j\otimes u_{ji}.$$ 
For compact quantum groups this notion was  introduced
by Podles in
\cite{Podles},  as a generalization of
the classical notion
for an action of a compact group on a $C^*$--algebra \cite{GLR}.
We show some simple properties of the spectrum.\medskip

\noindent{\bf 2.1 Proposition} {\sl
\begin{description} 
\item{\rm a)} If $u\in sp(\eta)$ and if
 $z$ is a unitary representation of $G$ such that $(z, u)$ contains an
isometry, 
then $z\in sp(\eta)$,
\item{\rm b)} if $u\in sp(\eta)$ and $\overline{u}$ is a unitary
representation
equivalent to the complex conjugate ${u}_*$ then
$\overline{u}\in sp(\eta)$. Here $u_*$ denotes the representation, in general not 
unitary, whose matrix elements are the adjoints of those of $u$.
\end{description}
}\medskip

\noindent{\bf Proof} a) If $S\in(z, u)$ is an isometry and $T:H_u\to{\cal
C}$
is a faithful intertwining map then $T\circ S: H_z\to {\cal C}$ is a
faithful intertwining map:
$$\eta\circ(TS)=(\eta\circ T)S=(T\otimes\iota)\circ u S=(TS)\otimes\iota
z.$$
b) If $(c_1,\dots, c_d)$ is a linearly independent multiplet transforming
like the unitary representation $u$, then $(c_1^*,\dots, c_d^*)$ is a
linearly independent multiplet transforming like the complex conjugate 
representation $u_*$, which, in general, is just an invertible
representation, but equivalent to a unitary representation
\cite{Wcmp}.
Let
$\mu=(\mu_{ij})\in(\overline{u}, u_*)$ be an invertible intertwiner with 
a unitary representation $\overline{u}$.
Set
$f_j:=\sum_p\mu_{pj}c_p^*$. The multiplet $(f_1,\dots, f_d)$ is linearly
independent since $\mu$ is invertible, and transforms like $\overline{u}$:
$$\eta(f_j)=\sum_{p,r}\mu_{pj}c_r^*\otimes u_{rp}^*=$$
$$\sum_{r} c_r^*\otimes ({u}_*\mu)_{r,j}=\sum_r c_r^*\otimes
(\mu
\overline{u})_{r,j}=$$
$$\sum_{r,s}\mu_{rs}c_r^*\otimes \overline{u}_{s,j}=\sum_s f_s\otimes
\overline{u}_{s,j}.$$

We denote  the linear span of all spectral multiplets by ${\cal C}_{sp}$. 
Part a) of the previous proposition tells us 
that ${\cal C}_{sp}$ is generated, as a linear space, by those
nonzero multiplets 
transforming according to unitary irreducible $G$--representations
(such multiplets are automatically linearly independent by
irreducibility).
\medskip

\noindent{\bf 2.2 Proposition} {\sl
\begin{description}
\item{\rm a)} If $T: H_u\to {\cal C}$ is any linear map
satisfying  $\eta\circ T=T\otimes \iota\circ u$
then the image of $T$ lies in ${\cal C}_{sp}$.
\item{\rm b)} 
${\cal C}_{sp}$ is a unital $^*$--subalgebra of ${\cal C}$ invariant under
the
$G_\infty$--action: $\eta({\cal C}_{sp})\subset{\cal C}_{sp}\odot{\cal
A}_\infty$.
\end{description}}\medskip

\noindent{\bf Proof} a) We can assume $T\neq0$.
Let us write $u$ as a direct sum of unitary irreducible 
subrepresentations $u_1,\dots, u_p$.
Let $H_i$ be the subspace
corresponding
to $u_i$, for $i=1,\dots,p$.
Since each $u_i$ is irreducible, $T$ is either faithful or zero on $H_i$.
We can assume that there is a $q\leq p$, $q\geq1$ such that the restriction
$T_i$
of $T$ to $H_i$ is faithful 
for $i=1,\dots, q$ and $T=0$ on $H_i$ for $i>q$. 
Thus any element $T\psi$ in the image of $T$ can be written
as a sum $\sum_{k\leq q} T_k\psi_k$ with  $T_k\psi_k\in{\cal C}_{sp}$ now.  

b)The trivial representation is clearly in the spectrum 
with spectral subspace ${\mathbb C}I$, so ${\cal C}_{sp}$ contains the
identity.
 If
$T: H_u\to{\cal C}$ and $S:H_v\to{\cal C}$ are faithful maps 
intertwining $u$ and $v$ respectively with $\eta$, then
the map 
$$H_u\otimes H_v\to{\cal C}, \quad\psi\otimes\phi\to T(\psi)S(\phi),\quad
\psi\in H_u, \phi\in H_v$$
intertwines the tensor product $u\otimes v$ with $\eta$. By a)
$T(\psi)S(\phi)$ lies in ${\cal C}_{sp}$, so ${\cal C}_{sp}$ is an
algebra. Furthermore the proof of part b) of the previous proposition
shows that ${\cal C}_{sp}$ is a $^*$--subalgebra.
\bigskip

Consider an action
 $\eta:{\cal B}\to{\cal
B}\otimes{\cal A}$ 
of 
$G=({\cal A},\Delta)$
on a unital $C^*$--algebra ${\cal B}$ (with $\otimes$ the minimal tensor product):
a unital $^*$--homomorphism such that $\eta\otimes\iota\circ\eta=\iota\otimes\Delta\circ\eta$.

In the $C^*$--algebraic case, one can similarly define,  the
spectrum of the action, $sp(\eta)$, and the spectral $^*$--subalgebra
${\cal B}_{sp}$, a unital $^*$--subalgebra invariant under
the action of $G_\infty$. 
In the case where ${\cal B}={\cal A}$ and $\eta=\Delta$, 
${\cal A}_{sp}={\cal A}_\infty$.

 There is an alternative way of introducing the notion of spectrum. 
   Let $u$ and $v$ be corepresentations of a Hopf $C^*$--algebra on Hilbert
spaces $H$ and $K$ respectively. Then we may define a coaction $\eta$ on
the space
$(H,K)$ of linear mappings from $H$ to $K$ by setting 
$$\eta(b):=vb\otimes 1u^*,\,\, b\in(H,K).$$ 

   The above formulae may be used to define a coaction on $(H,K)\otimes \cC$ 
for a $C^*$--algebra $\cC$ with the variables in $\cC$ being spectators. 
If we are now given in addition a $C^*$--algebra $\cC$ carrying a coaction 
$\eta$ of $\cA$, then we can let it act on $(H,K)\otimes\cC$ with the 
variables in $(H,K)$ being spectators. Combining it with the coaction of 
$\cA$ on $(H,K)$ defined by $u$ and $v$ as above we get a coaction $\beta$ 
on $(H,K)\otimes\cC$. Since $\cA$ is not, in general, commutative, this 
coaction must be spelled out in detail. We have 
$$\beta(a):=\hat v\eta(a)\hat u^*,\quad a\in(H,K)\otimes\cC,$$ 
where $\hat v$ denotes $v$ with $1_\cC$ inserted between the tensor product 
factors in $v$. Choosing orthonormal bases $\psi_p,\chi_r$ in $H$ and $K$ and 
expressing $u$, $v$ and $a$ in terms of these bases, the coaction may be 
written 
$$\beta(a)=\sum_{p,q,r,s}\chi_r\psi_p^*\otimes(1_\cC\otimes v_{rs}
\eta(a_{sq})1_\cC\otimes u^*_{pq}).$$ 
Using this or the previous formula, it is easy to check that $\beta$ is 
a coaction.\smallskip 
 
  The set of fixed points will be denoted $(u\otimes\eta,v\otimes\eta)$. 
Restricting ourselves to finite-dimensional unitary corepresentations and 
allowing $\cC$ to be just a $^*$--algebra we define the spectral category 
Sp$(\cC,\eta)$ by letting $(u\otimes\eta,v\otimes\eta)$ be the set 
of arrows from the object $u\otimes\eta$ to the object $v\otimes\eta$
with the obvious law of composition. After choosing orthonormal bases, 
and setting
$$a=\sum_{s,t}\psi_t\varphi_s^*\otimes a_{ts},$$
$a\in(u\otimes\eta,v\otimes\eta)$ reads in coordinates 
$$\eta(a_{nm})=\sum_{t,s}I_\cC\otimes v_{nt}^*a_{ts}\otimes I_\cA
I_\cC\otimes u_{sm}.$$

  There is also a tensor product of a restricted nature in Sp$(\cC,\eta)$. 
If $a\in(u\otimes\eta,m\iota\otimes\eta)$, 
$a=\sum_{s,t}\varphi_s\psi_t^*\otimes a_{st}$, and 
$b\in(v\otimes\eta,n\iota\otimes\eta)$, 
$b=\sum_{p,q}\lambda_p\mu_p^*\otimes b_{pq}$,
where $m$ and $n$ are positive integers,
 then 
$$a\top{}b:=\sum_{p,q,s,t}\varphi_s\lambda_p\psi_t^*\mu_q^*\otimes
a_{st}b_{pq}$$ 
is in $(u\top{}v\otimes\eta,mn\iota\otimes\eta)$. In fact,
$\eta(a_{st})=\sum_fa_{sf}\otimes u_{ft}$ and 
$\eta(b_{pq})=\sum_gb_{pg}\otimes v_{gq}$ so that 
$$\eta(a_{st}b_{pq})=\sum_{f,g}a_{sf}b_{pg}\otimes u_{ft}v_{gq},$$ 
as required. Notice that if $a$ and $b$ are unitaries then so is
$a\top{}b$.\smallskip

 One can check that if $S\in(u',u)$, $T\in(v,v')$ and
$X\in(u\otimes\eta,v\otimes\eta)$ then
$T\otimes 1XS\otimes 1\in(u'\otimes\eta,v'\otimes\eta)$.

 An element of $\cC_u:=(u\otimes\eta,\iota\otimes\eta)$, where $\iota$
denotes the trivial corepresentation, will be called a $u$--multiplet
and may be thought of as a multiplet transforming according to the
corepresentation $u$. More precisely, we have
$$\eta(c_i)=\sum_jc_j\otimes u_{ji}.$$
Thus these multiplets coincide with those introduced above when
defining the spectrum. In fact, if we define the coordinates of
$T:H_u\to\cC$ with respect to an orthonormal basis by $T\phi_i:=T_i$
and the coordinates of $S\in(u\otimes\eta,\iota\otimes\eta)$ by
$S:=\sum_i\phi^*\otimes S_i$ then setting $S_i:=T_i$ defines a
canonical isomorphism from the set of intertwining maps to
$(u\otimes\eta,\iota\otimes\eta)$.

 If $\eta:\cC\to\cC\otimes\cA$ and $\beta:\cD\to\cD\otimes\cA$ are two 
coactions of $\cA$ and $k:\cC\to\cD$ is a morphism commuting with the coactions, 
i.e.\ $k\otimes I_\cA\circ\eta=\beta\circ k$ then there is an induced 
functor $k_*:$Sp$\eta\to$Sp$\beta$, $k_*(u\otimes\eta):=(u\otimes\beta)$ 
and $k_*(C)=1_{(H_u,H_v)}\otimes k(C)$ for $C\in(u\otimes\eta,v\otimes\eta)$. 
Note that $k_*$ being a $^*$--functor of $C^*$--categories will map 
unitaries to unitaries and isometries to isometries.\bigskip

\noindent{\it 2.3 Multiplicities of spectral representations}
\bigskip

We call an action $\eta: {\cal B}\to{\cal B}\otimes{\cal A}$ of
$G=({\cal A},\Delta)$ on a unital $C^*$--algebra ${\cal B}$ {\it
nondegenerate} if $\eta({\cal B})I\otimes{\cal A}$ is dense in ${\cal
B}\otimes{\cal A}$.
In \cite{Podles} the following result is proven.

\medskip

\noindent{\bf 2.3 Theorem} {\sl \cite{Podles} 
Let $\eta:{\cal B}\to{\cal
B}\otimes{\cal A}$ be a nondegenerate action of a compact quantum group
$G$ on  ${\cal B}$. Then for any irreducible $u\in
sp(\eta)$ there is a subspace $W_u\subset{\cal B}$ containing any spectral
multiplet tranforming like  $u$, such that
\begin{description} 
\item{\rm a)} $W_u$ splits 
$W_u=\oplus_{i\in I_u} W_u^i$
into an algebraic direct sum of 
subspaces $W_u^i$, each of them corresponding to $u$.
\item{\rm b)} If $\hat{\eta}$ is a complete set of inequivalent
irreducibles in $sp(\eta)$ then the linear span  of
$\{W_u,
u\in\hat{\eta}\}$, which coincides with the spectral
$^*$--subalgebra ${\cal B}_{sp}$, is dense in ${\cal
B}$.
\end{description}
}\medskip

The cardinality of the set $I_u$ is called the {\it multiplicity} of
the irreducible $u$
in $\eta$, and it will be denoted $\text{mult}(u)$. 
\bigskip

\noindent{\it 2.4 Quantum subgroups and quotient spaces}
\bigskip

A compact quantum subgroup $K=({\cal A}',\Delta')$ of $G$, as
introduced in \cite{Podles}, is a  compact quantum group for which
there
exists a surjective $^*$--homomorphism $\pi:{\cal A}\to{\cal A}'$ such
that
$$\pi\otimes\pi\circ\Delta=\Delta'\circ\pi.$$
A {\it closed bi--ideal} ${\cal I}$ of ${\cal A}$ is a norm closed
two--sided ideal of
${\cal A}$ such that
$$\Delta({\cal I})\subset {\cal
A}\otimes{\cal I}+{\cal I}\otimes{\cal A}.$$ 
There is a surjective correspondence from 
  quantum subgroups of $G$ 
and closed bi--ideals of ${\cal A}$, which
associates to  a subgroup defined by the
surjection $\pi$, the kernel of $\pi$.
Subgroups defined by surjections with the same kernel are 
isomorphic  as Hopf $C^*$--algebras \cite{P2}.

The subgroup $K$ acts (on the left) on the $C^*$--algebra ${\cal A}$ via
$$\delta:=\pi\otimes\iota\circ\Delta: {\cal A}\to{\cal A}'\otimes{\cal
A}.$$
The fixed point algebra 
$${\cal A}^\delta:=\{T\in{\cal A}: \delta(T)=I\otimes T\}$$
is defined to be the quantum quotient space of right cosets.
This algebra
has
a natural right action of $G$:
$$\eta_K:=\Delta\upharpoonright_{{\cal A}^\delta}: {\cal A}^\delta\to{\cal
A}^\delta\otimes{\cal A},$$
known to be nondegenerate and ergodic \cite{Wang}, see also \cite{P2}.
We set:
$$K\backslash G:= ({\cal A}^\delta,\eta_K).$$ 
One can similarly consider the action of $K$ on the right on ${\cal A}$:
$$\rho:=\iota\otimes\pi\circ\Delta: {\cal A}\to{\cal A}\otimes{\cal
A}',$$
with fixed point algebra 
$${\cal A}^\rho:=\{T\in{\cal A}: \delta(T)=I\otimes T\},$$
called  the quantum quotient space of left cosets. This algebra carries
a left action of $G$:
$$\eta^K:=\Delta\upharpoonright_{{\cal A}^\rho}: {\cal A}^\rho\to {\cal
A}\otimes{\cal A}^\rho,$$
and we set
$$G/K:=({\cal A}^\rho, \eta^K).$$

As in the group case,  unitary  representations of $G$ can be
restricted to unitary representations of 
$K$ on the same Hilbert space:
  $$u\upharpoonright_K:=\iota\otimes\pi\circ u: H_u\to H_u\otimes{\cal
A}'.$$
\medskip

\noindent{\bf 2.4 Proposition} {\sl \cite{P2} 
If $T\in(u,v)$  then $T\in(u\upharpoonright_K, v\upharpoonright_K)$
as well. So
the map 
$\text{Rep}(G)\to\text{Rep}(K)$,
taking 
$u\to u\upharpoonright_K$, and acting trivially on the arrows
defines a faithful tensor $^*$--functor. The smallest full tensor
$^*$--subcategory of $\text{Rep}(K)$ with subobjects and direct sums
containing the $u\upharpoonright_K$, for $u\in\text{Rep}(G)$, is
$\text{Rep}(K)$. 
}\medskip

We shall consider 
the subspace $K_u$ of $H_u$  of $K$--invariant vectors for the restricted 
representation $u\upharpoonright_K$: this is the set
of all $k\in H_u$ for which
$u\upharpoonright_K(k)=k\otimes I$.
The dimension of $K_u$ is the multiplicity of the trivial
representation $\iota_K$ of $K$ in $u\upharpoonright_K$. 
\medskip

\noindent{\bf 2.5 Proposition} {\sl If $K$ is a compact quantum subgroup
of
$G$ defined by the surjection $\pi$ then for any representation $u$ of $G$
and vectors $\psi\in H_u$, $k\in K_u$, the elements 
$u_{\psi,k}-(\psi, k)I$ and $u_{k,\psi}-(k,\psi)I$ belong to the closed
bi--ideal
$\text{ker}\pi$.}\medskip

\noindent{\bf Proof}
Notice then that for any $\psi\in H_u$ and $k\in K_u$
$$(u\upharpoonright_K)_{\psi,k}=\pi(u_{\psi,k})=(\psi,k)I,$$
so $u_{\psi,k}-(\psi,k)I\in\text{ker}\pi$.
On the other hand, picking an orthonormal basis $(\psi_i)$ of $H_u$,
the condition for $k$ to be a fixed vector for $u\upharpoonright_K$ can be
written 
$$\sum_i\pi(u_{ji})(\psi_i, k)=(\psi_j, k)I,$$
where $u_{ji}:=u_{\psi_j,\psi_i}$.
Since the matrix $(u_{ji})$ is unitary, this condition can be rewritten in
terms of $u^*$ and reads 
$$\sum_j\pi(u^*_{rj})(\psi_j,k)=(\psi_r,k)I,$$
or, taking the adjoint,
$$\sum_j\pi(u_{rj})(k,\psi_j)=(k,\psi_r)I.$$
Thus we also have, for
$k\in K_u$, $\psi\in H_u$, $u\in\text{Rep}(G)$,
$$(u\upharpoonright_K)_{k,\psi}=\pi(u_{k,\psi})=(k,\psi)I.$$
\medskip

As in the group case, we can construct elements of the coset spaces 
using invariant vectors for the subgroup: 
for 
any representation $u$ of $G$, 
if we pick vectors $k\in K_u$, $\psi\in H_u$,
the coefficient
$u_{k,\psi}$ lies in ${\cal A}^\delta$ and $u_{\psi,k}$ lies in ${\cal
A}^\rho$. We only show the former: if $(\psi_j)$ is an orthonormal basis
of $H_u$:
$$\delta(u_{k,\psi_j})=\pi\otimes\iota\circ\Delta(u_{k,\psi_j})=$$
$$\sum_i\pi(u_{k,\psi_i})\otimes
u_{ij}=\sum_{i}I\otimes (k,\psi_i)u_{ij}=I\otimes u_{k,\psi_j}.$$
Fix a complete set $\hat{G}$ of unitary irreducible representations of
$G$, and set
$$\hat{G}_K:=\{u\in\hat{G}: K_u\neq0\}.$$
\medskip

\noindent{\bf 2.6 Proposition} {\sl 
The 
linear
space generated  
by the matrix coefficients
$\{u_{k,\psi}\},$
(resp. $\{u_{\psi,k}\},$)
as $u\in\hat{G}_K$, $k\in K_u$, $\psi\in H_u$, vary,
coincides with  ${\cal A}^\delta_{sp}$ (resp. ${\cal A}^\rho_{sp})$).
}\medskip

\noindent{\bf Proof} Let $V$ denote the  linear space defined in the statement. 
We show, for completeness,  that $V$ is dense.
Consider the conditional expectation $E:{\cal A}\to{\cal A}^{\delta}$
onto the fixed point algebra, obtained by averaging the
$K$--action: $E(T):=h'\otimes\iota_{\cal A}\circ\delta(T)$, with $h'$ the
Haar measure of $K$. Since ${\cal A}$ is generated, as a Banach space,
by the coefficents of its irreducible unitary representations
$u_{\psi',\psi}$, $\psi,\psi'\in H_u$,  the set $\{E(u_{\psi',\psi}):
u\in\hat{G}, \psi,\psi'\in H_u\}$,
is total in ${\cal A}^{\delta}$. 
A computation shows that $E(u_{\psi',\psi})=\sum_k
h'\circ\pi(u_{\psi',\psi_k})u_{\psi_k,\psi}$, with $\{\psi_k\}$ any
orthonormal basis of $H_u$.
By \cite{Wcmp},
$h'\pi(u_{\psi',\phi})=h'((u\upharpoonright_K)_{\psi',\phi})=0$
for all $\psi',\phi\in H_u$, unless
$u\upharpoonright_K$
contains the identity representation. Assume then that this is the case.
Replace, if necessary, that orthonormal basis with another one such that
 $u\upharpoonright_K$, when represented as a   matrix with
entries in
${\cal
A'}$, becomes diagonal of the form 
$\text{diag}(I,\dots,I, v)$, with $v$ a unitary representation of
$K$ which does not contain the trivial representation, so that
$h'\pi(u_{\psi_i,\psi_k})=0$
unless $i\leq \text{dim}(K_u)$ and $k=i$. In that case,
$h'\circ\pi(u_{\psi_i,\psi_i})=1$. Therefore 
$E(u_{\psi',\psi})=u_{\psi',\psi}$ if $\psi'\in K_u$, and
$E(u_{\psi',\psi})=0$ if $\psi'\in K_u^\perp$.
We have thus shown that $V$ is norm dense in ${\cal A}^\delta$.
One can easily check that for any $u\in\hat{G}_K$, $\phi\in K_u$,  the map
$\psi\in H_u\to u_{\phi,\psi}\in{\cal A}^\delta$ intertwines $u$ with $\eta_K$, so
$V$  is contained in ${\cal A}^\delta_{sp}$. Podles shows in \cite{Podles}
that the subspace
$W_u$ described in Theorem 2.3 coincides with the linear span of
$\{u_{\phi,\psi}\}$, $\phi\in
K_u$, $\psi\in H_u$, so by Theorem 2.3, b) $V={\cal A}^\delta_{sp}$.
\medskip

As a consequence, the above   subspace 
is a unital  $^*$--algebra endowed with the restricted action
 of the Hopf $^*$--algebra $G_\infty:=({\cal
A}_\infty,\Delta)$, still
denoted by $\eta_K:
 {\cal A}^{\delta}_{sp}
\to
 {\cal A}^{\delta}_{sp}
\odot{\cal A}_\infty$.
\bigskip

\noindent{\it 2.5 Stabilizer and kernel}
\bigskip

  In this subsection we examine situations that should give rise to a 
compact quantum subgroup.

  Let $\alpha$ be an action of a compact quantum group $G:=(\cA,\Delta)$ on 
a $C^*$--algebra $\cB$ and $\varphi$ a state of $\cB$. We look for an 
appropriate notion of the stabilizer of $\varphi$ under the action. Consider 
$$S_\varphi:=\{\pi\in\text{Rep}\cA:\varphi\otimes\pi(\alpha(B))=\varphi(B)
\otimes 1_\pi,\,\, B\in\cB.\}$$ 
Note that if $G$ is a compact group, the irreducible representations of $\cA$ 
are equivalent to characters and labelled by the elements of $G$. Thus the 
irreducible part of $S_\varphi$ corresponds exactly to the stabilizer of $\varphi$ 
under the action. We now investigate the stability properties of $S_\varphi$. 
If $W\in(\pi',\pi)$ is an isometry and $\pi\in S_\varphi$ then 
$$\varphi\otimes\pi'(\alpha(B))=W^*\varphi\otimes\pi(\alpha(B))W=
W^*\varphi(B)\otimes 1_\pi W=\varphi(B)\otimes 1_{\pi'},$$ 
so that $\pi'\in S_\varphi$. Similarly, if $W_i\in(\pi_i,\pi)$ are isometries 
with $\sum_iW_iW_i^*=1_\pi$ and $\pi_i\in S_\varphi$ for all $i$, then
$$\varphi\otimes\pi(\alpha(B))=\sum_iW_i\varphi\otimes\pi_i(\alpha(B))W_i^*
=\sum_i\varphi(B)\otimes W_iW_i^*=\varphi(B)1_\pi,$$ 
so that $\pi\in S_\varphi$ and $S_\varphi$ is closed under subobjects and 
direct sums. If $\pi,\pi'\in S_\varphi$ then 
$$\varphi\otimes(\pi*\pi')(\alpha(B))=\varphi\otimes\pi\otimes\pi'(\iota\otimes
\Delta\circ\alpha(B))=\varphi\otimes\pi\otimes\pi'(\alpha\otimes\iota\circ
\alpha)(B).$$ 
For any $B\in\cB$ and $A\in\cA$, we have 
$$\varphi\otimes\pi\otimes\pi'(\alpha\otimes\iota)(B\otimes A)=
\varphi(B)\otimes 1_\pi\otimes\pi'(A)=\varphi\otimes\iota\otimes\pi'
(B\otimes 1_\pi\otimes A).$$ 
Thus for any $C\in\cB\otimes\cA$, we shall have 
$\varphi\otimes\pi\otimes\pi'\circ\alpha\otimes\iota(C)=
\varphi\otimes\iota\otimes\pi'(\hat C)$, where $\hat C$ denotes $C$ with a 
$1_\pi$ inserted in the middle. Taking $C=\alpha(B)$, we get 
$$\varphi\otimes(\pi*\pi')(\alpha(B))=\varphi\otimes\iota\otimes\pi'(
\widehat{\alpha(B)})=\varphi(B)\otimes 1_{\pi*\pi'},$$ 
so that $\pi*\pi'\in S_\varphi$.

   Before pursuing this line of reasoning, we ask when a surjective 
morphism $\pi:\cA\to \pi(\cA)$ defines a quantum subgroup of $\cA$. 
We must be able to equip $\pi(\cA)$ with a comultiplication $\Delta'$ 
such that $\pi\otimes\pi\Delta=\Delta'\pi$. This equation will have a 
solution if and only if $\pi(A)=0$ implies $\pi\otimes\pi\Delta(A)=0$. 
Any solution is unique since $\pi$ is surjective and for the same reason 
$\Delta'$ will be a comultiplication because $\Delta$ is a comultiplication. 
Finally, $\pi(\cA)\otimes I\Delta'\pi(\cA)$ and $I\otimes\pi(\cA)\Delta'\pi(\cA)$ 
are dense in $\pi(\cA)\otimes\pi(\cA)$ since $\cA\otimes I\Delta(\cA)$ and 
$I\otimes\cA\Delta(\cA)$ are dense in $\cA\otimes\cA$. Thus we have a simple 
necessary and sufficient condition for a surjective morphism 
$\pi:\cA\to \pi\cA$ to define a quantum subgroup which may be written 
$$\text{ker}\pi*\pi\supset\text{ker}\pi.$$ 

  We now pick one representative $\pi_i$ from each equivalence class of 
finite-dimensional irreducibles in $S_\varphi$ and form the direct sum,
say, 
$\pi=\oplus_{i\in I}\pi_i$. To show that $\pi$ defines a quantum subgroup 
we must show that $\pi_i(A)=0$ for $i\in I$ implies $\pi_j*\pi_k(A)=0$. 
But this is the case since, by the results above, $\pi_j*\pi_k$ is a finite 
direct sum of the $\pi_i$, $i\in I$.

As an example of an infinite dimensional representation
belonging to a stabilizer,
consider the ergodic action $\eta_K$ of $G=({\cal A},\Delta)$ on the
quantum
quotient space ${\cal A}^\delta$
 by a quantum subgroup $K$, and
let $e$ be an everywhere
defined counit of 
${\cal A}$, which we restrict to a state $\varphi$ of ${\cal A}^\delta$.
Let
$\pi:{\cal
A}\to{\cal A}'$ be the surjection  defining the subgroup, regarded as a 
representation of ${\cal A}$.
We compute 
$\varphi\otimes\pi(\eta_K(B))$ 
on the elements $B=u_{k,\psi}$,
with $k\in K_u$, $\psi\in H_u$:
$$\varphi\otimes\pi(\eta_K(u_{k,\psi}))=\varphi\otimes\pi(\sum_j
u_{k,\psi_j}\otimes
u_{\psi_j,\psi})=$$
$$\sum_j
(k,\psi_j)\otimes\pi(u_{\psi_j,\psi})=I\otimes\pi(u_{k,\psi})=I\otimes
(k,\psi)1_\pi$$
by prop. 2.5, and this expression equals in turn 
$\varphi(u_{k,\psi})\otimes 1_\pi$. Since the elements $u_{k,\psi}$ span a
dense subspace of ${\cal A}^\delta$, we can conclude that $\pi\in
S_\varphi$.

  Similarly, we can look for the kernel of the action of a compact quantum group. 
Letting $\alpha:\cB\to\cB\otimes\cA$ be the action of a compact quantum group 
on a $C^*$--algebra $\cB$, let 
$$K_\alpha:=\{\pi\in\text{Rep}\cA:1_\cB\otimes\pi\alpha(B)=B\otimes 1_\pi,\,\,
B\in\cB\}.$$ 

  Arguing as above, we may verify that $K_\pi$ is closed under subobjects, 
direct sums but there seems to be no reason for it to be closed under tensor 
products.

  Having failed to define the kernel of an action $\alpha$, we can at least 
define an action to be faithful when the irreducible part of $K_\alpha$ reduces 
to the trivial representation $\iota$.

  We may similarly talk of the kernel of a representation $u:H\to H\otimes\cA$ 
of a quantum group by defining 
$$M_u:=\{\pi\in\text{Rep}\cA:1_H\otimes\pi u(\psi)=\psi\otimes 1_\pi,\,\,
\psi\in H\}.$$ 
Arguing as above, we may verify that $M_u$ is closed under subobjects, 
direct sums and tensor products and proceed to define the kernel of $u$ 
as the quantum subgroup associated with $M_u$. $u$ is faithful when 
the irreducible part of $M_u$ reduces to the trivial representation
$\iota$.
\medskip

\end{section}

\begin{section} {Quasitensor functors and a  finiteness
theorem}

\noindent{\bf 3.1 Definition}
Let   ${\cal T}$ and ${\cal R}$ be  strict tensor $C^*$--categories
\cite{DRinventiones}. 
We shall always
assume that the tensor units are irreducible:
$(\iota,\iota)={\mathbb C}$. A (covariant)  $^*$--functor ${\cal F}: {\cal
T}\to {\cal R}$
will be called
{\it quasitensor} if 
$${\cal F}(\iota)=\iota,\eqno(3.1)$$
if for objects $\rho$, $\sigma\in{\cal T}$ there is an
isometry 
$$S_{\rho,\sigma}: {\cal F}(\rho)\otimes{\cal F}(\sigma)\to{\cal
F}(\rho\otimes\sigma)\eqno(3.2)$$
such that
$$S_{\rho,\iota}=S_{\iota,\rho}=1_{{\cal F}(\rho)},\eqno(3.3)$$
$$S_{\rho\otimes\sigma,\tau}\circ S_{\rho,\sigma}\otimes 1_{{\cal
F}(\tau)}=S_{\rho,\sigma\otimes\tau}\circ 1_{{\cal F}(\rho)}\otimes
S_{\sigma,\tau}=:S_{\rho,\sigma,\tau},\eqno(3.4)$$
$$E_{\rho\otimes\sigma,\tau}\circ E_{\rho,\sigma\otimes\tau}\leq
E_{\rho,\sigma,\tau},
\eqno(3.5)$$
with $E_{\rho,\sigma}\in({\cal F}(\rho\otimes\sigma), {\cal
F}(\rho\otimes\sigma))$ the range projection  of
$S_{\rho,\sigma}$ and $E_{\rho,\sigma,\tau}\in({\cal
F}(\rho\otimes\sigma\otimes\tau),{\cal F}(\rho\otimes\sigma\otimes\tau))$
the range projection of $S_{\rho,\sigma,\tau}$, and if 
$${\cal F}({S\otimes T})\circ
S_{\rho,\sigma}=S_{\rho',\sigma'}\circ{\cal
F}(S)\otimes{\cal F}(T),\eqno(3.6)$$ 
 for any other pair of objects $\rho'$, $\sigma'$ and arrows
$S\in(\rho, \rho')$,
$T\in(\sigma,\sigma')$.

In this paper we shall only deal with the case where 
   ${\cal R}$ is 
the tensor $C^*$--category
${\cal H}$
with objects Hilbert spaces and arrows 
from a Hilbert space $H$ to a Hilbert space $H'$, 
the set 
$(H, H')$ 
of all
bounded linear mappings from $H$ to $H'$.
We shall assume that ${\cal H}$ is strictly tensor, namely that the tensor
product between  the objects has been realized in a 
 strictly associative way.

In order to economize on brackets, we evaluate tensor products of arrows
before composition. Moreover in the sequel, for the sake of simplicity, we
shall occasionally identify, with a
less precise
but lighter  notation,
${\cal F}(\rho)\otimes{\cal F}(\sigma)$ with a subspace of
${\cal F}(\rho\otimes\sigma)$, the  
image of
$S_{\rho,\sigma}$. Equations $(3.2)$--$(3.6)$ shall then be  written:
$${\cal F}(\rho)\otimes {\cal F}(\sigma)\subset {\cal
F}(\rho\otimes\sigma),\eqno(3.7)$$
$$E_{\rho,\iota}=E_{\iota,\rho}=1_{{\cal F}(\rho)},\eqno(3.8)$$
$$E_{\rho,\sigma}\otimes 1_{{\cal F}(\tau)}\circ E_{\rho\otimes\sigma,\tau}=
 1_{{\cal F}(\rho)}\otimes E_{\sigma,\tau}\circ E_{\rho,\sigma\otimes\tau}=:E_{\rho,\sigma,\tau},
\eqno(3.9)$$
$$E_{\rho\otimes\sigma,\tau}({\cal F}(\rho)\otimes{\cal F}(\sigma\otimes\tau))
\subset{\cal F}(\rho)\otimes{\cal F}(\sigma)\otimes{\cal
F}(\tau),\eqno(3.10)$$
$${\cal F}({S\otimes T})\upharpoonright_{{\cal F}(\rho)\otimes
{\cal F}(\sigma)}={\cal F}(S)\otimes{\cal F}(T).\eqno(3.11)$$ 
Equation (3.10)
 combined with (3.9) requires that the projection  onto
${\cal F}(\rho\otimes\sigma)\otimes{\cal F}(\tau)$ actually takes the subspace
${\cal F}(\rho)\otimes{\cal F}(\sigma\otimes\tau)$ onto 
${\cal F}(\rho)\otimes{\cal F}(\sigma)\otimes{\cal F}(\tau)$.
Therefore we necessarily have  
$$E_{\rho,\sigma,\tau}=E_{\rho\otimes\sigma,\tau}\circ E_{\rho,
\sigma\otimes\tau}=E_{\rho, \sigma\otimes\tau}\circ
E_{\rho\otimes\sigma,\tau}.\eqno(3.12)$$ \medskip

Notice that
any tensor $^*$--functor from ${\cal T}$ to ${\cal H}$ is quasitensor.
\medskip

\noindent{\bf 3.2 Example}
Assume that  ${\cal T}=\text{Rep(G)}$, the representation category of a
compact quantum group $G$. Then 
the embedding
functor
$H:\text{Rep}(G)\to{\cal H}$
associating to each representation $u$ its Hilbert space $H_u$ and acting
trivially on the arrows, is
tensor, and therefore quasitensor.
\medskip

Quasitensor $^*$--functors arise naturally in abstract tensor
$C^*$--categories.
\medskip

\noindent{\bf 3.3 Proposition} {\sl Let ${\cal T}$ be a tensor
$C^*$--category with $(\iota,\iota)=\mathbb C1_\iota$.
For any object $\rho$ of ${\cal T}$, consider the
Hilbert space
$\hat{\rho}:=(\iota,\rho)$, with inner product
$$(\phi,\phi')1_\iota:=\phi^*\circ\phi',\quad \phi,\phi'\in
(\iota,\rho).$$
For $T\in(\rho,\sigma)$ define a bounded linear map
$\hat{T}:\hat{\rho}\to\hat{\sigma}$
by $\hat{T}(\phi)=T\circ\phi$. 
This is a quasitensor $^*$--functor.}\medskip

\noindent{\bf Proof} It is easy to check that, for $T\in(\rho,\sigma)$, $S
\in(\sigma,\tau)$, $\widehat{S\circ
T}=\hat{S}\circ\hat{T}$ and that $\widehat{T^*}=\hat{T}^*$. This shows
that we
have a $^*$--functor such that $\hat{\iota}=(\iota,\iota)={\mathbb C}$.
and $(3.1)$ is satisfied.
For $\phi\in \hat{\rho}=(\iota,\rho)$, $\psi\in
\hat{\sigma}=(\iota,\sigma)$, the
map $\phi,\psi\to
\phi\otimes\psi=\phi\otimes 1_\sigma\circ\psi\in
\widehat{\rho\otimes\sigma}=(\iota,\rho\otimes\sigma)$
defines
an isometric map
from $\hat{\rho}\otimes \hat{\sigma}$ to $\widehat{{\rho\otimes\sigma}}$. 
The
copy of
$\hat{\rho}\otimes \hat{\sigma}\otimes \hat{\tau}$ sitting inside
$\widehat{\rho\otimes\sigma}\otimes \hat{\tau}$ is the subspace generated
by 
elements 
$\phi\otimes1_{\sigma\otimes\tau}\circ\psi\otimes 1_{\tau}\circ\eta$, for $\phi\in \hat{\rho}$, $\psi\in
\hat{\sigma}$, $\eta\in\hat{\tau}$, and it coincides with the copy of the
same
Hilbert space sitting inside $\hat{\rho}\otimes \widehat{\sigma\otimes
\tau}$,
hence
$(3.9)$ holds. 
We now check $(3.10)$. First notice that if $\phi_i$ is an orthonormal
basis of $\hat{\sigma}$ then every finite sum  $\sum_{\text{finite}} 
1_\rho\otimes
(\phi_i\circ\phi_i^*)$ is an element of $(\rho\otimes\sigma,
\rho\otimes\sigma)$, defining, by composition, a projection map
from $\widehat{{\rho\otimes \sigma}}$ onto a subspace of $\hat{\rho}\otimes
\hat{\sigma}$. The strong limit of this net converges, in the strong
topology
defined by $\widehat{{\rho\otimes\sigma}}$, to the projection map
$E_{\rho,\sigma}$. This shows that the
projection map of $E_{\rho\otimes\sigma,\tau}$ is the strong limit of 
$\sum_i \widehat{1_{\rho\otimes\sigma}\otimes\psi_i\circ\psi_i^*}$, with
$\psi_i\in
\hat{\tau}$ an orthonormal basis. Thus $E_{\rho\otimes\sigma,\tau}$ 
 takes $\hat{\rho}\otimes \widehat{{\sigma\otimes\tau}}$ into
$\hat{\rho}\otimes
\hat{\sigma}\otimes \hat{\tau}$,
and the proof of  
$(3.10)$ is complete.
\medskip

The following proposition is of help in constructing 
more examples.
\medskip

\noindent{\bf 3.4 Proposition} {\sl Let ${\cal S}$, ${\cal T}$ be tensor
$C^*$--categories. If ${\cal G}: {\cal S}\to {\cal T}$ is a tensor
$^*$--functor and ${\cal F}:{\cal T}\to{\cal H}$ is quasitensor then
${\cal F}\circ {\cal G}$ is quasitensor.}\medskip

\noindent{\bf Proof} Property $(3.1)$ is obvious.  Set ${\cal E}:={\cal
F}\circ{\cal G}$. Since 
${\cal G}$ is a tensor and ${\cal F}$ is quasitensor,  $${\cal
E}(\rho)\otimes{\cal
E}(\sigma)\subset {\cal F}({\cal G}(\rho)\otimes{\cal G}(\sigma))=
{\cal F}({\cal G}(\rho\otimes\sigma))={\cal E}(\rho\otimes\sigma).$$
Hence $(3.7)$ and $(3.11)$ hold. Set $E^{\cal E}_{\rho,\sigma}:=E^{\cal
F}_{{\cal G}(\rho), {\cal G}(\sigma)}$. This is the orthogonal projection
from ${\cal E}(\rho\otimes\sigma)$ to ${\cal E}(\rho)\otimes{\cal
E}(\sigma)$, and it is easy to check that it satisfies properties
$(3.8)$, $(3.9)$ and $(3.10)$. \medskip

\noindent{\it Remark} One can similarly show that if ${\cal F}$ is tensor
and ${\cal G}$ is quasitensor then ${\cal F}\circ {\cal G}$ is
quasitensor.\medskip

\noindent{\bf 3.5  Example} Let $G$ be a compact quantum group, 
${\cal T}$ be a tensor $C^*$--category and 
$\rho: \text{Rep}(G)\to{\cal T}$ is a tensor $^*$--functor.
If we compose $\rho$  
with
the quasitensor $^*$--functor ${\cal T}\to{\cal H}$  associated
with ${\cal T}$ as in  Prop. 3.3, we obtain 
an interesting `abstract' quasitensor
$^*$--functor
${\cal F}_\rho: \text{Rep}(G)\to{\cal H}$ with Hilbert spaces 
${\cal F}_\rho(u)=(\iota, \rho_u)$ and defined on arrows by 
$${\cal
F}_\rho(T)\phi:=\hat{\rho(T)}\phi=\rho(T)\circ\phi\in(\iota,\rho_v),\quad
\phi\in(\iota,\rho_u),
T\in(u,v).$$
\medskip

\noindent{\bf 3.6 Example} 
Apply the construction described in the previous example to the tensor
$C^*$--category ${\cal T}=\text{Rep}(K)$, where $K$ is a 
compact quantum subgroup of $G$, with $\rho$ given by the canonical 
tensor $^*$--functor 
$\text{Rep}(G)\to \text{Rep}(K)$
described in Prop. 2.4. We now obtain a `concrete' quasitensor
$^*$--functor, that we shall denote, with abuse of notation,  by
$K:\text{Rep}(G)\to{\cal H}$, where now $K_u=(\iota_K,
u\upharpoonright_K)$. We shall call $K$ the {\it invariant vectors
functor}. 
\medskip

In Sec.\  7  we shall exhibit examples of quasitensor
$^*$--functors associated to ergodic actions of compact quantum groups.
The rest of this section is devoted to showing the following theorem.
\medskip

\noindent{\bf 3.7 Theorem} {\sl
Let ${\cal T}$ be a strict tensor $C^*$--category and 
${\cal F}:{\cal T}\to{\cal H}$ a quasitensor $^*$--functor. If $\rho$ 
has a conjugate $\overline{\rho}$ in ${\cal T}$ then
 ${\cal F}(\rho)$ is finite dimensional and 
$$\text{dim}{\cal F}(\rho)=\text{dim}{\cal F}(\overline{\rho}).$$
Furthermore if ${\cal F}(\rho)\neq 0$ and if  
$R\in(\iota,\overline{\rho}\otimes\rho)$
and $\overline{R}\in(\iota,\rho\otimes\overline{\rho})$ is a solution of 
the 
conjugate equations for $\rho$ in ${\cal T}$ then
 $\hat{R}:=S_{\overline{\rho},\rho}^*\circ {\cal F}(R)\in{\cal
F}(\overline{\rho})
\otimes{\cal F}(\rho)$ and
$\hat{\overline{R}}:=S_{\rho,\overline{\rho}}^*\circ{\cal
F}(\overline{R})\in{\cal F}(\rho)
\otimes{\cal F}(\overline{\rho})$ is a solution of the conjugate equations
for ${\cal F}(\rho)$ in
${\cal H}$.}\medskip

\noindent{\bf Proof} 
For the sake of simplicity, we shall use the simplified notation
$(3.7)$--$(3.11)$. This amounts to replacing $S_{\rho,\sigma}$
by the
identity and $S_{\rho,\sigma}^*$ by $E_{\rho,\sigma}$.
 Apply the functor
${\cal F}$ to
the relation $$\overline{R}^*\otimes 1_\rho\circ 1_\rho\otimes R=1_\rho$$
and get 
$${\cal F}(\overline{R}\otimes 1_\rho)^*\circ {\cal F}(1_\rho\otimes R)=
1_{{\cal F}(\rho)}.$$
Using successively $(3.11)$, $(3.8)$, $(3.9)$, $(3.10)$, we get
$${\cal F}(1_\rho\otimes R)(\psi)=\psi\otimes{\cal F}(R),
\quad\psi\in {\cal F}(\rho),$$
$${\cal F}(\overline{R}\otimes 1_\rho)(\phi)=
{\cal F}(\overline{R})\otimes\phi,
\quad\phi\in{\cal F}({\rho}).$$
So
$$(\phi,\psi) =
(\phi, {\cal F}(\overline{R}\otimes 1_\rho)^*\circ 
{\cal F}(1_\rho\otimes R)\psi)
=$$
$$({\cal F}(\overline{R})\otimes\phi, 
\psi\otimes{\cal F}(R))=$$
$$({\cal F}(\overline{R})\otimes\phi, 
E_{\rho\otimes\overline{\rho}, \rho}(\psi\otimes{\cal F}(R))=$$
$$({\cal F}(\overline{R})\otimes\phi, E_{\rho,\overline{\rho},\rho}\circ
E_{\rho\otimes\overline{\rho}, \rho}
(\psi\otimes{\cal F}(R))=$$
$$(E_{\rho,\overline{\rho}}\circ{\cal F}(\overline{R})\otimes\phi, 
\psi\otimes E_{\overline{\rho},\rho}\circ{\cal F}(R)),$$ so
$$\hat{\overline{R}}^*\otimes 1_{{\cal F}(\rho)}\circ
1_{{\cal F}(\rho)}\otimes\hat{R}=1_{{\cal F}(\rho)}.$$
At this point we can start over with
$\overline{\rho}$ and obtain the relation
$$\hat{R}^*\otimes 1_{{\cal F}(\overline{\rho})}\circ 
1_{{\cal F}(\overline{\rho})}\otimes\hat{\overline{R}}=
1_{{\cal F}(\overline{\rho})},$$
which implies that ${\cal F}(\overline{\rho})$ is finite dimensional as
well with the same dimension as ${\cal F}(\rho)$.
\medskip

\noindent Recall that in a tensor $C^*$--category with conjugates, the infimum of all
the
$$d_{R,\overline{R}}(\rho):=\|R\|\|\overline{R}\|$$ is the {\it intrinsic dimension}
of $\rho$, denoted $d(\rho)$ \cite{LongoRoberts}. If $\rho$ is irreducible
(in the sense that
$(\rho,\rho)={\mathbb C}$)
the spaces $(\iota,\overline{\rho}\otimes\rho)$ and $(\iota,\rho\otimes\overline{\rho})$ are one
dimensional, so any solution of the conjugate equations is of the form $\lambda R$,
$\mu\overline{R}$ with $\overline{\mu}\lambda=1$. Therefore in this case
for any solution $(R,\overline{R})$ of the conjugate equations,
$$d(\rho)= \|R\|\|\overline{R}\|.$$
\medskip

\noindent{\bf 3.8 Corollary} {\sl 
If ${\cal F}:{\cal T}\to{\cal H}$ is a quasitensor $^*$--functor and if $\rho$
is an  object of ${\cal T}$ with 
 a conjugate defined by $R$ and $\overline{R}$ then
$$\text{dim}({\cal F}(\rho))\leq d_{\hat{R},\hat{\overline{R}}}({\cal
F}(\rho))\leq   
d_{R,\overline{R}}(\rho).$$
Furthermore 
$d_{\hat{R},\hat{\overline{R}}}({\cal F}(\rho))=d_{R,\overline{R}}(\rho)$
if and only if 
${\cal F}(R)\in
\text{Image}S_{\overline{\rho},\rho}$
and 
${\cal F}(\overline{R})\in\text{Image}S_{\rho,\overline{\rho}}$.
}
\medskip

\noindent{\bf Proof}
Note that
 $$\|{\cal F}(R)\|^2=\|{\cal F}(R)^*{\cal F}(R)\|=\|{\cal
F}(R^*R)\|=
\|R\|^2,$$
so $\|\hat{R}\|\leq\|R\|$ and, similarly,
$\|\hat{\overline{R}}\|\leq\|\overline{R}\|$. Thus the last inequality 
follows. This also shows that
$d_{\hat{R},\hat{\overline{R}}}({\cal F}(\rho))=d_{R,\overline{R}}(\rho)$
if and only if $\|\hat{R}\|=\|{\cal F}(R)\|$ and 
$\|\hat{\overline{R}}\|=\|{\cal F}(\overline{R})\|$, 
and the last statement follows.

Let $J:{\cal F}(\rho)\to{\cal F}(\overline{\rho})$ be the antilinear
invertible associated to $\hat{R}$,
$\hat{\overline{R}}$ by $J\psi=r_\psi^*\circ\hat{R}$. Then
$$d_{\hat{R},\hat{\overline{R}}}^2({\cal F}(\rho))=
\text{Trace}(JJ^*)\text{Trace}((JJ^*)^{-1}).$$
The first inequality is now a consequence of the elementary  fact
that
for any positive invertible matrix $Q\in M_n$, 
$$n^2\leq\text{Trace}(Q)\text{Trace}(Q^{-1}).$$\medskip

In Sect.\ 7 we shall relate this result to the work of
\cite{BRV}.\medskip

  Let $\cF:\cS\to\cT$ be a quasitensor functor and $R,\bar R$ be a solution 
of the conjugate equations for an object $\rho$ of $\cS$. Then, as we have 
seen, $\hat R:=E_{\bar\rho,\rho}\circ\cF(R)$ and 
$\hat{\bar R}:=E_{\rho.\bar\rho}\circ\cF({\bar R})$ is a solution of the
conjugate equations 
for $\cF(\rho)$. We show that this construction has certain functorality 
properties. Given $T\in(\rho,\rho')$ define $T^\bullet$ by 
$$T^\bullet\otimes 1_\rho\circ R_\rho:=1_{\bar\rho'}\otimes T^*\circ R_{\rho'}.$$ 
Then 
$$\cF(T^\bullet\otimes 1_\rho)\circ\cF(R_\rho)=\cF(1_{\bar\rho'}\otimes
T^*)\circ\cF(R_{\rho'}),$$ 
so 
$$E_{\bar\rho',\rho}\circ\cF(T^\bullet\otimes 1_\rho)\circ\cF(R_\rho)=
E_{\bar\rho',\rho}\circ\cF(1_{\bar\rho'}\otimes T^*)\circ\cF(R_{\rho'}),$$ 
and then
$$\cF(T^\bullet)\otimes 1_{\cF(\rho)}\circ \hat R_\rho=
1_{\cF(\bar\rho')}\otimes\cF(T^*)\circ\hat R_{\rho'}.$$ 
It follows that 
$$\cF(T^\bullet)\otimes 1_{\cF(\rho)}\circ\hat R_\rho=
\cF(T)^\bullet\otimes 1_{\cF(\rho)}\circ\hat R_\rho,$$
from which we get 
$$\cF(T)^\bullet=\cF(T^\bullet).$$ 
Now suppose that $R,\bar R$ is a standard solution of the conjugate equations. 
Then $R=\sum_i\bar W_i\otimes W_i\circ R_i$, where 
$R_i\in(\iota,\bar\rho_i\rho_i)$ and $\bar R_i\in(\iota,\rho_i\bar\rho_i)$
are 
normalized solutions of the conjugate equations for the irreducible $\rho_i$. 
Thus $\cF(R)=\sum_i\cF(\bar W_i\otimes W_i)\circ\cF(R_i)$ giving 
$\hat R=E_{\bar\rho,\rho}\circ\cF(R)=\sum_i\cF(\bar W_i)\otimes
\cF(W_i)\circ \hat R_i$. We similarly get 
$$\hat{\bar R}=\sum_i\cF(W_i)\otimes\cF(\bar W_i)\circ\hat{\bar R_i}.$$ 
Consequently, if the $\cF(\rho_i)$ are irreducible, $\hat R,\hat{\bar R}$ is 
a standard solution of the conjugate equations for $\cF(\rho)$.\smallskip 

 For future use we note that $S_{\rho, \sigma}$ continues to be a natural 
transformation when antilinear intertwiners are allowed, as follows 
from the next result.\medskip 

\noindent{\bf 3.9 Lemma} {\sl If $\cF:\cT\to\cH$ is a quasitensor functor 
then for the antilinear operators $J$ associated with solutions 
$\hat R,\hat{\bar R}$ of the conjugate equations we have 
$$J_{\rho\otimes\sigma}S_{\rho,\sigma}=S_{\bar\sigma,\bar\rho}
J_\sigma\otimes J_\rho\circ\theta_{\rho,\sigma},$$ 
provided $R_{\rho\otimes\sigma}:=1_{\bar\sigma}\otimes R_\rho\otimes 
1_\sigma\circ R_\rho$.}\smallskip 

\noindent
{\it Remark} For a discussion of antilinear arrows see \cite{PR1}, 
where it is pointed out that 
$J_\sigma\otimes J_\rho\circ\theta_{\rho,\sigma}$ is a natural tensor 
product to use for antilinear arrows.\smallskip 

\noindent
{\bf Proof.} Identifying $\cF(\rho)\otimes\cF(\sigma)$ with a subspace of 
$\cF(\rho\otimes\sigma)$, we have to show that 
$$J_{\rho\otimes\sigma}\psi\otimes\phi=J_\sigma\phi\otimes J_\rho\psi,
\quad \psi\in\cF(\rho),\,\,\phi\in\cF(\sigma).$$
If $\chi\in\cF(\rho\otimes\sigma)$ then 
$r_\chi^*\circ\hat R_{\rho\otimes\sigma}=
J_{\rho\otimes\sigma}\chi$. Hence
$$J_{\rho\otimes\sigma}\psi\otimes\phi=
r_{\psi\otimes\phi}^*\circ
\hat R_{\rho\otimes\sigma}=$$ 
$$r_{\psi\otimes\phi}^*\circ
E_{\bar\sigma\otimes\bar\rho,\rho\otimes\sigma}
{\cal F}(1_{\bar\sigma}\otimes R_\rho\otimes 1_\sigma\circ R_\sigma)=$$
$$r_{\psi\otimes\phi}^*
\circ 1_{{\cal
F}(\bar\sigma\otimes\bar\rho)}\otimes E_{\rho,\sigma}\circ
E_{\bar\sigma\otimes\bar\rho,\rho\otimes\sigma}
{\cal F}(1_{\bar\sigma}\otimes R_\rho\otimes 1_\sigma\circ R_\sigma).$$
On the other hand by $(3.9)$,
$$1_{{\cal
F}(\bar\sigma\otimes\bar\rho)}\otimes E_{\rho,\sigma}\circ
E_{\bar\sigma\otimes\bar\rho,\rho\otimes\sigma}=
E_{\bar\sigma\otimes\bar\rho,\rho}\otimes 1_{{\cal F}(\sigma)}\circ
E_{\bar\sigma\otimes\bar\rho\otimes\rho,\sigma},$$
hence the last term above equals
$$r_{\psi\otimes\phi}^*\circ
E_{\bar\sigma\otimes\bar\rho,\rho}\otimes 1_{{\cal F}(\sigma)}\circ
E_{\bar\sigma\otimes\bar\rho\otimes\rho,\sigma}\circ{\cal
F}(1_{\bar\sigma}\otimes R_\rho\otimes 1_\sigma)\circ {\cal
F}(R_\sigma).\eqno(3.13)$$
Now by $(3.11)$,
$${\cal F}(1_{\bar\sigma}\otimes R_\rho^*\otimes 1_{\sigma})
\circ E_{\bar\sigma\otimes\bar\rho\otimes\rho,\sigma}=
{\cal F}(1_{\bar\sigma}\otimes R_\rho^*)\otimes 1_{{\cal
F}(\sigma)}\circ E_{\bar\sigma\otimes\bar\rho\otimes\rho,\sigma}=$$
$$E_{\bar\sigma,\sigma}\circ{\cal F}(1_{\bar\sigma}\otimes
R_\rho^*)\otimes 1_{{\cal
F}(\sigma)}\circ E_{\bar\sigma\otimes\bar\rho\otimes\rho,\sigma},$$
hence taking the adjoint,
$$
E_{\bar\sigma\otimes\bar\rho\otimes\rho,\sigma}\circ
{\cal F}(1_{\bar\sigma}\otimes R_\rho\otimes 1_{\sigma})=
{\cal F}(1_{\bar\sigma}\otimes R_\rho)\otimes 1_{{\cal
F}(\sigma)}\circ
E_{\bar\sigma,\sigma}=$$
$$1_{{\cal F}(\bar\sigma)}\otimes{\cal F}(R_\rho)\otimes 1_{{\cal
F}(\sigma)}\circ E_{\bar\sigma,\sigma},$$
and $(3.13)$ becomes
$$r_{\psi\otimes\phi}^*\circ
E_{\bar\sigma\otimes\bar\rho,\rho}\otimes 1_{{\cal F}(\sigma)}\circ
1_{{\cal F}(\bar\sigma)}\otimes {\cal F}(R_\rho)\otimes 1_{{\cal
F}(\sigma)}\circ
E_{\bar\sigma,\sigma}\circ {\cal
F}(R_\sigma).\eqno(3.14)$$
Now $$E_{\bar\sigma\otimes\bar\rho,\rho}\circ
E_{\bar\sigma.\bar\rho\otimes\rho}=1_{{\cal F}(\bar\sigma)}\otimes
E_{\bar\rho,\rho}\circ E_{\bar\sigma,\bar\rho\otimes\rho}$$
and $$E_{\bar\sigma,\bar\rho\otimes\rho}\circ 1_{{\cal
F}(\bar\sigma)}\otimes{\cal F}(R_\rho)=1_{{\cal
F}(\bar\sigma)}\otimes {\cal F}(R_\rho)$$
so $(3.14)$ equals 
$$r_{\psi\otimes\phi}^*\circ 1_{{\cal F}(\bar\sigma)}\otimes
\hat{R_\rho}\otimes 1_{{\cal F}(\sigma)}\circ \hat{R_\sigma}=
J_\sigma\phi\otimes J_\rho\psi
$$
as required.

\end{section}

\begin{section}{Concrete quasitensor functors}

Let $G$ be a compact quantum group. 
In this section we construct 
quasitensor $^*$--functors $\text{Rep}(G)\to{\cal H}$ which associate
a subspace $K_u$ of the representation Hilbert
space
$H_u$ with the representation $u$, and generalize the invariant vectors functor
associated with a compact quantum subgroup $K$ of a compact quantum group
$G$ described in Example 3.6.

For each unitary representation $u$ of $G$, we suppose assigned 
a subspace
$K_u$ of the representation Hilbert space $H_u$, and we look 
for sufficient  conditions on the projection maps $E_u:H_u\to K_u$.
\medskip

\noindent{\bf 4.1 Lemma} {\sl
 For $u\in\text{Rep}(G)$, let $E_u:
H_u\to H_u$ be an orthogonal projection such that
$$E_\iota=1_{\mathbb C},\eqno(4.1)$$ 
$$TE_u=E_vT,\quad T\in(u,v),\eqno(4.2)$$
$$E_u\otimes E_v=I\otimes E_v\circ
E_{u\otimes v}.\eqno(4.3)$$
Then 
the $^*$--functor:
$$K_u:=E_uH_u,$$
$$K_T:=T\upharpoonright_{K_u}\in(K_u, K_v)$$
is quasitensor.}\medskip

\noindent{\bf Proof} Property $(3.1)$ 
follows from $(4.1)$. Multiplying  both sides of $(4.3)$ on the right by
$E_{u\otimes v}$ and taking the adjoint gives
$$E_u\otimes E_v\circ E_{u\otimes v}=E_u\otimes E_v=E_{u\otimes v}\circ
E_u\otimes E_v,$$ so the subspace $K_u\otimes K_v$ of $H_u\otimes H_v$ 
is contained in $K_{u\otimes v}$, and this shows $(3.7)$. 
We show $(3.11)$.
For $S\in(u,u')$, $T\in(v,v')$:
 $$K_{S\otimes T}\upharpoonright_{K_u\otimes K_v}=((S\otimes
T)\upharpoonright_{K_{u\otimes v}})\upharpoonright_{K_u\otimes K_v}=$$
$$(S\otimes T)\upharpoonright_{K_u\otimes K_v}=K_S\otimes K_T.$$
In this case the projection $E_{u,v}: K_{u\otimes v}\to K_u\otimes K_v$
is given by the restriction of $E_u\otimes E_v$ to $K_{u\otimes v}$, so
$(3.8)$ follows easily. Since the copy of $K_u\otimes K_v\otimes K_z$
sitting inside  $K_u\otimes K_{v\otimes z}$ and $K_{u\otimes v}\otimes
K_z$ is, in both cases, the subspace $K_u\otimes K_v\otimes K_z$ of
$H_u\otimes H_v\otimes H_z$, $(3.9)$ is satisfied. It remains to check
$(3.10)$.
If $H$ and $H'$ are
Hilbert spaces, we consider the operators of tensoring on the right 
 of $H$ by vectors $\phi\in H'$:
$$r_\phi: \psi\in H\to\psi\otimes \phi\in H\otimes H'.$$
Now $(4.3)$ implies that $E_u\otimes E_v\circ E_{u\otimes v}=I\otimes
E_v\circ E_{u\otimes v}$. Thus 
if $\psi\in K_{u\otimes v}$ and $\phi\in K_v$
then 
$r_\phi^*(\psi)\in K_u$.
If then $\psi'\in
K_z$ then
$\psi'\otimes(r_\phi^*\psi)\in K_{z}\otimes K_u$. Now choose
$\phi=\phi_k$,
an orthonormal basis of $K_v$, apply the operator $r_{\phi_k}$,
 use the fact that
on $K_{z\otimes u\otimes v}$, $\sum_k r_{\phi_k}r_{\phi_k}^*$ 
converges strongly
to the projection map $E_{z\otimes u,v}$ onto $K_{z\otimes u}\otimes K_v$
and obtain the desired relation.
\medskip

The property that the isometric inclusion map $K_u\otimes K_v\subset
K_{u\otimes v}$ be 
simply the identity map is well described by the notion of quasitensor
natural transformation.\medskip

\noindent{\bf 4.2 Definition} Let ${\cal F}$, ${\cal G}$ be quasitensor
$^*$--functors  from a tensor $C^*$--category ${\cal T}$ to the
tensor
$C^*$--category of
Hilbert spaces ${\cal H}$. Let $S^{\cal F}_{\rho,\sigma}$ and
$S^{\cal G}_{\rho,\sigma}$ 
be the defining set of isometries for ${\cal F}$ and ${\cal G}$
respectively.
A
natural
transformation $\eta:{\cal F}\to{\cal
G}$ will be called {\it quasitensor} if for objects $\rho,\sigma\in{\cal
T}$,
$$\eta_\iota:{\cal F}(\iota)={\mathbb C}\to{\cal G}(\iota)={\mathbb
C}\text{
is
the identity map }$$
$$\eta_{\rho\otimes\sigma}\circ
S^{\cal F}_{\rho,\sigma}=S^{\cal
G}_{\rho,\sigma}\circ\eta_\rho\otimes\eta_\sigma.$$
\medskip

\noindent{\bf 4.3 Proposition} {\sl Let $K:\text{Rep}(G)\to{\cal H}$ be 
the quasitensor $^*$--functor obtained from projections $(E_u)$ satisfying
properties $(4.1)$--$(4.3)$. Then
the inclusion map
$W_u: K_u\to H_u$ defines a  quasitensor natural transformation
from the functor $K$ to the embedding functor
$H:\text{Rep}(G)\to{\cal H}$.}\medskip

We next show that the invariant vectors functors associated with quantum
subgroups, fit into this description.
Let  $K$ be a compact quantum subgroup of $G$ with Haar measure $h'$, and,
for an invertible
representation $u$ of $G$,
let $E_u^K: H_u\to H_u$ be 
the  idempotent onto the subspace of 
$u\upharpoonright_K$--fixed vectors obtained averaging over the
$K$--action: $E_u^K(\psi)=\iota\otimes h'\circ u\upharpoonright_K(\psi)$.
If $u$ is unitary, $E_u^K$ is
a selfadjoint projection. Indeed, 
$$(E_u^K(\psi_i),
\psi_j)=\overline{h'(\pi(u_{ji}))}=h'(\pi(u_{j i}^*))=$$
$$h'(\kappa'(\pi(u_{i j})))=h'(\pi(u_{i j}))=(\psi_i, E_u^K(\psi_j)),$$
since the Haar measure $h'$ of $K$ is left invariant by the coinverse
$\kappa'$ \cite{Wcmp}.
 \medskip

\noindent{\bf 4.4. Lemma} {\sl If $u$ and $v$ are invertible 
and finite dimensional $G$--representations, and $K$ is a compact quantum
subgroup  of $G$ then for any
$T\in(u\upharpoonright_K,v\upharpoonright_K)$, thus in particular for any
$T\in(u,v)$,
$$T\circ E_u^K=E_{v}^K\circ T.$$}\medskip

\noindent{\bf Proof} 
For $\psi\in H_u$, $$T\circ E_u^K(\psi) =T\iota\otimes h'\circ
u\upharpoonright_K(\psi)=$$
$$\iota\otimes
h'(T\otimes Iu\upharpoonright_K(\psi))=\iota\otimes
h'(v\upharpoonright_K(T(\psi)))=E_{v}^K\circ T(\psi).$$
\medskip

We thus obtain another proof of the fact that the invariant vectors
functor is quasitensor. \medskip

\noindent{\bf 4.5 Theorem} {\sl If 
$K$ is a compact quantum subgroup of $G$,
the projections 
$u\to E^K_u$ satisfy properties $(4.1)$--$(4.3)$.
Therefore the  associated invariant vectors functor 
is a  quasitensor $^*$--functor and the inclusion map
$K_u\to H_u$ is a quasitensor natural transformation from the functor $K$
to the
embedding functor $H$.}\medskip

\noindent{\bf Proof} 
We check properties $(4.1)$--$(4.3)$ on the projections $E^K_u$. Clearly
$(4.1)$ is verified, and $(4.2)$ follows from the previous lemma. 
Since the tensor product of $K$--invariant vectors is a $K$--invariant
vector
for the tensor product representation, we have: $E^K_u\otimes E^K_v\leq
E^K_{u\otimes v}$, and   $(4.3)$ 
becomes equivalent to
$$E^K_u\otimes E^K_v E^K_{u\otimes
v}=I\otimes E^K_v\circ E^K_{u\otimes v},$$
or, in other words, to the 
fact that if
$\psi\in K_{v}$, $\eta\in K_{u\otimes v}$ then 
$r_\psi^*(\eta)\in K_u$. In order to verify this property, we use  the
fact that
the restriction to $K$ of a tensor product representation is the tensor
product of the restrictions, so
$$\eta\in K_{u\otimes v}=(\iota, (u\otimes v)\upharpoonright_K)=(\iota,
u\upharpoonright_K\otimes v\upharpoonright_K).$$
Using 
 $\psi\in K_v=(\iota, v\upharpoonright_K)$, too, gives,
$$r_\psi^*(\eta)=1_{u\upharpoonright_K}\otimes\psi^*\circ\eta\in(\iota,
u\upharpoonright_K)=K_u,$$
where the operations are to be understood in the tensor $C^*$--category 
Rep$(K)$.\medskip

\noindent{\it  Remark}
Notice that a similar argument shows that the projections $E^K_u$ also satisfy
$$E_u^K\otimes E_v^K=E_u^K\otimes I\circ E_{u\otimes v}^K.\eqno(4.4)$$

\noindent{\it Remark}
Notice that if $G$ is a  {\it maximal} quantum group (i.e. obtained
from its smooth part ${\cal A}_\infty$ by completing with 
respect to the maximal $C^*$--seminorm) and $K$ is
the trivial subgroup (corresponding to the counit $e:{\cal A}\to{\mathbb
C}$ of $G$) then the associated functor $K$ coincides with the embedding 
functor $H: \text{Rep G}\to{\cal H}$, which is tensor. 
At the other extreme, if $K=G$, $K_u=(\iota, u)$, and this is {\it not}
a tensor functor,
as  $K_u=K_{\overline{u}}=0$, for example, if $u$ is irreducible,
but $K_{\overline{u}\otimes u}=(\iota,\overline{u}\otimes u)\neq0$.
\medskip

We next construct certain maps which will be useful later on when describing
multiplicities of spectral 
representations in quantum
quotient spaces.
Given a unitary representation $u:H_u\to H_u\otimes{\cal A}$, we set
$$u^K:=E^K_u\otimes I\circ u: H_u\to K_u\otimes{\cal A}.$$
For any pair of vectors $\phi,\psi\in H_u$, consider the coefficients
of $u^K$:
$$u^K_{\phi,\psi}:=\ell_{\phi}^*\circ
u^K(\psi)=\ell_{E^K_u(\phi)}^*\circ u(\psi)=u_{E^K_u(\phi),\psi},$$
which belong to ${\cal A}^\delta_{sp}$.
So  the range of $u^K$ is actually contained in 
$K_u\odot{\cal A}^\delta_{sp}$:
$$u^K:H_u\to K_u\odot{\cal A}^\delta_{sp}.$$
The map $u\to u^K$ satisfies the following properties.
\medskip

\noindent{\bf 4.6 Proposition} {\sl Let $u$ and $v$ be unitary
representations
of $G$. 
\begin{description}
\item {\rm a)}
For any intertwiner $T\in(u, v)$,
$$K_{T}\otimes I \circ u^K=v^K\circ T,$$
\item{\rm b)}
For $\phi\in K_u$, $\phi'\in K_{v}$, $\psi\in H_u$, $\psi'\in H_v$,
$$(u\otimes v)^K_{{\phi\otimes\phi'},
\psi\otimes\psi'}=u^K_{\phi,\psi}u^K_{\phi',\psi'}.$$
\end{description}

}\medskip

\noindent{\bf Proof} a)
$$v^K\circ T=E^K_v\otimes I\circ v\circ T=E^K_v\otimes I\circ T\otimes
I\circ u=$$
$$T\otimes I\circ E^K_u\otimes I\circ u=K_{T}\otimes I u^K.$$
Property b) follows from the fact that $K_u\otimes K_v\subset K_{u\otimes
v}$, so
$$(u\otimes v)^K_{\phi\otimes\phi',\psi\otimes\psi'}=u\otimes 
v_{\phi\otimes\phi',\psi\otimes\psi'}=$$
$$u_{\phi,\psi}v_{\phi',\psi'}=
u^K_{\phi,\psi}v^K_{\phi',\psi'}.$$
\medskip

\end{section}

\begin{section} {Characterizing the invariant vectors functor}

For each $u\in\text{Rep}(G)$, let $E_u$ be an orthogonal projection 
on the representation Hilbert space $H_u$ satisfying properties
$(4.1)$--$(4.3)$. In terms of the Hilbert spaces $K_u=E_uH_u$, these
conditions can be written
$$K_\iota=H_\iota={\mathbb C},\eqno(5.1)$$
$$TK_u\subset K_v.\quad T\in(u,v),\eqno (5.2)$$
$$r_k^*K_{u\otimes v}\subset K_u,\quad k\in K_v,\eqno(5.3)$$
$$K_u\otimes K_v\subset K_{u\otimes v}.\eqno(5.4)$$
As a consequence of these conditions, one gets
$$(\iota,u\otimes v)\subset
K_{u\otimes v}.\eqno(5.5)$$ In fact, pick
$T\in(\iota, u\otimes v)$.
Properties $(5.1)$ and $(5.2)$ show that for any $\lambda\in
K_\iota={\mathbb
C}$, $T\lambda\in K_{u\otimes v}$, so
$T=T1\in K_{u\otimes v}$.
In particular, thanks to $(5.3)$,
$$r_k^*T\in K_u,\quad T\in (\iota,u\otimes v),\quad  k\in K_v.\eqno(5.6)$$
Now let $\overline{u}$ be a conjugate of $u$ defined by
$R\in(\iota,\overline{u}\otimes u)$, $\overline{R}\in(\iota,
u\otimes\overline{u})$, and let $j: H_u\to H_{\overline{u}}$ be the
associated antilinear invertible intertwiner defined by:
$$R=\sum_i j\phi_i\otimes \phi_i,$$
$$\overline{R}=\sum_j j^{-1}\psi_j\otimes \psi_j,$$
with $(\phi_i)$ and $(\psi_j)$ orthonormal bases of $H_u$ and
$H_{\overline{u}}$
respectively, that we choose to complete the orthonormal bases of the
corresponding
subspaces $K_u$ and $K_{\overline{u}}$. Property $(5.6)$ applied to 
$R$ and $\overline{R}$ shows that
$$jK_u=K_{\overline{u}}.\eqno (5.7)$$
When $G$ is a group, conditions $(5.1)$, $(5.2)$, $(5.4)$ and $(5.7)$ 
are known to characterize a subgroup of
$G$ with the property that each $K_u$ is 
the invariant subspace of $u$ in restriction to the subgroup
\cite{Roberts}.

   We next show that the antilinear invertible 
operator $j^*$ defined by a solution of the conjugate equations may be 
regarded as an antilinear intertwiner from $\bar u$ to $u$. The 
intertwining property of $\overline R$, when expressed in coordinate 
form, reads 
$$\sum_{i,p}u_{ji}\bar u_{qp}\overline R_{(ip)}=\overline R_{(jq)}.$$ 
Expressed in terms of $j$ and using the unitarity of $u$, we get 
$$\sum_p\bar
u_{qp}(\psi_p,{j^{-1}}^*\varphi_k)=\sum_ju^*{}_{kj}(\psi_q,{j^{-1}}^*\varphi_j).$$ 
Now, if $A\in\cA$, 
$$(\psi_q\otimes
A,{j^{-1}}^*\otimes^*u\varphi)=\sum_{i,j}(\psi_q,{j^{-1}}^*\varphi_j)A^*u^*_{ji}
(\varphi,\varphi_i)=\sum_pA^*\bar u_{qp}(\psi_p,{j^{-1}}^*\varphi)$$ 
$$=\sum_p((\psi_q\otimes A,\bar
u\psi_p)(\psi_p,{j^{-1}}^*\varphi)=(\psi_q\otimes A,\bar 
u{j^{-1}}^*\varphi).$$ 
Thus ${j^{-1}}^*\otimes^*u=\bar u{j^{-1}}^*$ and therefore
$j^*\otimes^*\bar u=uj^*$,
which is the
form of the intertwining 
relation for an antilinear operator. Note that this intertwining 
relation could alternatively 
be deduced as above by starting with $R$ rather than $\overline R$.\smallskip

Now, in view of what happens for compact groups, the question naturally
arises of whether, given a compact quantum group $G$ and 
subspaces $K_u\subset H_u$ for each finite--dimensional unitary
representation $u$ of $G$, conditions $(5.1)$, $(5.2)$, $(5.4)$ and
$(5.7)$ are still
sufficient for the existence of a
unique compact quantum subgroup of $G$ whose 
the subspaces of invariant vectors are the $K_u$.

Uniqueness does not hold in general for compact quantum
groups.
In fact, one
can have a proper subgroup $K$ of a compact quantum group $G$
with $K_u=(\iota,u)$ for any representation $u$ of $G$:
consider a  group $G$ with a
nonfaithful Haar measure $h$ and form
the reduced 
group $G_{\text{red}}=({\cal A}_{\text{red}},\Delta_{\text{red}})$
obtained
completing the dense Hopf $^*$--subalgebra
 in the norm
defined by the GNS representation $\pi_h$. 
Since $\pi_h:{\cal A}\to{\cal A}_{\text{red}}$ is a surjection
intertwining
the corresponding coproducts, $G_{\text{red}}$ becomes a subgroup of $G$.
Since $h$ is not faithful,  $\pi_h$ has a nontrivial kernel. In this sense,
$G_{\text{red}}$ is a proper subgroup of $G$.
Furthermore $G_{\text{red}}$ has a faithful Haar measure, while $G$ has 
not, so $G$ and $G_{\text{red}}$ are not isomorphic as Hopf
$C^*$--algebras. 
However, $G$ and $G_{\text{red}}$
have the same
representation categories, and therefore
the same  spaces of invariant
vectors. 

We could have also considered
the maximal compact quantum group,
$G_{\text{max}}=({\cal A}_{\text{max}},\Delta_{\text{max}})$ obtained
completing the dense Hopf $^*$--subalgebra with
respect to the maximal $C^*$--seminorm. There is again a surjection
$\pi:{\cal A}_{\text{max}}\to{\cal A}$ intertwining the coproducts,
so $G$ is a subgroup of $G_{\text{max}}$, and we still have
$$\text{Rep}(G_{\text{max}})=\text{Rep}(G)=\text{Rep}(G_{\text{red}}),$$
so $G_{\text{red}}$ is also a proper subgroup of $G_{\text{max}}$ with
the same spaces of invariant vectors as $G_{\text{max}}$.
This example suggests that one possible way to get uniqueness is to
restrict attention to maximal compact quantum groups.

As far as existence is concerned, we start by showing the following result.
\medskip

\noindent{\bf 5.1 Theorem} {\sl Let $G=({\cal A},\Delta)$ be a
compact quantum group,
and, for each $u\in\text{Rep}(G)$, let $K_u$
be a subspace of the representation Hilbert space $H_u$ satisfying 
conditions $(5.1)$, $(5.2)$, $(5.4)$ and $(5.7)$. Then 
there exists
a 
compact quantum subgroup $K$ of $G$ such that $K_u\subset (\iota_K,
u\upharpoonright_K)$ for $u\in\text{Rep}(G)$ and such that
the linear span of 
$\{u_{k,\phi}, k\in K_u,\phi\in H_u, u\in\text{Rep}(G)\}$ is a unital
$^*$--subalgebra of the
quantum quotient space  ${\cal
A}^\delta$.}\medskip

\noindent{\bf Proof} Consider a complete set $A$ of irreducible
representations $u$ of $G$ such that $K_u\neq0$. Let $M$ denote the linear
span generated by the set
$$\{x^u_{\phi,k}:=u_{\phi,k}-(\phi,k)I, u\in A, k\in K_u, \phi\in H_u\},$$
in the Hopf $C^*$--algebra ${\cal A}$. It is easy to check that
$$\Delta(x^u_{\phi,k})=\sum_i u_{\phi,\phi_i}\otimes
x^u_{\phi_i,k}+x^u_{\phi,k}\otimes I,$$
with $(\phi_i)$ an orthonormal basis of $H_u$.
Therefore $\Delta(M)\subset {\cal A}_\infty\odot M+M\odot {\mathbb C}I$.
Let ${\cal J}$ be the closed two--sided ideal of ${\cal A}$ generated by 
$M$. Then $\Delta({\cal J})\subset {\cal A}\otimes{\cal J}+{\cal
J}\otimes{\cal A}$, so ${\cal J}$ is a closed bi--ideal. Consider the
associated 
compact quantum subgroup $K=({\cal A}/{\cal J}, \Delta')$ of $G$
with coproduct $\Delta'(q(a))=q\otimes q\circ\Delta(a)$, where $q:{\cal
A}\to{\cal A}/{\cal J}$ is the canonical surjection.
We show that $K_u\subset (\iota_K, u\upharpoonright_K)$. For $k\in K_u$:
$$u\upharpoonright_K(k)=\iota\otimes q\circ u(k)=$$
$$\iota\otimes
q(\sum_i\phi_i\otimes u_{\phi_i,k})=\sum_i\phi_i\otimes q(u_{\phi_i,k})=$$
$$\sum_i \phi_i\otimes (\phi_i, k)I=k\otimes I.$$ We are left to show 
that the linear span $V$ of  all the $u_{k,\phi}$ is a
unital $^*$--subalgebra 
of ${\cal A}^\delta$. Since $K_u\subset (\iota, u\upharpoonright_K)$, 
$V$ is contained in ${\cal A}^\delta$, and $I\in V$, as $K_\iota={\Bbb
C}$. Therefore it suffices to show that
$V$ is a $^*$--subalgebra of ${\cal A}$. On the other hand the
$^*$--algebra structure of ${\cal A}$ recalled at the end of subsection
2.1 and properties $(5.4)$ and $(5.7)$ 
show that $V$ is a $^*$--subalgebra. 
\medskip

\noindent{\it Remark}
Under conditions $(5.1)$, $(5.2)$, $(5.4)$ and $(5.7)$ alone,
one
can not, in general, 
identify the subspace $K_u$ with the space of {\it all} the invariant
vectors 
$(\iota_K, u\upharpoonright_K)$ for some quantum subgroup $K$ of $G$. 
In fact there is an example, due to Wang \cite{Wang} of an ergodic 
action $\delta$ on a commutative $C^*$--algebra ${\cal C}$ which is not a
quotient
action.
Now the commutativity of ${\cal C}$ allows a faithful
embedding of $({\cal C}_{sp},\delta)$ into a quantum quotient space by a 
quantum subgroup and a construction of subspaces $K_u$ satisfying the
above equations (see Theorem 11.4). Such an embedding, though, can not
extend to an isomorphism of the whole of ${\cal C}$ and therefore  the
$K_u$ can not be the spaces of all the invariant vectors.

We shall be able, though,  to give a positive answer if we replace $(5.4)$
by a stronger, still necessary, condition, motivated by the following
argument.
\medskip

In the group case, the representation category of $G$ contains,
among its intertwiners, the {\it permutation symmetry}: the 
operators $\theta_{u,v}\in(u\otimes v,v\otimes u)$
permute the order of factors in the tensor product.
Consequently, for $u,v,z\in\text{Rep}(G)$, 
$$1_u\otimes\phi\otimes1_z\circ k
\in K_{u\otimes v\otimes z},\quad
k\in K_{u\otimes z},
\phi\in K_v.\eqno(5.8)$$
In fact, by $(5.4)$, $\phi\otimes k\in K_{v\otimes u\otimes z}$, so
$1_u\otimes\phi\otimes1_z\circ k=(\vartheta_{v,u}\otimes 1_z)\phi\otimes
k$
and this is an element of $K_{u\otimes v\otimes z}$ thanks to $(5.2)$.

On the other hand, 
$(5.8)$ is still a 
necessary condition
for the $K_u$ to be the invariant subspaces 
of the restriction of $u$ to a quantum subgroup, 
 as one can easily show using the same argument
as the one used,
in Theorem 4.5, to show the necessity of condition $(4.3)$.
Therefore, in the quantum group case, it seems natural to replace $(5.4)$
by the stronger
condition $(5.8)$.

Assume then that  we have a functor $K$
associating to any representation $u\in\text{Rep}(G)$ a subspace
$K_u\subset H_u$ satisfying $(5.1)$, $(5.2)$, $(5.7)$ and $(5.8)$.
We  consider, for $u,v\in\text{Rep}(G)$, the subspace of $(H_u, H_v)$
defined by
$$<H_u,H_v>:=\{\overline{R}^*\otimes 1_v\circ1_u\otimes\phi,
\phi\in
K_{\overline{u}\otimes v}\},$$
where $\overline{R}\in(\iota,u\otimes \overline{u})$ is an
intertwiner arising from
a solution of the conjugate equations for $u$ in $\text{Rep}(G)$.
\medskip

\noindent{\bf 5.2 Proposition} {\sl The space $<H_u,H_v>$ is independent
of the
choice of $\overline{u}$ and $\overline{R}$.}\medskip

\noindent{\bf Proof} If $R'\in(\iota, u\otimes \tilde{u})$ arises
from
another solution to the conjugate equations then there exists a unitary 
$U\in(\overline{u},\tilde{u})$ such that $R'=1_u\otimes U\circ
\overline{R}$
\cite{LongoRoberts}. Since $U\otimes 1_v\in (\overline{u}\otimes v,
\tilde{u}\otimes v)$, 
$U\otimes 1_vK_{\overline{u}\otimes v}=K_{\tilde{u}\otimes v}$ by $(5.2)$.
Therefore any $\phi'\in K_{\tilde{u}\otimes v}$ is of the form 
$U\otimes 1_v\phi$ with $\phi\in K_{\overline{u}\otimes v}$. This  
implies
${R'}^*\otimes 1_v\circ1_u\otimes\phi'=
{\overline{R}}^*\otimes 1_v\circ1_u\otimes\phi$.
\medskip

\noindent{\bf 5.3 Proposition} {\sl For $u,v\in\text{Rep}(G)$,
$(u,v)\subset<H_u,H_v>$.}\medskip

\noindent{\bf Proof} 
Fix a pair ${R}\in(\iota,\overline{u}\otimes {u})$ 
$\overline{R}\in(\iota,
u\otimes\overline{u})$ 
solving the conjugate equations.
Pick $T\in(u,v)$ and set 
$\phi_T:=1_{\overline{u}}\otimes T
\circ{R}$, 
So $\phi_T\in (\iota, \overline{u}\otimes v)\subset
K_{\overline{u}\otimes v}$.
 Since
$$\overline{R}^*\otimes 1_v\circ 1_u\otimes\phi_T=\overline{R}^*\otimes
1_v\circ
1_{u\otimes\overline{u}}\otimes T\circ1_u\otimes R=$$
$$T\circ \overline{R}^*\otimes
1_u\circ 1_u\otimes R=T,$$
we can conclude that $T\in<H_u,H_v>$.\medskip

\noindent{\bf 5.4 Proposition} {\sl The Hilbert spaces $H_u$ and the 
subspaces $<H_u,H_v>$, as $u$ and $v$
vary in $\text{Rep}(G)$, form, respectively, the objects and  the
 arrows of a tensor
$^*$--subcategory
of the category of Hilbert spaces. This category contains the image of
$\text{Rep}(G)$ under the embedding functor $H$, therefore it has
conjugates.
}\medskip

\noindent{\bf Proof} Since $(u,v)\subset <H_u,H_v>$, $1_u\in<H_u,H_u>$ for
$u\in\text{Rep}(G)$. We show that if $T\in<H_u,H_v>$ and $S\in<H_v,H_z>$
then
$S\circ T\in <H_u,H_z>$. In fact, writing 
$T= \overline{R}_u^*
\otimes 1_v
\circ 1_u\otimes \phi$ 
and
$S= 
\overline{R}_v^*
\otimes 1_z
\circ 1_v\otimes \psi$ 
then 
$$S\circ T=
\overline{R}_v^*
\otimes 1_z
\circ 1_v\otimes \psi\circ \overline{R}_u^*
\otimes 1_v
\circ 1_u\otimes \phi=$$
$$\overline{R}_v^*\otimes 1_z\circ\overline{R}_u^*\otimes
1_{v\otimes\overline{v}\otimes
z}\circ
1_{u\otimes\overline{u}\otimes v}\otimes\psi\circ 1_u\otimes\phi=$$
$$\overline{R}_u^*\otimes 1_z\circ
1_{u\otimes\overline{u}}\otimes\overline{R}_v^*\otimes1_z\circ 1_u\otimes
(\phi\otimes\psi)=$$
$$\overline{R}_u^*\otimes 1_z\circ 1_u\otimes(
1_{\overline{u}}\otimes\overline{R}_v^*\otimes
1_z(\phi\otimes\psi)
)\in<H_u,H_z>$$
as $\phi\otimes\psi\in K_{\overline{u}\otimes v\otimes \overline{v}\otimes
z}$
and 
$(1_{\overline{u}}\otimes\overline{R}_v\otimes
1_z)^*(\phi\otimes\psi)\in K_{\overline{u}\otimes z}$ by $(5.2)$.

We next show that $<H_u,H_v>^*=<H_v,H_u>$. Pick $\phi\in
K_{\overline{u}\otimes{v}}$
and set 
$T=\overline{R}_u^*\otimes 1_v\circ1_u\otimes \phi$. Then
$T^*=1_u\otimes \phi^*
\circ \overline{R}_u\otimes 1_v$.
We use the solution  $R_{\overline{u}\otimes v}:= 1_{\overline{v}}\otimes
\overline{R}_u\otimes 1_v\circ R_v\in(\iota, \overline{v}\otimes
u\otimes\overline{u}\otimes
v)$ and 
$\overline{R}_{\overline{u}\otimes v}:=1_{\overline{u}}\otimes
\overline{R}_v\otimes
1_u\circ
R_u\in(\iota, \overline{u}\otimes v\otimes\overline{v}\otimes u)$ of the
conjugate
equations for $\overline{u}\otimes v$. 
 Thanks to $(5.7)$,
$\psi:=1_{\overline{v}\otimes u}\otimes \phi^*R_{\overline{u}\otimes
v}=j\phi\in
K_{\overline{v}\otimes u}$.
We have: 
$$\overline{R}_v^*\otimes 1_u\circ 1_v\otimes\psi=\overline{R}_v^*\otimes
1_u\circ
1_{v\otimes\overline{v}\otimes u}\otimes\phi^*\circ
1_{v\otimes\overline{v}}\otimes
\overline{R}_u\otimes 1_v\circ 1_v\otimes R_v=$$
$$1_u\otimes\phi^*\circ \overline{R}_u\otimes 1_v\circ
\overline{R}_v^*\otimes
1_v\circ
1_v\otimes R_v=1_u\otimes\phi^*\circ\overline{R}_u\otimes 1_v=T^*.$$
Therefore $T^*\in<H_v,H_u>$.

We are left to show that if $S=\overline{R}_u^*\otimes
1_{u'}\circ 1_u\otimes\phi\in<H_u,H_{u'}>$, $T=\overline{R}_v^*\otimes 
1_{v'}\circ 1_v\otimes\psi\in<H_v,H_{v'}>$ then $S\otimes T\in<H_{u\otimes
v},
H_{u'\otimes v'}>$.
Consider the following solution to the conjugate equations for $u\otimes
v$:
${R}_{u\otimes v}=1_{\overline{v}}\otimes{R}_u\otimes
1_{{v}}\circ {R}_v$,
$\overline{R}_{u\otimes v}=1_u\otimes\overline{R}_v\otimes
1_{\overline{u}}\circ \overline{R}_u$.
Since $\eta:=1_{\overline{v}}\otimes\phi\otimes1_{v'}\circ \psi\in
K_{\overline{v}\otimes\overline{u}\otimes u'\otimes v'}$ by 
$(5.8)$ then 
$\overline{R}_{u\otimes v}^*\otimes 1_{u'\otimes v'}\circ 1_{u\otimes
v}\otimes\eta
\in <H_{u\otimes v}, H_{u'\otimes v'}>$. 
On the other hand,
$$\overline{R}_{u\otimes v}^*\otimes 1_{u'\otimes v'}\circ 1_{u\otimes
v}\otimes\eta=$$
$$\overline{R}_u^*\otimes 1_{u'\otimes v'}\circ 
1_u\otimes
\overline{R}_v^*\otimes 1_{\overline{u}\otimes u'\otimes v'}\circ
1_{u\otimes v\otimes \overline{v}}\otimes\phi\otimes 1_{v'}\circ
1_{u\otimes v}\otimes \psi=$$
$$\overline{R}_u^*\otimes 1_{u'\otimes v'}\circ 1_u\otimes\phi\otimes
1_{v'}\circ 1_u\otimes\overline{R}_v^*\otimes 1_{v'}\circ
1_{u\otimes v}\otimes \psi=$$
$$S\otimes 1_{v'}\circ 1_u\otimes T=S\otimes
T.$$
\medskip

A compact quantum subgroup  of $G$ will be called maximal if it is
maximal as a
compact quantum group (recall that this means that the norm on the dense
Hopf $^*$--subalgebra coincides with the maximal $C^*$--seminorm).
The above results
lead to the following theorem.\medskip

\noindent{\bf 5.5 Theorem} {\sl Consider the tensor $C^*$--category of
finite
dimensional unitary representations of a maximal compact quantum group
$G$.
Suppose for each representation $u$ of $G$ on a Hilbert space $H_u$ there
is given a subspace $K_u\subset H_u$ satisfying properties $(5.1)$,
$(5.2)$, $(5.7)$ and $(5.8)$. Then the $K_u$ are the subspaces of
invariant vectors for a unique  maximal compact quantum subgroup of $G$.}
\medskip

\noindent{\bf Proof} Consider the category ${\cal T}$ obtained
completing the spaces 
$<H_u,H_v>$, for $u,v\in\text{Rep}(G)$, with respect to 
subobjects and direct sums. This is still a tensor $^*$--category of
Hilbert spaces with conjugates, now with subobjects and direct sums.
Thanks to Woronowicz's 
Tannaka--Krein Theorem, we can construct  a Hopf $^*$--algebra 
$({\cal C},\Delta')$. Since $\text{Rep}(G)\subset{\cal T}$, this Hopf
$^*$--algebra is a model, in the sense of \cite{Wtk}, for the Hopf
$^*$--algebra constructed
from $\text{Rep}(G)$, which, in turn, is isomorphic to 
$({\cal A}_{\infty},\Delta)$.
 Therefore we can find a $^*$--epimorphism $\pi:{\cal A}_{\infty}\to{\cal
C}$
such that $\Delta'\circ\pi=\pi\otimes\pi\circ\Delta$. Set ${\cal
J}:=\text{ker}(\pi)$. This is obviously a $^*$--ideal, but also an
algebraic coideal, as if $a\in{\cal J}$ then
$b:=\iota\otimes\pi(\Delta(a))$ must belong to the kernel of
$\pi\otimes\iota_{\cal C}$. Since $\text{ker}(\pi\otimes\iota_{\cal C})=
{\cal J}\odot{\cal C}=\text{Image}(\iota_{\cal J}\otimes\pi)$, we can
find $c\in{\cal J}\odot{\cal A}$ such that $b=\iota_{\cal
J}\otimes\pi(c)$. Set $d:=\Delta(a)-c$, which is easily checked to lie in 
$\text{ker}(\iota_{\cal A}\otimes\pi)$, that in turn, equals ${\cal
A}\otimes{\cal J}$. Then $\Delta(a)=c+d\in{\cal J}\odot{\cal A}+{\cal
A}\odot{\cal J}$.
Consider the completion ${\cal A}'$ of ${\cal C}$ in the maximal
$C^*$--seminorm, so
$\Delta'$ extends to the completion. We thus obtain a
compact quantum group $({\cal
A}',\Delta')$ and a $^*$--homomorphism $\tilde{\pi}:{\cal
A}_{\infty}\to{\cal
A}'$
intertwining $\Delta$ with $\Delta'$ and with the same kernel as $\pi$,
because the inclusion ${\cal C}\subset {\cal A}'$ is faithful. Therefore 
we can extend $\tilde{\pi}$ to a  $^*$--homomorphism $q$
from ${\cal A}$ to ${\cal
A}'$ interwining the coproducts. This $^*$--homomorphism is a surjection,
as its range contains ${\cal C}$, which is dense in ${\cal A}'$. 
Thus $({\cal A}',\Delta')$ is a maximal compact quantum subgroup of $G$,
that we
denote, with abuse of notation, by $K$. Since
$q({\cal A}_{\infty})={\cal C}$, and for any quantum subgroup one always
has
$q({\cal A}_{\infty})={\cal
A}'_{\infty}$, we deduce that ${\cal C}={\cal A}'_{\infty}$. Therefore the
representation category of the subgroup $K$ coincides with the category we
started out
with: $\text{Rep}(K)={\cal T}$. In particular  for any representation $u$
of $G$,  $(\iota_K, u\upharpoonright_K)=<H_\iota, H_u>=K_u$.
If $K_1=({\cal A}_1,\Delta_1)$ is  another maximal quantum subgroup of $G$
with
spaces of invariant
vectors given by the $K_u$, then by Frobenius reciprocity, for any pair of
representations $u$ and $v$ of $G$,
$(u\upharpoonright_{K_1}, v\upharpoonright_{K_1})$ 
is canonically linearly
isomorphic to 
$$(\iota_K, \overline{u}\upharpoonright_{K_1}\otimes
v\upharpoonright_{K_1})=(\iota_K, (\overline{u}\otimes
v)\upharpoonright_{K_1})=K_{\overline{u}\otimes v}.$$
But 
$(u\upharpoonright_{K}, v\upharpoonright_{K})$ is also linearly
isomorphic,
according to the same isomorphism, to $K_{\overline{u}\otimes v}$, so
$(u\upharpoonright_{K_1}, v\upharpoonright_{K_1})=(u\upharpoonright_K,
v\upharpoonright_K)$. Hence $\text{Rep}(K_1)=\text{Rep}(K)$. 
It follows that there is a $^*$--isomorphism 
$\eta:{\cal A}_{1 \infty}\to{\cal A}'_{\infty}$ intertwining the
corresponding
coproducts, which extends to a $^*$--isomorphism of the completions.
\medskip

 We now give 
two constructions of subspaces $K_u\subset H_u$ satisfying our equations
making no mention of a compact quantum group. Let $\cT$ be a tensor 
$C^*$--subcategory with conjugates  of a tensor category of Hilbert spaces. The objects of 
$\cT$ will be denoted by $u,v,\dots$ as in the category of finite dimensional 
representations of a compact quantum group whilst the arrows will be written 
for example as $T\in(u,v)$. We suppose that $\cT$ comes equipped with a 
permutation symmetry $\varepsilon$. We now define 
$$K_u=\{\varphi\in H_u:\varepsilon(u,v)\varphi\otimes 1_v=1_v\otimes\varphi\,\, 
\text{for all objects}\, v\}.$$ 
Note that the defining condition can be replaced by 
$\varepsilon(v,u) 1_v\otimes\varphi=\varphi\otimes 1_v$ for all objects $v$. 
But more is true: we have 
$\varepsilon(u,v)\varphi\otimes\psi=\psi\otimes\varphi$ provided 
either $\varphi\in K_u$ and $\psi\in H_v$ or $\varphi\in H_u$ and $\psi\in K_v$. 
\smallskip 

  We check the validity of our equations $(5.1)$, $(5.2)$, $(5.7)$,  
$(5.4)$ for this definition of $K$. $(5.1)$ is 
true by definition. If $T\in (u,v)$ and $\varphi\in K_u$ then 
$$\varepsilon(v,w)T\varphi\otimes 1_w=1_w\otimes T\varepsilon(u,w)\varphi\otimes 1_w 
=1_w\otimes T\varphi.$$ 
Thus $(u,v)K_u\subset K_v$, which is $(5.2)$. If $\varphi\in K_u$ and 
$\psi\in K_v$ then 
$$\varepsilon(u\otimes v,w)\varphi\otimes\psi\otimes 1_w=\varepsilon(u,w)\otimes 1_v 
1_u\otimes \varepsilon(v,w)\varphi\otimes\psi\otimes 1_w=1_w\otimes\varphi\otimes\psi.$$ 
Thus $K_u\otimes K_v\subset K_{u\otimes v}$, which is $(5.4)$. Let 
$R\in(\iota,\bar u\otimes u)$ and $\overline R\in(\iota,u\otimes\bar u)$
be a solution 
of the conjugate equations and define the corresponding antilinear 
operator $j$ by: $j\varphi:=r_\varphi^*\circ R=1_{\bar
u}\otimes \varphi^* R$.  
If $\varphi\in K_u$ then 
$$j\varphi\otimes 1_v=
1_{\bar u}\otimes\varphi^*\otimes
1_v\circ
R\otimes 1_v=$$
$$1_{\bar u}\otimes\varphi^*\otimes
1_v\circ
\varepsilon(v, \bar u\otimes u)\circ 1_v\otimes R=
1_{\bar u}\otimes\varphi^*\otimes
1_v\circ
1_{\bar u}\otimes\varepsilon(v,u)\circ\varepsilon(v,\bar u)\otimes 1_u
\circ 1_v\otimes R=$$
$$1_{\bar u}\otimes 1_v\otimes \varphi^*
\circ\varepsilon(v,\bar u)\otimes 1_u
\circ 1_v\otimes R.$$
Acting on the left by $\varepsilon(\bar u,v)$, we get 
$$\varepsilon(\bar u,v)j\varphi\otimes 1_v=
1_v\otimes 1_{\bar u}\otimes\varphi^*\circ \varepsilon(\bar u,v)\otimes
1_u\circ \varepsilon(v,\bar u)\otimes 1_u\circ 1_v\otimes R=$$
$$1_v\otimes
j\varphi.$$
This proves $(5.7)$.\smallskip 

In fact there is a further property, $(5.8)$, valid here.
\smallskip 

$$1_u\otimes\phi\otimes 1_z\circ k\in K_{u\otimes v\otimes z},\quad k\in K_{u\otimes v},
\,\,\phi\in K_v.$$ 
In fact, 
$$\varepsilon(u\otimes v\otimes z,w)
(1_u\otimes\phi\otimes 1_z\circ k)\otimes 1_w
=$$
$$\varepsilon(u\otimes v,w)\otimes 1_z\circ 1_{u\otimes
v}\otimes\varepsilon(z,w)
\circ 1_u\otimes\phi\otimes 1_z\otimes 1_w\circ k\otimes 1_w=$$ 
$$\varepsilon(u\otimes v,w)\otimes 1_z\circ 1_u\otimes\phi\otimes 1_w\otimes 
1_z\circ 1_u\otimes\varepsilon(z,w)\otimes k\otimes 1_w=$$
$$1_w\otimes 1_u\otimes\phi\otimes 1_z\circ\varepsilon(u,w)\otimes
1_z\circ
1_u\otimes\varepsilon(z,w)\circ k\otimes 1_w=$$
$$1_w\otimes 1_u\otimes\phi\otimes 1_z\circ\varepsilon(u\otimes v,w)k\otimes 1_w=
1_w\otimes(1_u\otimes\phi\otimes 1_z\circ k),$$ 
as required.\smallskip 

 A second example of defining subspaces $K_u\subset H_u$ is the following.\smallskip

Consider an inclusion $\cA\subset\cF$ of $C^*$--algebras and suppose that $\cT$ is 
a tensor $C^*$--category of endomorphisms with conjugates of $\cA$, where 
each object $\rho$ of the category is induced by a Hilbert space $H_\rho$ 
of support $I$ 
in $\cF$. We suppose that $\cA'\cap\cF=\Bbb C$ then this Hilbert space is 
unique. We thus have an embedding of $\cT$ in a category of Hilbert spaces 
in $\cF$.\smallskip 

  Now suppose that $\cB$ is a $C^*$--algebra with $\cA\subset\cB\subset\cF$ 
and set $K_\rho:=H_\rho\cap\cB$. The $K_\rho$ are Hilbert spaces in $\cB$. 
In this case $(5.1)$ and $(5.2)$ are obvious. Pick a solution $R,\bar R$ 
of the conjugate equations of $\rho$. Note that 
$j_\rho(\psi):={\bar\rho}(\psi)^* R\in H_{\bar\rho}$, where $\psi\in
H_\rho$. The inverse
$j_{\bar\rho}$ of $j_\rho$ is 
defined by $j_{\bar\rho}\bar\psi:=\rho(\bar\psi)^*\bar R$. 
It is easily checked that these operators are inverses of one another. 
If $\psi\in K_\rho$ then $j_\rho\psi={\bar\rho}(\psi)^* R\in\cB$
and $(5.7)$
follows. If $\phi\in K_\sigma$ and $\chi\in K_{\rho\otimes\sigma}$, then 
$\rho(\phi^*)\chi F=\rho(\phi^*)\rho\sigma(F)\chi=\rho(F)\chi$. Thus 
$\rho(\phi^*)\chi\in H_\rho\cap\cB=K_\rho$, verifying
$(5.3)$.\smallskip
 
  Note that we have another functor $K$ from $\cT$ to the category of Hilbert 
spaces, where $K_T:=T\rest K_\rho$ for $T\in(\rho,\sigma)$. It is not clear 
whether $K$ satisfies $(5.8)$. The problem is that we do not know whether 
$\rho(\phi)\in\cB$. On the other hand we do have 
$$\psi^*K_{\rho\sigma}\subset K_\sigma,$$ 
since if $\chi\in K_{\rho\sigma}$ then $\psi^*\chi$ induces $\sigma$ and 
is thus in $H_\sigma$. But since $\psi$ and $\chi$ are in $\cB$ their 
product is actually in $K_\sigma$.\smallskip

\end{section}

\begin{section}{A Galois correspondence}

We set up a correspondence between functors $K:u\mapsto K_u$ satisfying 
$(5.1)$, $(5.2)$, $(5.7)$ and $(5.8)$  above and closed
bi-ideals $\cJ$. Given $K$ we define 
$K^\perp$ to be the closed ideal generated by the $(u_{ij}-\delta_{ij}I)$ 
for all $i$ and all $j\leq m_u$ and all corepresentations $u$, where $m_u$
is the
Hilbert space dimension of $K_u$, and we have chosen an orthonormal basis
of $H_u$ whose first $m_u$ vectors are an orthonormal basis of $K_u$.
 Obviously, if 
$K_1\subset K_2$, then $K_1^\perp\subset K_2^\perp$. Given $\cJ$, we 
define the functor $\cJ^\perp$ by $\varphi\in\cJ_u^\perp$ 
if and only if $(u-\iota(u))\varphi\in H_u\otimes\cJ$. Again, if 
$\cJ_1\subset\cJ_2$, then $\cJ_{1,u}^\perp\subset\cJ_{2,u}^\perp$. Given $K$ 
and $u$, if $j\leq m_u$ then
$$(u-\iota(u))\varphi_j=\sum_i\varphi_i\otimes(u_{ij}-\delta_{ij}I)\in
 H_u\otimes K^\perp.$$
Thus $K_u\subset K_u^{\perp\perp}$. Given $\cJ$, then if 
$J\in\cJ^{\perp\perp}$, 
$J$ is in the closed ideal generated by the $(u_{ij}-\delta_{ij}I)$ for 
all $i$, $j\leq m_u$ and all $u$. But if $j\leq m_u$ then 
$u\varphi_j=\sum_i\varphi_i\otimes(u_{ij}-\delta_{ij}I)\in H_u\otimes\cJ$, 
so $(u_{ij}-\delta_{ij}I)\in\cJ$ and $\cJ^{\perp\perp}\subset\cJ$. From 
these relations, we can derive the usual closure relations: 
$K^\perp=K^{\perp\perp\perp}$ and $\cJ^\perp=\cJ^{\perp\perp\perp}$.\smallskip 

Up to this point we have made no use of relations $(5.1)$, $(5.2)$,
$(5.7)$ and $(5.8)$ but we now want to show that that 
$K_u=K^{\perp\perp}_u$ and, by Theorem 5.5, it follows that the $K_u$ are
the spaces of invariant vectors for a unique maximal compact quantum
subgroup $G'$ of $G$, where $G'=({\cal A}/\cJ,\Delta')$. We let $\pi:{\cal
A}\to{\cal A}/\cJ$ denote the canonical surjection. Now, by definition, 
$\varphi\in {\cJ}^\perp_u$ if $(u-\iota(u))\varphi\in H_u\otimes \cJ$.
But then $1\otimes\pi u(\varphi)=\varphi\otimes I'$ so $\varphi\in K_u$.
Conversely, if $\varphi\in K_u$, then $1\otimes\pi
u\varphi=1\otimes\pi\iota(u)\varphi$. But
$\text{ker}(1\otimes\pi)=H_u\otimes\cJ$ so $\varphi\in {\cJ}^\perp_u$.
Thus $K_u={\cJ}^\perp_u$ and $K^{\perp\perp}_u={\cJ}^{\perp\perp\perp}_u=
{\cJ}^\perp_u=K_u$ as required.
\end{section}

\begin{section}{Functors arising from ergodic actions}

Consider a nondegenerate action $\eta:{\cal B}\to{\cal B}\otimes{\cal A}$
of a compact quantum group $G=({\cal A},\Delta)$ on a unital $C^*$--algebra ${\cal B}$.
If $\eta$
is  {\it ergodic} (${\cal B}^\eta={\mathbb C}I$),  much can be said about
the spectrum of the action.

Boca showed that any irreducible representation in the spectrum of $\eta$
has a finite multiplicity
$\text{mult}(u)$ bounded above
by the {\it
quantum
dimension}  q--dim$(u)$ of $u$ \cite{Boca}. 

In \cite{BRV} it was shown that 
$\text{mult}(u)$
can be bigger than the Hilbert space dimension of $u$. 
We shall outline, with a different notation, the main tools used in
\cite{BRV}.

For any unitary 
 (not necessarily
irreducible)
representation $u$ of $G$ of finite dimension,
  consider the  vector space $L_u$ of linear maps 
$T: H_u\to{\cal B}$ intertwining $u$ with $\eta$.
 For $S, T\in L_u$,
 $\sum T(e_i)S(e_i)^*$
is independent of the choice of  
orthonormal basis  $(e_i)$ of $H_u$, and defines an element of the fixed
point algebra, which, though, reduces to the complex numbers.
Thus $L_u$ is a Hilbert space with inner product:
$$(S, T)I:=\sum T(e_i)S(e_i)^*.$$
Note that if we regard $S$ and $T$ as elements of 
$(u\otimes\eta,\iota\otimes\eta)$, then 
$(S,T)I=T\circ S^*$.\medskip

\noindent{\bf 7.1 Proposition} {\sl
The Hilbert space $L_u$ is non-trivial if and only if
 $u$ contains a subrepresentation $v\in sp(\eta)$. In particular, 
if $u$ is irreducible then $u\in sp(\eta)$ if and only if $L_u\neq0$.}
\medskip

For any intertwiner $A\in(u,v)$ let $L_A: L_v\to L_u$ be the linear map
defined by:
$$L_A(T)=T\circ A.$$
It is easy to check that $L$ is  a contravariant $^*$--functor.

We obtain a covariant functor passing to the dual space.
Consider, for any $u\in sp(\eta)$ another Hilbert space, 
 $\ov{L}_u$,
which, as a vector
space, is the
complex conjugate
of $L_u$,  with  inner product:
$$(\overline{S},\overline{T}):=\sum S(e_i)T(e_i)^*=(T,S)$$
where $(e_i)$ is an orthonormal basis  of $H_u$.
Identify the complex conjugate of $L_u$ with the dual of $L_u$, and set,
for $A\in (u,v)$,
$$\ov{L}_A: \phi\in\ov{L}_u\to\phi\circ L_A\in\ov{L}_v.$$
This is now a covariant $^*$--functor.
The functors $L$ and $\ov{L}$ are related by:
$$\ov{L}_A(\overline{T})=\overline{L_{A^*}(T)},\quad A\in(v,u), T\in
L_v.$$
\medskip

\noindent{\bf 7.2 Definition} The functor $\ov{L}$ will be called the
{\it spectral functor} associated with the ergodic action $\eta:{\cal
B}\to{\cal
B}\otimes{\cal A}$, 
and $\ov{L}_u$ the {\it spectral subspace} corresponding to the
representation
 $u$.\medskip

\noindent{\bf 7.3 Theorem} {\sl The $^*$--functor $\ov{L}:
\text{Rep}(G)\to{\cal H}$ is
quasitensor. Therefore for any $u\in\text{Rep}(G)$, $\ov{L}_u$ is
finite
dimensional and
$\text{dim}(\ov{L}_{u})=\text{dim}(\ov{L}_{\overline{u}})$.}\medskip

\noindent{\bf Proof} The space $L_\iota$ is, by definition,   the
fixed point algebra,
which reduces to the complex numbers, so $L_\iota={\mathbb
C}=\ov{L}_\iota$, and $(3.1)$ follows. If $S\in L_u$, $T\in L_v$ then
the map $S\otimes T: \psi\otimes\phi\in H_u\otimes H_v\to S(\psi)T(\phi)$
intertwines
$u\otimes v$ with $\eta$, so $S\otimes T\in L_{u\otimes v}$.
If $S'\in L_u$, $T'\in L_v$ is another pair, and $(e_i)$ and $(f_j)$ are
orthonormal bases of $H_u$ and $H_v$ respectively, then
$$(S'\otimes T', S\otimes T)I=\sum_{i,j}  S(e_i)T(f_j)(S'(e_i)T'(f_j))^*=
(S',S)(T',T)I,$$
so the linear span of all $S\otimes T$ is just a copy of $L_u\otimes L_v$
sitting inside $L_{u\otimes v}$. 

Moreover, if $A\in(u',u)$, $B\in(v',v)$,
$S\in L_{u}$, $T\in L_{v}$, $e\in H_{u'}$, $f\in H_{v'}$
$$L_{A\otimes B}(S\otimes T)(e\otimes f)=(S\otimes T)\circ(A\otimes
B)(e\otimes f)=$$
$$S\otimes T(Ae\otimes
Bf)=S(Ae)T(Bf)=$$
$$L_A(S)(e)L_B(T)(f)=L_A(S)\otimes L_B(T)(e\otimes f).$$
Hence 
$$L_{A\otimes B}\upharpoonright_{L_u\otimes L_v}=L_A\otimes L_B.$$
Next we define an  inclusion 
$$S_{u,v}: \ov{L}_u\otimes \ov{L}_v\to\ov{L}_{u\otimes v}$$
by taking $\overline{S}\otimes\overline{T}$ to $\overline{S\otimes T}$,
for $S\in L_u$, $T\in L_v$ and hence $S\otimes T\in L_{u\otimes v}$.
As for $L$, one can check that this inclusion is an isometry from the
tensor product Hilbert space $\ov{L}_u\otimes\ov{L}_v$ to the
Hilbert space
$\ov{L}_{u\otimes v}$, and $(3.2)$ follows.
The relation between $L$ and $\ov{L}$ on arrows shows that 
$$\ov{L}_{A\otimes
B}\upharpoonright_{\ov{L}_{u'}\otimes
\ov{L}_{v'}}=\ov{L}_A\otimes\ov{L}_B,$$
and this shows $(3.6)$.
Relation $(3.3)$ is obvious for $L$ and $\ov{L}$ as well. Equation 
$(3.4)$ follows
from the fact that if $R\in L_u$, $S\in L_v$, $T\in L_z$ then $(R\otimes
S)\otimes T=R\otimes(S\otimes T)$.
It remains to show property $(3.5)$.
We define, for any $S\in L_{v}$, $T\in L_{u\otimes v}$,
 a linear map $\tilde{S}(T): H_u\to{\cal B}$ by
$$\tilde{S}(T)(\psi)=\sum_k T(\psi\otimes f_k)S(f_k)^*,$$
which is independent of the choice of the orthonormal basis $(f_k)$ of
$H_v$. We check that $\tilde{S}(T)\in L_u$. In fact, on an orthonormal
basis  $(\psi_i)$ of $H_u$,
$$\eta\circ\tilde{S}(T)(\psi_i)=\eta(\sum_k T(\psi_i\otimes
f_k)S(f_k)^*)=$$
$$\sum_k\eta(T(\psi_i\otimes f_k))\eta(S(f_k))^*=\sum_{r,s,k,p}
(T(\psi_r\otimes f_s)\otimes u_{ri}v_{sk})(S(f_p)\otimes v_{pk})^*=$$
$$\sum_{r,s,k,p} 
T(\psi_r\otimes f_s)S(f_p)^*\otimes
u_{ri}v_{sk}v_{pk}^*=
\sum_{r,s}
T(\psi_r\otimes f_s)S(f_s)^*\otimes
u_{ri}=$$
$$\sum_r\tilde{S}(T)(\psi_r)\otimes u_{ri}.$$
We have thus defined a linear map:
$$\tilde{S}: L_{u\otimes v}\to L_u.\eqno(7.1)$$
Conjugating $\tilde{S}$ gives 
a linear map
$$\hat{{S}}:\overline{T}\in\ov{L}_{u\otimes
v}\to\overline{\tilde{S}(T)}\in\ov{L}_u.$$
Consider the operator 
$$r(\overline{S}): \ov{L}_u\to\ov{L}_{u}\otimes\ov{L}_v$$
which tensors on the right by $\overline{S}$ and the previously
defined isometric inclusion
map
$$S_{u,v}:\ov{L}_u\otimes\ov{L}_v\to\ov{L}_{u\otimes v}.$$
We claim that
$$\hat{{S}}^*=S_{u,v}\circ r(\overline{S}).$$
In fact, for $T\in L_{u\otimes v}$, $T'\in L_u$,
$$(\overline{T'}, \hat{S}^*\overline{T})=(\hat{S}\overline{T'},
\overline{T})=$$
$$(\overline{\tilde{S}(T')},\overline{T})=\sum_{i,j}
T'(\psi_i\otimes f_j)S(f_j)^*T(\psi_i)^*=$$
$$(\overline{T'},
\overline{T\otimes S})=(\overline{T'},
S_{u,v}\overline{T}\otimes\overline{S})=$$
$$(\overline{T'},
S_{u,v}\circ r(\overline{S})\overline{T}),$$
and the claim is proved. Taking the adjoint,
$$\hat{{S}}=r(\overline{S})^*\circ S_{u,v}^*,$$
which implies 
$$S_{u,v}^*=\sum_p r(\overline{S_p})\hat{{S_p}},$$
with $(\overline{S_p})$ an orthonormal basis of $\ov{L}_v$.
Replace $u$ with a tensor product representation $u\otimes z$,
and compute, for $\overline{T}\in\ov{L}_u$,
$\overline{T'}\in\ov{L}_{z\otimes v}$
$$S_{u\otimes z, v}^*\circ
S_{u,z\otimes v}(\overline{T}\otimes\overline{T'})=\sum_p
r(\overline{S_p})\circ\hat{S_p}(\overline{T\otimes T'})=$$
$$\sum_p r(\overline{S_p})(\overline{\tilde{S_p}(T\otimes T')}).$$
Now for $\psi\in H_u$, $\phi\in H_z$,
$$\tilde{S_p}(T\otimes T')(\psi\otimes\phi)=\sum_k (T\otimes
T')(\psi\otimes\phi\otimes f_k)S_p(f_k)^*=$$
$$ T(\psi)\sum_k T'(\phi\otimes
f_k)S_p(f_k)^*=T\otimes\tilde{S_p}(T')(\psi\otimes\phi)$$ so
$$E_{u\otimes z, v}\circ S_{u, z\otimes
v}(\overline{T}\otimes\overline{T'})=\sum_p\overline{T\otimes\tilde{S_p}(T')\otimes
S_p}$$
and therefore
 $$E_{u\otimes z,
v}(\text{Image}S_{u, z\otimes v})\subset
\text{Image}S_{u,z,v}$$
and the proof is complete.
\medskip

Let $\overline{u}$ be a conjugate of $u$ defined by $R\in(\iota,
\overline{u}\otimes u)$ and $\overline{R}\in(\iota, u\otimes\overline{u})$
solving the conjugate equations.  We shall identify
the solution $\hat{R}$, $\hat{\overline{R}}$ of the conjugate equations
for $\ov{L}_u$ constructed in Theorem 3.7 and we shall relate it to the 
notion of {\it quantum multiplicity} introduced in \cite{BRV}.

We can write
$R=\sum_i j e_i\otimes e_i$ and $\overline{R}=\sum _k j^{-1}f_k\otimes
f_k$, with $j: H_u\to H_{\overline{u}}$ an antilinear invertible
map and $(e_i)$ and $(f_k)$ orthonormal bases of $H_u$ and
$H_{\overline{u}}$, respectively.

For $T\in L_u$, the multiplet $T(e_i)^*$ transforms like the complex
conjugate  invertible representation $u_*$, so the linear
map
$$\tilde{T}:=\overline{\psi}\in\overline{H_u}\to T(\psi)^*\in{\cal
B},\quad\psi\in
H_u$$
intertwines $u_*$ with $\eta$.

Consider the
linear map $$Q:
H_{u_*}=\overline{H_u}\to H_{\overline{u}}$$
obtained composing the (antiunitary) complex conjugation
$\overline{H_u}\to
H_u$ with ${j^*}^{-1}$. Regard $\overline{H_u}$ as the Hilbert space for
$u_{*}$.
The condition that $R$ is an intertwiner can be written:
$Q\in (u_*,\overline{u})$.

Thus $\tilde{T}\circ Q^{-1}: H_{\overline{u}}\to {\cal C}$
intertwines $\overline{u}$ with $\eta$.
We have therefore  defined an 
antilinear invertible  map
$$J: L_u\to L_{\overline{u}},$$
$$J(T)(\psi)=\tilde{T}\circ Q^{-1}(\psi)=T(j^{*}(\psi))^*,$$
with inverse
$J^{-1}:L_{\overline{u}}\to L_u,$
$$J^{-1}(S)(\phi)=S({j^{*}}^{-1}(\phi))^*, \quad \phi\in H_u.$$
If we conjugate $J$ with the antilinear invertible
map
$T\in L_u\to\overline{T}\in\ov{L}_u$, we obtain 
an antilinear invertible
$$\ov{J}:\ov{L}_u\to\ov{L}_{\overline{u}}.$$
If $u$ is irreducible, the  {\it quantum multiplicity} of
$u$ is defined in \cite{BRV} by
$$\text{q--mult}(u)^2:=
\text{Trace}(\ov{J}\ov{J}^*)\text{Trace}((\ov{J}\ov{J}^*)^{-1}).$$
\medskip

\noindent{\bf 7.4 Theorem} {\sl If $u$ is irreducible, one has
$$\hat{R}=\sum \ov{J} \overline{T_i}\otimes\overline{T_i},$$
$$\hat{\overline{R}}=\sum\ov{J}^{-1}\overline{S_j}\otimes\overline{S_j},$$
where $(\overline{T_i})$ and $(\overline{S_j})$ are orthonormal bases of
$\ov{L}_u$ and $\ov{L}_{\overline{u}}$, respectively. 
Therefore 
$$\text{q--mult}(u)=d_{\hat{R},\hat{\overline{R}}}(\ov{L}_u).$$}\medskip

\noindent{\bf Proof} We need to show that for $\overline{T}\in
\ov{L}_u$,
$r(\overline{T})^*\circ\hat{R}=
\ov{J}\overline{T}$. Recall that
$\hat{R}=S_{\overline{u},u}^*\ov{L}_R$, so, with the same notation as
in
the previous theorem,
$$r(\overline{T})^*\circ\hat{R}=
r(\overline{T})^*\circ S_{\overline{u}, u}^*\circ\ov{L}_R=
\hat{T}\circ\ov{L}_R=\hat{T}(\ov{L}_R).$$
In the last equation we have identified 
$\ov{L}_R\in(\ov{L}_\iota,\ov{L}_{\overline{u}\otimes u})$
with the element
$\ov{L}_R(\overline{1})$,
 with $1\in {L}_\iota$ the intertwiner $1\in
H_\iota\to
I\in{\cal B}$.
Now
$$\ov{L}_R(\overline{1})=\overline{L_{R^*}(1)}=\overline{1\circ
R^*},$$
so
$$\hat{T}(\ov{L}_R)=\overline{\tilde{T}(1\circ R^*)}.$$
But $1\circ R^*\in L_{\overline{u}\otimes u}$ takes $\xi\in
H_{\overline{u}}\otimes H_u$ to $(R,\xi)I$, hence 
$\tilde{T}(1\circ R^*)\in L_{\overline{u}}$ is 
the map taking $\psi\in H_{\overline{u}}$ to
$$\sum_i (1\circ R^*)(\psi\otimes e_i)T(e_i)^*=\sum_i (R, \psi\otimes
e_i)T(e_i)^*=$$
$$\sum_i(je_i,\psi)T(e_i)^*=\sum_i(j^*\psi,
e_i)T(e_i)^*=T(j^*\psi)^*=J(T)(\psi).$$
So $\tilde{T}(1\circ R^*)=J(T)$, and the proof is complete.
\medskip

As a consequence, we obtain the following result, proved in \cite{BRV}, in
turn
extending  Boca's result.
\medskip

\noindent{\bf 7.5 Corollary} {\sl Let $\eta$ be an ergodic
nondegenerate 
action of a compact quantum group $G$ on a unital $C^*$--algebra ${\cal
B}$. Then for any irreducible representation $u$ of $G$, 
$$\text{mult}(u)\leq q\text{--mult}(u)\leq q\text{--dim}(u).$$
Furthermore $q\text{--mult}(u)=q\text{--dim}(u)$ if and only if 
$\ov{L}_R\in\text{Image}S_{\overline{u}, u}$ and
$\ov{L}_{\overline{R}}\in\text{Image}S_{u,\overline{u}}$
for some (and hence any) solution $(R, \overline{R})$ to the conjugate
equations for $u$.
}\medskip

\noindent{\bf Proof} We are stating, thanks to the previous theorem,  
that
$$\text{dim}(\ov{L}_u)\leq
d_{\hat{R},\hat{\overline{R}}}(\ov{L}_u)\leq 
d_{R,\overline{R}}(u),$$
and this follows from Corollary 3.8  applied to the functor $\ov{L}$.
\medskip

In the last part of this section we relate
the functor $\ov{L}$ associated to the $G$--action on a quantum
quotient space with the functor $K$ of invariant vectors.
\medskip

\noindent{\bf 7.6 Proposition} {\sl Let $K$ be a compact quantum subgroup
of a compact quantum group $G$.
A unitary  $G$--representation $u$ contains a  subrepresentation of
$sp({\eta}_K)$ if and only if 
$K_u\neq0$. 
In particular, if $u$ is irreducible,
 then $u\in
sp({\eta}_K)$ if and only if $K_u\neq0$. In this case,
$\text{mult}(u)=\text{dim}(K_u)$.
}\medskip 

\noindent{\bf Proof} 
Proposition 2.6 shows that any unitary irreducible 
$G$--representation $v$ for which $K_v\neq0$  lies in
$sp({\eta}_K)$. If $u$ is any representation such that $K_u\neq0$ then
$u$ contains an irreducible subrepresentation $v$ such that $K_v\neq0$.
So $v\in sp(\eta_K)$. 
Conversely, if $u$ contains a spectral subrepresentation $u'$ and
if $z$ is an irreducible component of $u'$, then by Prop 2.1, $z\in
sp(\eta_K)$. By Prop. 2.6 $K_z\neq 0$, so  $K_u\neq0$.
\medskip

For any $\psi\in K_u$, we have a map $T_\psi: H_u\to{\cal
A}^\delta_{sp}$ defined by
$$T_\psi(\phi)=\ell_{\psi}^*u(\phi),$$
which actually lies in $L_u$, as on an orthonormal basis of $H_u$:
$$\eta_K\circ T_\psi(e_i)=\Delta(\ell_\psi^*u(e_i))=$$
$$\sum_k
\ell_\psi^*u(e_k)\otimes \ell_{e_k}^*u(e_i)=T_\psi\otimes\iota\circ
u(e_i).$$
\medskip

\noindent{\bf 7.7 Theorem} {\sl Let $K$ be a compact quantum subgroup of
$G$.
The map
$$U_u:\psi\in K_u\to\overline{T_\psi}\in\ov{L}_u$$
is a  quasitensor natural transformation from the functor $u\to
K_u$ to the
functor
$u\to\ov{L}_u$ such that $U_u$ is a unitary whenever $K_u\neq0$.}\medskip

\noindent{\bf Proof} 
We first show naturality, namely that for $A\in(u,v)$,
$$\ov{L}_A\circ U_u=U_v\circ K_{A}.$$
For $\psi\in K_u$ the right hand side computed on $\psi$ equals
$\overline{T_{A\psi}}$, and the left hand side equals
$$\ov{L}_A(\overline{T_\psi})=\overline{L_{A^*}(T_\psi)}=
\overline{T_\psi\circ A^*}.$$
Now for $\phi\in H_v$,
$$T_\psi\circ
A^*(\phi)=\ell_\psi^*u(A^*\phi)=$$
$$\ell_{\psi}^*A^*\otimes\iota
v(\phi)=\ell_{A\psi}^* v(\phi)=T_{A\psi}(\phi),$$
so $$\overline{T_\psi\circ A^*}=\overline{T_{A\psi}}.$$
We show that $U_u$ is isometric:
$$(\overline{T_\psi}, \overline{T_\phi})I=\sum_k
T_\psi(e_k)T_\phi(e_k)^*=$$
$$\sum_k\ell_\psi^*u(e_k)((\ell_\phi)^*u(e_k))^*=
(\psi,\phi) I$$
by the unitarity of $u$.
Now, if $u$ is irreducible, then $U_u$ is known to be surjective.
More generally, if $v_i$ are irreducible subrepresentations of $u$ and
$S_i\in(v_i, u)$ are isometries such that $\sum_i S_iS_i^*=I$
then by naturality
$U_{u}K_{S_i}=\ov{L}_{S_i}\circ U_{v_i},$
which shows that $U_uU_u^*=I$.
Finally, restricting $U_{u\otimes v}$ to $K_u\otimes K_v$ gives,
for $\psi\in K_u$, $\psi'\in K_v$,
$$U_{u\otimes v}(\psi\otimes\psi')=\overline{T_{\psi\otimes\psi'}}.$$
On the other hand for $\phi\in H_u$, $\phi'\in H_v$,
$$T_{\psi\otimes\psi'}(\phi\otimes\phi')=\ell_{\psi\otimes\psi'}^*(u\otimes
v)(\phi\otimes\phi')=$$
$$\ell_\psi^*u(\phi)\ell_{\psi'}^*v(\phi')=(T_\psi\otimes
T_{\psi'})(\phi\otimes\phi')$$
so
$$U_{u\otimes v}(\psi\otimes\psi')=\overline{T_\psi\otimes
T_{\psi'}}=$$
$$\overline{T_{\psi}}\otimes\overline{T_{\psi'}}=U_u(\psi)\otimes
U_v(\psi').$$
 \medskip

\end{section}

\begin{section} {Multiplicity maps in ergodic  actions}

Let $\eta: {\cal B}\to{\cal B}\otimes{\cal A}$ be a nondegenerate ergodic
action of a
compact quantum group $G=({\cal A},\Delta)$ on a unital $C^*$--algebra
${\cal B}$.  
Consider the Hilbert space $\ov{L}_u$ associated with the
unitary representation $u$, and 
an orthonormal basis $(\overline{T_k})$ in $\ov{L}_u$. Let
  $$c_u: H_u\to\ov{L}_u\otimes{\cal B}$$ denote  the linear map
defined by 
$$c_u(\phi)=\sum_k \overline{T_k}\otimes T_k(\phi).$$
This map does not depend on the choice of the orthonormal basis of
$\ov{L}_u$.
We shall call the map $c_u$ the {\it multiplicity map} of $u$ in
$\eta$.

\medskip

\noindent{\bf 8.1 Proposition} {\sl The map $c_u$ is nonzero if and only
if
$L_u$ is nonzero, i.e. if and only if $u$ contains a spectral
subrepresentation.}\medskip

\noindent{\bf Proof}
We need to show that $c_u=0$ implies $L_u=0$.
Indeed, if  $L_u$ were $\neq0$ then 
for any orthonormal basis $(\overline{T}_k)$ of $\ov{L}_u$, $T_k=0$.
Since $(T_k)$ is a linear basis of $L_u$, we must have $L_u=0$. A
contradiction.\medskip

We can also represent
$c_u$ as the rectangular matrix, still denoted by $c_u$,
whose $k$--th row is $(T_k(e_1),\dots, T_k(e_d))$, with
$(e_i)$ an orthonormal basis of $H_u$,
$d=\text{dim}(u)$, $k=1,\dots,\text{dim}(L_u)$.
\medskip

\noindent{\bf 8.2 Proposition} {\sl Let $u$ contain a 
subrepresentation in the spectrum of $\eta$.
Then the matrix $c_u$ satisfies the following properties.
\begin{description}
\item{\rm a)}
$c_uc_u^*=I$,
\item{\rm b)} each row of $c_u$  transforms like $u$,
\item{\rm c)} for any   multiplet $c$ transforming
like $u$, $cP_u=c$, with $P_u$ the domain projection of $c_u$.
\end{description}
Conversely, if  $c'_u\in M_{p, \text{dim}(u)}({\cal B})$ 
satisfies a)--c) then 
$p=\text{dim}(L_u)$ and $c'_u$ is of the form of $c_u$.
}\medskip

\noindent{\bf Proof} Property 
b) is obviously satisfied, a) and c) are a consequence of the
fact that $(\overline{T_k})$
is an orthonormal basis of $\ov{L}_u$. We show uniqueness. If 
$c'_u$ is as in the statement then
the rows of $c'_u$
lie in $L_u$ by b). Furthermore these rows must be an orthonormal basis
of $L_u$ by $a)$ and $c)$. In particular,
$p=\text{dim}(L_u)$.
\medskip

We next investigate the relationship between $u\to c_u$ and 
the functor $u\to\ov{L}_u$.
Define general coefficients of $c_u$: for $\phi\in\ov{L}_u$,
$\psi\in H_u$, set
$$c^u_{\phi,\psi}:=\ell_{\phi}^*c_u(\psi)\in{\cal B}.$$
\medskip

\noindent{\bf 8.3 Proposition} {\sl The map $c_u:
H_u\to\ov{L}_u\otimes{\cal B}$ satisfies the following properties.
\begin{description}
\item{\rm a)}
$\ov{L}_A\otimes I\circ c_u=c_v\circ A,\quad    A\in(u, v),$
\item{\rm b)} for $\phi\in\ov{L}_u$, $\phi'\in\ov{L}_v$, $\psi\in
H_u$, $\psi'\in H_v$,
$$c^{u\otimes v}_{\phi\otimes\phi',
\psi\otimes\psi'}=c^u_{\phi,\psi}c^u_{\phi',\psi'}.$$
\end{description}}
\medskip

Let $K$ be a compact quantum subgroup of $G$ and $\eta_K: {\cal
A}^\delta\to{\cal A}^\delta\otimes{\cal A}$ the ergodic action defining
the 
quantum quotient space $K\backslash G$. We can then associate to any
representation $u$ of $G$ the multiplicity map $c^K_u$.

In  the last part of Sect. 4 we have also defined, for any
representation $u$,
the map $u^K=E_u^K\otimes I\circ u: H_u\to K_u\odot{\cal
A}^\delta_{sp}$. 
Now consider the natural unitary transformation
$U_u: K_u\to\ov{L}_u$ between the
functors $u\to K_u$ and $u\to\ov{L}_u$ defined in the previous
section.
\medskip

\noindent{\bf 8.4 Proposition} {\sl If $U:K\to{\ov{L}}$ is the unitary
natural transformation defined in Theorem 7.7, then for any representation
$u$,
$$U_u\otimes\iota\circ u^K=c^K_u.$$}\medskip

\noindent{\bf Proof} 
If $(\psi_k)$ is an orthonormal basis of $K_u$, we can write
$u^K=\sum_k\psi_k\otimes T_{\psi_k}$. The rest of the proof is now clear.
\medskip

\noindent{\bf 8.5 Corollary} {\sl If $(\psi_k)$ and $(e_j)$ are
orthonormal 
bases of $K_u$ and $H_u$ respectively, the matrix
$(u^K_{\psi_k,e_j})=(\ell_{\psi_k}^*u(e_j))$ satisfies properties a)--c)
of Prop.\ 8.2.}
\medskip

We next show the linear independence of the coefficients of the $c_u$'s for
inequivalent irreducibles.
The proof is obtained by generalizing a result of Woronowicz \cite{Wcmp}
stating the linear independence of matrix coefficients of inequivalent
irreducible representations of a compact matrix pseudogroup.
\medskip

\noindent{\bf 8.6 Proposition} {\sl 
Let $\eta:{\cal C}\to{\cal C}\odot {\cal A}_\infty$ be an action of the
Hopf $^*$--algebra $G_\infty$
on a unital $^*$--algebra ${\cal C}$, and let 
$S$ be a set of unitary, irreducible, pairwise inequivalent, 
representations of $G$ in the spectrum of $\eta$.
For each $u\in S$, let 
$c_u=(c^u_{ij})\in M_{p_u,\text{dim}(u)}({\cal C})$
satisfy a) and b) of Prop.\ 8.2. 
It follows that the set of matrix coefficients $\{c^u_{ij}, 
i=1,\dots, p_u,j=1,\dots,\text{dim}(u), u\in S\}$ 
is linearly independent in ${\cal C}$.}\medskip

\noindent{\bf Proof} 
Let $F$ be a finite subset of $S$.
Consider
the linear subspace 
of $\oplus_{u\in F}M_{p_u, \text{dim}(u)}$,
$$M:=\{\oplus_{u\in F} c^u_\rho:=\oplus_{u\in F}(\rho(c^u_{ij})),
\rho\in{\cal C}'\},$$
with ${\cal C}'$ the dual of ${\cal C}$ as a vector space 
and also the subspace of $\oplus_{u\in F}M_{\text{dim}(u)}$,
$$B:=\{\oplus_{u\in
F}u_\sigma:=\oplus_{u\in F}(\sigma(u_{pq})),\sigma\in{\cal
A}_\infty'\}.$$
Notice that $MB\subset M$, since for $\rho\in{\cal C}'$,
$\sigma\in{\cal A}_\infty'$, $u\in S$, 
$c^u_\rho u_\sigma=c^u_{\rho*\sigma}$, where
$\rho*\sigma:=\rho\otimes\sigma\circ\eta$.
On the other hand,
by Lemma 4.8 in
\cite{Wcmp},  $B=\oplus_{u\in F}M_{d_u}$.
Thus,
for each fixed $u^0\in F$,
choosing $\sigma$ such that
$u_{\sigma}$ is zero for $u\neq u^0$ and 
$u^0_{\sigma}$ is a matrix unit
$\theta^{u^0}_{h,k}$, shows that 
for every $\rho\in{\cal C}'$, every $u^0\in F$ and every
$h,k=1,\dots, \text{dim}({u^0})$ there exists $\rho'\in{\cal C}'$
such that $c^u_{\rho'}$ is zero if $u\neq u^0$ and 
$c^{u^0}_{\rho'}$ is the matrix with
identically zero columns except
for the $k$-th one, which coincides with the $h$-th column of
$c^{u^0}_\rho$.
Assume now that we have a vanishing finite linear combination in ${\cal
C}$:
$\sum_{i,j,u} \lambda^u_{ij}c^u_{ij}=0$. 
Then for all $\rho'\in{\cal
C}'$,
$\sum \lambda^u_{ij}\rho'(c^u_{ij})=0$. 
Choose $F$ finite and large enough so that this sum runs over $F$.
By the previous
arguments, for every $\rho\in{\cal C}'$, every $u^0\in S$ and every choice
of $h,k$,
$\sum_{{i=1}}^m\lambda^{u^0}_{i,k}\rho(c^{u^0}_{i,h})=0$. It follows that
$$\sum_i\lambda^{u^0}_{i,k}c^{u^0}_{i,h}=0,\quad h,k=1,\dots,d_{u^0}.$$
Multiplying on the right by ${c^{u^0}}_{p,h}^*$, summing up over $h$ and
using
$c_{u^0}c_{u^0}^*=I$ gives $\lambda^{u^0}_{p,k}=0$ for all
$u^0\in S$, $p=1,\dots, p_{u^0}$,
$k=1,\dots,\text{dim}({u^0})$.
\medskip

\noindent{\bf 8.7. Theorem} {\sl Let $\eta:{\cal B}\to{\cal B}\otimes{\cal
A}$ be a nondegenerate, ergodic  $G$--action of a compact quantum group
on a unital
$C^*$--algebra ${\cal B}$.  Let $\hat{\eta}$
be a complete set of unitary irreducible elements of $sp(\eta)$.
Associate with any $u\in\hat{\eta}$ a corresponding multiplicity map
$c_u$.
Then
the set
of all matrix coefficients  
$$\{\ell_{\overline{T_i}}^*c^u(e_j)=T_i(e_j),
u\in\hat{\eta}, 
(\overline{T_i}) \text{ o.n.b. of } \ov{L}_u , (e_j) \text{ o.n.b. of }
H_u  
\}$$
 is a linear basis for the dense
spectral $^*$--subalgebra ${\cal B}_{sp}$.}\medskip

\end{section}

\begin{section} {A duality theorem for ergodic $C^*$--actions}

We have seen in Sect. 7 that an ergodic nondegenerate action 
$\eta$ of a compact quantum group $G$ on a unital $C^*$--algebra has an
associated quasitensor $^*$--functor $\ov{L}:\text{Rep}(G)\to{\cal H}$
to the category of, necessarily finite dimensional,
Hilbert spaces. In this section we
shall conversely construct from 
a quasitensor $^*$--functor ${\cal F}$
an ergodic nondegenerate action having ${\cal F}$ as its spectral functor.
\medskip

\noindent{\bf 9.1 Theorem} {\sl Let $G=({\cal A},\Delta)$ be a compact
quantum
group and
${\cal F}:\text{Rep}(G)\to{\cal H}$ be a quasitensor $^*$--functor. Then
there exists a 
unital $C^*$--algebra ${\cal B}_{\cal F}$ with an  ergodic nondegenerate
$G$--action   $\eta_{\cal F}:{\cal
B}_{\cal F}\to{\cal B}_{\cal F}\otimes{\cal A}$ 
and a quasitensor natural unitary transformation from the associated
spectral functor $\ov{L}$ to ${\cal F}$.}\medskip

The proof
of the above theorem is
inspired by 
the proof of the Tannaka--Krein duality theorem given by Woronowicz in
\cite{Wtk}. A similar result is obtained in \cite{BRV} for unitary fibre
functors.\medskip

\noindent{\bf Proof}
We start by considering a complete set $\hat{{\cal F}}$ 
of inequivalent, unitary, irreducible representations  $u$ of $G$ such
that ${\cal F}({u})\neq
0$, and form 
the algebraic direct sum
$${\cal C}_{\cal F}=\oplus_{u\in\hat{{\cal F}}}\overline{{\cal
F}(u)}\otimes H_u.$$
We shall make ${\cal C}_{\cal F}$ into a unital $^*$--algebra.
Consider, for each $u\in\hat{{\cal F}}$, orthonormal bases  $(T_k)$ 
and $(e_i)$
of ${\cal F}(u)$ and 
$H_u$ respectively, form the orthonormal basis   
$(\overline{T_k}\otimes e_i)$ of $\overline{{\cal F}(u)}\otimes H_u$.
The linear map
$$c_u: \phi\in H_u\to \sum_k T_k\otimes(\overline{T_k}\otimes\phi)\in
{\cal F}(u)\otimes {\cal C}_{\cal F}$$ is independent of the choice  of
orthonormal basis. 
Extend the definition of $c_u$ to all objects of $\text{Rep}(G)$: first
set $c_u=0$ if $u$ is irreducible but ${\cal F}(u)=0$. Then
 choose isometries $S_i\in(u_i, u)$, with $u_i\in\hat{{\cal F}}$, such
that
$\sum_i S_iS_i^*=I$, and set
$$c_u:=\sum_i{\cal F}({S_i})\otimes I\circ c_{u_i}S_i^*.$$ 
If $S'_j\in(v_j, u)$ is
another set of orthogonal isometries with ranges adding up to the
identity,
$S_i^*S'_j\in(v_j, u_i)$ are always scalars, by irreducibility, 
and nonzero only if $v_j=u_i$, so
$${\cal F}({S_i^*S'_j})\otimes I c_{v_j}=c_{u_i}S_i^*S'_j.$$ Multiplying
on the
left by ${\cal F}({S_i})$, on the right by ${S'}_j^*$, summing over $i$
and $j$ 
and using the fact that ${\cal F}$ is a $^*$--functor, shows that $c_u$
 is independent of the choice of the $S_i$'s, and one now has:
$${\cal F}(A)\otimes I\circ c_u=c_v\circ A,\quad
A\in(u,v),u,v\in\text{Rep}(G).$$
For $u,v\in\hat{{\cal F}}$, $\overline{T}\otimes\phi\in
\overline{{\cal F}(u)}\otimes
H_u$, 
$\overline{T'}\otimes\phi'\in\overline{{\cal F}(v)}\otimes H_v$, set
$$(\overline{T}\otimes\phi)(\overline{T'}\otimes\phi'):=\ell_{T\otimes 
T'}^*\circ c_{u\otimes
v}(\phi\otimes\phi').$$ 
In this way ${\cal C}_{\cal F}$ becomes an
associative algebra with
identity $I=\overline{1}\otimes 1\in \overline{{\cal F}(\iota)}\otimes
H_\iota.$
We next define the $^*$--involution on ${\cal C}_{\cal F}$ with the help
of the
conjugate representation. This representation is
 defined, up
to unitary equivalence, by intertwiners $R\in(\iota, \overline{u}\otimes
u)$ and $\overline{R}\in(\iota, u\otimes\overline{u})$ satisfying the
conjugate equations. If $u$ is irreducible, the spaces $(\iota,
\overline{u}\otimes u)$ and $(\iota, u\otimes\overline{u})$ are one
dimensional. For $u\in\hat{{\cal F}}$
consider orthonormal bases $(T'_i)$, $(e'_p)$ of ${\cal F}(\overline{u})$
and
$H_{\overline{u}}$ respectively, and set, for $T\in {\cal F}(u)$, $\phi\in
H_u$,
$$
(\ell_T^*c_u(\phi))^*
:=\sum_{i,p}({\cal F}(R), S_{\overline{u},u}T'_i\otimes
T)\ell_{T'_i}^*c_{\overline{u}}(e'_p)(\phi\otimes e'_p, \overline{R}).$$
Any other solution to the conjugate equations is of the form $(\lambda R,
\mu\overline{R})$ with $\lambda,\mu\in{\mathbb C}$
and $\overline{\mu}\lambda=1$, so the above definition is independent of
the
choice of $(R,\overline{R}).$ Recall that the pair $(R,\overline{R})$
defines a solution $(\hat{R},
\hat{\overline{R}})$ of the conjugate equations for ${\cal F}(u)$ as in
Theorem 3.7. So writing 
$$R=\sum_i je_i\otimes e_i=\sum_kf_k\otimes j^*f_k,$$
$$\overline{R}=\sum_kj^{-1}f_k\otimes f_k=\sum_ie_i\otimes {j^{-1}}^*
e_i,$$
 with $(e_i)$ and $(f_k)$ o.n.b. of $H_{u}$ and $H_{\overline{u}}$
respectively, and, again, with an analogous meaning of notation,
$$\hat{R}=\sum_i JT_i\otimes T_i=\sum_kT'_k\otimes J^*T'_k,$$
$$\hat{\overline{R}}=\sum_kJ^{-1}T'_k\otimes T'_k=\sum_iT_i\otimes
{J^{-1}}^*
T_i,$$ we obtain for the adoint the equivalent form
$$(\ell_T^*c_u(\phi))^*=\ell_{JT}^*c_{\overline{u}}({j^{-1}}^*\phi),$$
which is obviously independent of the choice of the orthonormal bases.
Replacing $u$ by $\overline{u}$, and therefore $R$ by $\overline{R}$
and $\hat{R}$ by $\hat{\overline{R}}$, and $j$ and $J$ in turn by
their inverses, shows that
$$(\ell_T^*c_u(\phi))^{**}=\ell_T^*c_u(\phi).$$
We have thus defined an antilinear involutive map $^*:{\cal
C}_{\cal F}\to{\cal
C}_{\cal F}$. 

We next show that 
for $u,v\in{\hat{\cal F}}$, $\phi\in H_u, T\in {\cal F}(u)$, $\phi'\in
H_v$, $T'\in{\cal F}(v)$,
$$((\bar T\otimes\phi)(\bar T'\otimes\phi'))^*=
(\bar T'\otimes\phi')^*(\bar T\otimes\phi)^*.$$
Choose
irreducible 
representations $(u_r)$ and 
 isometries $S_r\in (u_r, u\otimes v)$ 
with
pairwise orthogonal 
ranges, adding up to the identity. Then with an obvious meaning of 
notation, the previous equation becomes
$$\sum_r\ell_{J_r{\cal F}(S_r)^*T\otimes
T'}^*c_{\overline{u_r}}({j_r^{-1}}^*S_r^*\phi\otimes\phi')=
\ell_{J_vT'\otimes
J_uT}^* c_{\overline{v}\otimes\overline{u}}({j_v^{-1}}^*\phi'\otimes
{j_u^{-1}}^*\phi).$$
Now recall from \cite{LongoRoberts} that 
in any tensor $C^*$--category with conjugates ${\cal T}$ there is an
antilinear map
$$(\rho, \sigma)\to(\overline{\rho},\overline{\sigma})$$ given by
$$S\to
\overline{S}:=1_{\overline{\sigma}}\otimes{\overline{R}_\rho}^*\circ
1_{\overline{\sigma}}\otimes S^*\otimes 1_{\overline{\rho}}\circ
R_\sigma\otimes1_{\overline{\rho}}.$$
Applying this to categories of Hilbert spaces, with conjugates defined by
antilinear invertible maps $j_\rho$, $j_\sigma$, related in the usual way
to $R_\rho$ and $R_\sigma$, one finds that 
$$\overline{S}=j_{\sigma}\circ S\circ j_{\rho}^{-1}.$$
We claim, see Lemma 9.2, that ${\cal F}(\overline{S})=\overline{{\cal
F}(S)},$
where the latter is defined using the antilinear invertible maps $J_\rho$,
$J_\sigma$.
Also, if $\overline{\sigma}$ is a conjugate of $\sigma$ defined by 
$(R_\sigma, \overline{R}_\sigma)$ and if $\overline{\tau}$ is a conjugate
of $\tau$ defined by $R_\tau$ and $\overline{R}_\tau$ then
$\overline{\tau}\otimes\overline{\sigma}$ is a conjugate  of
$\sigma\otimes\tau$ defined by
$$R_{\sigma\otimes\tau}=1_{\overline{\tau}}\otimes R_\sigma\otimes
1_\tau\circ R_\tau$$
$$\overline{R}_{\sigma\otimes\tau}=1_{\sigma}\otimes
\overline{R}_\tau\otimes
1_{\overline{\sigma}}\circ \overline{R}_\sigma.$$
In conclusion, we obtain an antilinear map
$$S\in(\rho,
\sigma\otimes\tau)\to\overline{S}\in(\overline{\rho},
\overline{\tau}\otimes\overline{\sigma}).$$
Thus 
for an intertwiner $S$ from the
Hilbert
space $\rho$ to the Hilbert space $\sigma\otimes\tau$, 
$\overline{S}$ is the operator from $\overline{\rho}$ to
$\overline{\tau}\otimes\overline{\sigma}$ given by
$$\overline{S}=j_\tau\otimes j_\sigma\circ\vartheta_{\sigma,\tau}\circ
S\circ j_{\rho}^{-1},$$
with $\vartheta_{\sigma,\tau}: \sigma\otimes\tau\to\tau\otimes\sigma$ the
flip operator.
Therefore the left hand side of the equation we want to establish equals
$$\sum_r\ell_{J_r{\cal F}(S_r)^*T\otimes
T'}^*c_{\overline{u_r}}(\overline{S_r}^*
{j_v^*}^{-1}\phi'\otimes{j_u^*}^{-1}\phi)=$$
$$\sum_r\ell_{J_r{\cal F}(S_r)^*T\otimes
T'}^*{\cal
F}({\overline{S_r}}^*)\otimes I\circ
c_{\overline{v}\otimes\overline{u}}({j_v^*}^{-1}\phi'\otimes{j_u^*}^{-1}\phi)=$$
$$\sum_r\ell_{J_r\circ{\cal F}(S_r)^*\circ\vartheta_{v,u}\circ 
J_v^{-1}\otimes J_u^{-1}(J_v(T')\otimes
J_u(T))}^*{\cal
F}({\overline{S_r}}^*)\otimes I\circ
c_{\overline{v}\otimes\overline{u}}({j_v^*}^{-1}\phi'\otimes{j_u^*}^{-1}\phi)=$$
$$\sum_r\ell_{J_r\circ{\cal F}(S_r)^*\circ\vartheta_{v,u}\circ 
J_v^{-1}\otimes J_u^{-1}(J_v(T')\otimes
J_u(T))}^*
({J_r^{-1}}^*{\cal F}(S_r^*)J_{u\otimes v}^*\otimes I)
\circ
c_{\overline{v}\otimes\overline{u}}
({j_v^*}^{-1}\phi'\otimes{j_u^*}^{-1}\phi)=$$
$$\sum_r\ell_{J_{u\otimes v}{\cal F}(S_r)J_r^{-1}J_r{\cal
F}(S_r)^*\vartheta_{v,u} 
J_v^{-1}\otimes J_u^{-1}(J_v(T')\otimes
J_u(T))}^*
\circ
c_{\overline{v}\otimes\overline{u}}
({j_v^*}^{-1}\phi'\otimes{j_u^*}^{-1}\phi)=$$
$$\ell_{J_{u\otimes v}(T\otimes T')}^*
\circ
c_{\overline{v}\otimes\overline{u}}
({j_v^*}^{-1}\phi'\otimes{j_u^*}^{-1}\phi)=$$
$$\ell_{J_u(T)\otimes J_v(T')}^*
\circ
c_{\overline{v}\otimes\overline{u}}
({j_v^*}^{-1}\phi'\otimes{j_u^*}^{-1}\phi)$$
thanks to Prop. 3.9.
So far we have obtained a unital $^*$--algebra ${\cal C}_{\cal F}$.
Consider the linear map 
$$\eta_{\cal F}:{\cal C}_{\cal F}\to{\cal C}_{\cal F}\odot{\cal
A}_\infty,$$
$$\eta_{\cal F}(\overline{T}\otimes\phi)=\overline{T}\otimes u(\phi),\quad
T\in{\cal F}(u),\phi\in H_u, u\in\hat{\cal F}.$$
It is easy to check that $\eta_{\cal F}$ is multiplicative and that
$$\eta_{\cal
F}\otimes\iota\circ\eta_{\cal F}=\iota\otimes\Delta\circ\eta_{\cal
F}.$$
We show that 
$\eta_{\cal F}$ preserves the involutions. Recall that the relation 
$\overline{R}\in(\iota,u\otimes\overline{u})$ 
is equivalent to the fact that the linear operator $Q:\overline{\psi}\in
\overline{H_u}\to {j^*}^{-1}\psi\in H_{\overline{u}}$ lies in $(u_*,
\overline{u})$, so
$$\eta_{\cal F}((\ell_T^*c_u(e_i))^*)=\eta_{\cal F}
(\ell_{JT}^*c_{\overline{u}}(Q\overline{e_i}))=$$
$$JT\otimes\overline{u}(Q\overline{e_i})=
\sum_k JT\otimes
{j^*}^{-1}e_k\otimes u_{ki}^*=$$
 $$\eta_{\cal F}(\ell_T^*c_u(e_i))^*.$$
We introduce a state on ${\cal C}_{\cal F}$.
Consider the linear functional 
$h$ on ${\cal C}_{\cal F}$ defined by
$$h(\overline{T}\otimes\phi)=0 \quad T\in{\cal F}(u), \phi\in H_u,
u\in\hat{{\cal F}}, u\neq
\iota,$$
$$h(I)=1.$$
We claim that $h$  is  a positive faithful state on ${\cal
C}_{\cal F}$, so ${\cal
C}_{\cal F}$ has
a $C^*$--norm.

Consider the conditional expectation $E$ onto the fixed point
$^*$--subalgebra
obtained averaging over the group action. Since the Haar measure of $G$
annihilates the coefficients of all the irreducible representations $u$,
except for $u=\iota$, we see that
 $E=h$, and
we can conclude that $\eta_{\cal F}$ is an
ergodic algebraic action.

It is now evident that, for $u\in\hat{{\cal F}}$,
the maps, for  $T\in{\cal F}(u)$,
$$\gamma_T: \psi\in H_u\to \overline{T}\otimes\psi\in{\cal C}_{\cal F}$$
are elements of the Hilbert space $L_u$ associated with the algebraic 
ergodic space,
and that elements of this form span $L_u$.
Therefore we obtain  surjective linear maps, for $u\in\hat{\cal F}$,
$$V_u: T\in {\cal F}(u)\to\overline{\gamma_T}\in\ov{L}_u.$$
We show that $V_u$ is an isometry, and therefore a unitary.
For $S, T\in{\cal F}(u)$, and an orthonormal basis $(\psi_i)$ of $H_u$,
$$(V_u(T), V_u(S))I=\sum_i\gamma_T(\psi_i)\gamma_S(\psi_i)^*=$$
$$\sum_i(\overline{T}\otimes\psi_i)(\overline{S}\otimes\psi_i)^*=
\sum_i(\overline{T}\otimes\psi_i)(\overline{J S}\otimes
{j^*}^{-1}\psi_i).$$
Consider pairwise inequivalent irreducible representations $u_1,\dots
u_N$ of $G$ and, for $\alpha=1,\dots, N$,  isometries $s_{\alpha,
j}\in(u_\alpha,
u\otimes\overline{u})$, with ranges adding up to the identity. 
Then the last term above becomes
$$\sum_{i}\sum_{\alpha}\sum_j\overline{{\cal F}(s_{\alpha, j})^*(T\otimes
JS})\otimes s_{\alpha,j}^*(\psi_i\otimes {j^*}^{-1}\psi_i)=$$
$$\sum_{\alpha}\sum_j\overline{{\cal F}(s_{\alpha, j})^*(T\otimes
JS})\otimes s_{\alpha,j}^*\overline{R}_u.$$
Now 
$s_{\alpha, j}^*\overline{R}_u\in(\iota, u_\alpha)$, 
hence
$s_{\alpha, j}^*\overline{R}_u=0$ unless $u_\alpha=\iota$. In this case 
$N=1$ and $s_{\alpha_j}$ can be chosen to coincide with
$\overline{R}_u\|\overline{R}_u\|^{-1}$. We conclude that the last term
above equals
$$\|\overline{R}_u\|^{-2}\overline{{\cal F}(\overline{R}_u)^*(T\otimes
JS)}\otimes
\overline{R}_u^*\overline{R}_u=\overline{\hat{\overline{R_u}}^*T\otimes
JS}\otimes 1_\iota=$$
$$\overline{(S,T)1_\iota}\otimes 1_\iota=(T, S)I_\iota. $$
We extend $V$ to a quasitensor natural transformation from $\ov{L}$ to
${\cal F}$. First set $V_u=0$ if $u$ is irreducible but ${\cal F}(u)=0$.
Then
for any representation $u\in\text{Rep}(G)$, consider pairwise inequivalent 
irreducible representations $u_\alpha\in\hat\cF$ and isometries 
$s_{\alpha,j}\in(u_\alpha, u)$
decomposing $u$. Define a unitary map $V_u: {\cal F}(u)\to\ov{L}_u$ by 
$$V_u=\sum_\alpha\sum_j\ov{L}_{s_{\alpha,j}} V_{u_\alpha}{\cal
F}(s_{\alpha,j})^*.$$
It is easy to check that $V_u$ does not depend on the choice of the
isometries, it is a natural transformation and that   one now has, as a
consequence of the definition of multiplication in ${\cal C}_{\cal F}$,
that for $u,v\in\hat{{\cal F}}$,
$$V_{u\otimes v}\upharpoonright_{{\cal F}(u)\otimes{\cal F}(v)}=V_u\otimes
V_u.$$ From this relation one concludes, with routine computations, that
$V$ is a quasitensor unitary natural transformation between $\ov{L}$
and ${\cal F}$.
Recalling the definition of the maps $c_u$ at the
beginning of the proof, one also has that, for $\psi\in H_u$,
$u\in\hat{\cal F}$,
$$V_u\otimes I
c_u(\psi)=\sum_k\overline{\gamma_{T_k}}\otimes\gamma_{T_k}(\psi),$$
and $c_u$ is the multiplicity map associated with the functor $\ov{L}$ in
the representation $u$.
Therefore we can now say that
 ${\cal C}_{\cal F}$ is linearly spanned by entries 
of coisometries 
in matrix algebras over
${\cal C}_{\cal F}$,
the  maps $c_u$,
 so the maximal $C^*$--seminorm on ${\cal C}_{\cal F}$ is finite. 
 Completing ${\cal C}_{\cal F}$ in the maximal $C^*$--seminorm
yields a 
unital $C^*$--algebra ${\cal B}_{\cal F}$, a nondegenerate
ergodic $G$--action
$$\eta_{\cal F}:{\cal B}_{\cal F}\to{\cal B}_{\cal F}\otimes{\cal A},$$
with spectral functor $\ov{L}$,
and a natural unitary transformation from $\ov{L}$ to 
${\cal F}$.\medskip

We now claim the following.\medskip

\noindent{\bf 9.2 Lemma} {\sl If $S\in(u,v)$ and if 
$j_u: H_u\to H_{\overline{u}}$, $j_v: H_v\to H_{\overline{v}}$ are 
antilinear invertible maps defining conjugates of $u$ and $v$
respectively, with corresponding antilinear invertible maps
$J_u:{\cal F}(u)\to{\cal F}(\overline{u})$
$J_v:{\cal F}(v)\to{\cal F}(\overline{v})$, in the sense of Theorem 3.7,
then 
$${\cal F}( j_v S j_{u}^{-1})=J_v{\cal F}(S) J_{u}^{-1}.$$}\medskip

\noindent{\bf Proof} Let $R_u$, $\overline{R}_u$ be the solution to the
conjugate equations corresponding to $j_u$, and, similarly, $R_v$,
$\overline{R}_v$ the solution corresponding to $j_v$. We need to show that
$${\cal F}(1_{\overline{v}}\otimes \overline{R}_u^*\circ
1_{\overline{v}}\otimes S^*\otimes 1_{\overline{u}}\circ R_{v}\otimes
1_{\overline{u}})=
1_{{\cal F}(\overline{v})}\otimes \hat{\overline{R}_u}^*\circ
1_{{\cal F}(\overline{v})}\otimes{\cal
F}(S)^*\otimes 1_{{\cal F}(\overline{u})}\circ
\hat{R_v}\otimes 1_{{\cal F}(\overline{u})}.$$
We shall use the simplified notation $(3.7)$--$(3.11)$,  
replacing $S_{u,v}$ by the identity, $S_{u,v}^*$ by $E_{u,v}$,
$\hat{R}$ with $E_{\overline{u}, u}\circ R$ and 
$\hat{\overline{R}}$ by $E_{u, \overline{u}}\circ \overline{R}$.
Properties $(3.11)$ and $(3.8)$ show that, regarding ${\cal F}(R_v)$ as an
element of ${\cal F}(\overline{v}\otimes v)$, 
$${\cal F}(R_v\otimes
1_{\overline{u}})={\cal F}(R_v)\otimes 1_{{\cal F}(\overline{u})}.$$
We next show that
 $${\cal F}(1_{\overline{v}}\otimes \overline{R}_u^*)= 
1_{{\cal
F}(\overline{v})}\otimes {\cal F}(\overline{R_u})^*\circ E_{\overline{v},
u\otimes\overline{u}}.$$
In fact, if $\xi\in {\cal F}(\overline{v}\otimes u\otimes \overline{u})$
and $\eta\in {\cal F}(\overline{v})$,
$$(\eta, {\cal F}(1_{\overline{v}}\otimes \overline{R}_u^*)\xi)=
(\eta, {\cal F}(1_{\overline{v}}\otimes \overline{R}_u)^*\xi)=$$
$$({\cal F}(1_{\overline{v}}\otimes \overline{R}_u)\eta,\xi)=
(\eta\otimes{\cal F}(\overline{R}_u), \xi)=$$
$$(\eta\otimes{\cal F}(\overline{R}_u),E_{\overline{v},
u\otimes\overline{u}}\xi)=
(\eta,
1_{{\cal F}(\overline{v})}\otimes{\cal
F}(\overline{R}_u)^*\circ E_{\overline{v},
u\otimes\overline{u}}\xi).$$
Pick a vector $\zeta\in {\cal F}(\overline{u})$. Then
$${\cal F}(1_{\overline{v}}\otimes \overline{R}_u^*\circ
1_{\overline{v}}\otimes S^*\otimes 1_{\overline{u}}\circ R_{v}\otimes
1_{\overline{u}})\zeta=$$
$$1_{{\cal F}(\overline{v})}\otimes{\cal F}(\overline{R}_u^*)\circ
E_{\overline{v}, u\otimes\overline{u}}(
{\cal F}(1_{\overline{v}}\otimes S)^*({\cal
F}(R_v))\otimes\zeta)\eqno(9.1)$$
But thanks to $(3.12)$,
$$E_{\overline{v}, u\otimes\overline{u}}(
{\cal F}(1_{\overline{v}}\otimes S)^*(
{\cal F}(R_v))\otimes\zeta)=
E_{\overline{v}, u,\overline{u}}({\cal F}(1_{\overline{v}}\otimes
S)^*({\cal F}(R_v))\otimes\zeta)=$$
$$E_{\overline{v},u}({\cal F}(1_{\overline{v}}\otimes
S)^*{\cal F}(R_v))\otimes\zeta=$$
$$E_{\overline{v},u}\circ {\cal F}(1_{\overline{v}}\otimes S)^*({\cal
F}(R_v))\otimes \zeta=
1_{{\cal F}(\overline{v})}\otimes
{\cal F}(S)^*(E_{\overline{v},v}\circ
{\cal F}(R_v))\otimes\zeta=$$
$$(1_{{\cal F}(\overline{v})}\otimes
{\cal F}(S)^*(\hat{R_v}))\otimes\zeta=
1_{{\cal
F}(\overline{v})}\otimes{\cal F}(S)^*\otimes 1_{{\cal
F}(\overline{u})}\circ\hat{R_v}\otimes 1_{{\cal
F}(\overline{u})}(\zeta)=$$
$$1_{{\cal F}(\overline{v})}\otimes E_{u,\overline{u}}\circ
1_{{\cal
F}(\overline{v})}\otimes{\cal F}(S)^*\otimes 1_{{\cal
F}(\overline{u})}\circ\hat{R_v}\otimes 1_{{\cal
F}(\overline{u})}(\zeta)$$
Substituting back in $(9.1)$ gives
$$1_{{\cal F}(\overline{v})}\otimes \hat{\overline{R}}_u^*
\circ 1_{{\cal
F}(\overline{v})}\otimes{\cal F}(S)^*\otimes 1_{{\cal
F}(\overline{u})}\circ\hat{R_v}\otimes 1_{{\cal
F}(\overline{u})}(\zeta).$$
\medskip

\noindent{\bf 9.3 Lemma} {\sl The linear functional $h$ is a faithful
state 
on the
$^*$--algebra ${\cal
C}_{\cal F}$.}\medskip

\noindent{\bf Proof} For $u,v\in\hat{{\cal F}}$,
pick $T\in{\cal F}(u)$, $S\in{\cal F}(v)$, $\phi\in H_u$, $\psi\in H_v$.
Choose isometries $s_\alpha\in (u_\alpha, \overline{u}\otimes v)$,
with $u_\alpha$ irreducible, with orthogonal ranges and summing up to the
identity. Then
$$(\overline{T}\otimes\phi)^* (\overline{S}\otimes\psi)=
\ell_{J_uT\otimes
S}^*\circ c_{\overline{u}\otimes v}({j_u^*}^{-1}\phi\otimes\psi)=$$
$$\sum_\alpha\ell_{J_uT\otimes
S}^*\circ {\cal F}(s_\alpha)\otimes I\circ c_{
u_{\alpha}}(s_\alpha^*{j_u^*}^{-1}\phi\otimes\psi).$$
Let $A$ be the subset of all $\alpha$ for which $u_\alpha=\iota$.
Then for $\alpha\in A$, $(\iota, \overline{u}\otimes v)$ is always
zero unless $u=v$. If this is the case, there is, up to a phase, just one
$s_\alpha$,
which is of the form
$\lambda_u R_u$, with $\lambda_u=\|R_u\|^{-1}$, because $u$ is
irreducible and the space $(\iota,
\overline{u}\otimes u)$ is one--dimensional.
Thus 
$$h((\overline{T}\otimes\phi)^*
(\overline{S}\otimes\psi))=
\delta_{u,v}|\lambda_u|^2({\cal
F}(R_u)^* J_uT\otimes S)^*R_u^*{j_u^*}^{-1}\phi\otimes\psi)=$$
$$\delta_{u,v}|\lambda_u|^2(\hat{R_u}^* J_uT\otimes
S)^*(R_u^*{j_u^*}^{-1}\phi\otimes\psi))=
\delta_{u,v}|\lambda_u|^2(JT, JS)(\phi,\psi).$$
If $a\in{\cal C}_{\cal F}$ is written in the form
$a=\sum_{u\in F,i,j}\mu^u_{ij}\overline{T^u_i}\otimes\phi^u_j$ with 
$F$ a finite set, $T^u_i\in{\cal F}(u)$, $\phi^u_j\in H_u$ orthonormal
bases, the previous computation gives
$$h(a^*a)=\sum_{u\in
F}\sum_{i,r,j}\text{Trace}(j_u^*j_u)^{-1}\overline{\mu^u_{ij}}\mu^u_{rj}
(J_uT_i, J_uT_r)\geq0.$$
If $h(a^*a)=0$ then we can conclude that for all $j$, $\sum_i
\mu^u_{ij}J_uT_i=0$, hence $\mu^u_{ij}=0$ for all $i,j,u$, so $a=0$.
\medskip

One can induce $^*$--isomorphisms between two induced $C^*$--systems 
$({\cal B}_{\cal F},\eta_{\cal F})$ and $({\cal B}_{\cal G},\eta_{\cal
G})$ using natural transformations.
\medskip

\noindent{\bf 9.4 Proposition} {\sl Let ${\cal F},{\cal
G}:\text{Rep}(G)\to{\cal H}$ be two quasitensor $^*$--functors
and let $U:{\cal F}\to{\cal G}$ be a unitary quasitensor natural
transformation. Then there exists a unique $^*$--isomorphism
$\alpha_u:{\cal B}_{\cal F}\to {\cal B}_{\cal G}$ intertwining the 
corresponding $G$--actions, such that
$$\alpha_U(\overline{T}\otimes\phi)=\overline{U_uT}\otimes\phi,\quad
T\in{\cal F}(u), \phi\in H_u, u\in\hat{\cal F}.$$}\medskip

\noindent{\bf Proof} It is easy to check that that formula defines a
linear multiplicative map
 $\alpha_u$ commuting with the actions.
It is also easy to check that, for $T\in{\cal F}(u)$,
$$U_{\overline{u}}r_T^*\circ E^{\cal F}_{\overline{u},u}=r_{U_uT}^*\circ
E^{\cal G}_{\overline{u},u}\circ U_{\overline{u}\otimes u}.$$
Thus $$U_{\overline{u}}J^{\cal F}_u(T)=
U_{\overline{u}}r_T^*\circ E^{\cal F}_{\overline{u},u}\circ {\cal
F}(R_u)=$$
$$r_{U_uT}^*\circ E^{\cal G}_{\overline{u},u}\circ U_{\overline{u}\otimes
u}{\cal F}(R_u)=$$
$$r_{U_uT}^*\circ E^{\cal G}_{\overline{u},u}\circ 
{\cal G}(R_u)=J_u^{\cal G}U_uT,$$
which shows that $\alpha_u$ is $^*$--invariant. Thus $\alpha_U$ extends to
the completions in the maximal $C^*$--seminorms.
\medskip

\end{section}

\begin{section} {Applications to abstract duality theory}

As a first application of Theorem 9.1, we consider a tensor
$C^*$--category ${\cal T}$ and a faithful tensor $^*$--functor
$\rho:\text{Rep}(G)\to{\cal T}$ from the representation
category of a compact quantum group $G$ to ${\cal T}$. 
Then, as shown in Example 3.5, one has a quasitensor $^*$--functor
${\cal F}_\rho:\text{Rep}(G)\to{\cal H}$
which associates to any representation $u$ of $G$ the Hilbert space
$(\iota,\rho_u)$. 
Therefore we can apply Theorem 9.1 to ${\cal F}_\rho$ and obtain an
ergodic $G$--space canonically associated with the inclusion 
$\rho$. We summarize this in the following theorem.
\medskip

\noindent{\bf 10.1 Theorem} {\sl Let $\rho:\text{Rep}(G)\to{\cal T}$
be a tensor $^*$--functor from the representation category
of a compact quantum group $G$ to an abstract tensor $C^*$--category
${\cal T}$.  Then  there is a canonically associated ergodic
nondegenerate action of $G$ on a unital $C^*$--algebra ${\cal B}$ 
whose spectral functor can be identified with ${\cal F}_\rho$.}\medskip

The above $G$--space plays a central role in abstract duality theory 
for compact quantum groups. This matter will be developed elsewhere
\cite{DPR}. Here we  notice some  consequences of the previous
corollary.

The notion of {\it permutation symmetry} for an 
abstract  tensor $C^*$--category
has been introduced in \cite{DRinventiones}.
If $G$ is a compact group, $\text{Rep}(G)$ has permutation symmetry
defined by the intertwiners $\vartheta_{u,v}\in (u\otimes v, v\otimes u)$ exchanging the
order of factors in the tensor product.
 The previous theorem allows 
us to recover the following result.
\medskip

\noindent{\bf 10.2 Theorem} {\sl If $G$ is a compact group, ${\cal T}$ 
is a tensor $C^*$--category with permutation symmetry $\varepsilon$ and 
$\rho:\text{Rep}(G)\to{\cal T}$ is a tensor $^*$--functor such that
$\rho(\vartheta_{u,v})=\varepsilon(\rho_u,\rho_v)$
then
there exists a compact subgroup $K$ of $G$, unique up to conjugation,
and an isomorphism of the $G$--ergodic system associated with $\rho$ and
the ergodic
$C^*$--system induced by 
 the homogeneous 
space $K\backslash G$ over $G$.}\medskip

\noindent{\bf Proof} We show that the $C^*$--algebra associated with the 
quasitensor $^*$--functor $u\to(\iota,\rho_u)$ is commutative. 
This will suffice, as it is well known that any ergodic action 
of a compact group on a commutative $C^*$--algebra arises from the
transitive $G$--action 
on a quotient space $K\backslash G$ by a point stabilizer subgroup,
unique up to conjugation.
We need to show 
that if $u$ and $v$ are irreducible then for $k\in(\iota,\rho_u)$,
$k'\in(\iota,\rho_v)$, $\psi\in H_u$, $\psi'\in H_v$ then
$\overline{k}\otimes\psi$ and $\overline{k'}\otimes\psi'$ commute, as 
these elements span a dense $^*$--subalgebra. Choose inequivalent irreducible 
representations $(u_\alpha)$ of $G$ and isometries 
$s_{\alpha,j}\in (u_\alpha, u\otimes v)$ such that 
$\sum_\alpha\sum_js_{\alpha,j}s_{\alpha,j}^*=1_{u\otimes v}$.
Thus 
$$(\overline{k}\otimes\psi)(\overline{k'}\otimes\psi')=
\sum_\alpha\sum_j\overline{\widehat{\rho(s_{\alpha,j}^*)}(k\otimes 1_{\rho_v}\circ k')}
\otimes s_{\alpha,j}^*(\psi\otimes \psi').\eqno(10.1)$$
Now $$\vartheta_{v,u}\psi'\otimes\psi=\psi\otimes\psi'$$ and
$$\varepsilon(\rho_v,\rho_u)(k'\otimes 1_{\rho_u}\circ k)=
\varepsilon(\rho_v,\rho_u)(k'\otimes k)=k\otimes k'=
k\otimes 1_{\rho_v}\circ k'.$$
Thus the right hand side of $(10.1)$ can be written
$$\sum_\alpha\sum_j\overline{\widehat{\rho(s_{\alpha,j}^*)}
(\varepsilon(\rho_v,\rho_u)k'\otimes 1_{\rho_u}\circ k})\otimes 
s_{\alpha,j}^*\vartheta_{v,u}(\psi'\otimes\psi).\eqno(10.2)$$
On the other hand $\varepsilon(\rho_u,\rho_v)\circ \widehat{\rho(s_{\alpha,j})}=
\rho(\vartheta_{u,v})\circ \widehat{\rho(s_{\alpha,j})}$, and this is 
the map that takes an element $\xi\in (\iota, \rho_\alpha)$ to
the element $$\rho(\vartheta_{u,v})\circ\rho(s_{\alpha,j})\circ\xi=
\rho(\vartheta_{u,v}\circ s_{\alpha,j})\circ\xi\in(\iota, \rho_v\otimes\rho_u).$$
Thus 
$$\varepsilon(\rho_u,\rho_v)\circ \widehat{\rho(s_{\alpha,j})}=
\widehat{\rho(\vartheta_{u,v}\circ s_{\alpha,j})}.$$
Set
$t_{\alpha,j}:=\vartheta_{u,v}\circ s_{\alpha,j}\in
(u_\alpha, v\otimes u)$. Then 
 $(10.2)$ equals 
$$\sum_\alpha\sum_j\overline{\widehat{\rho(t_{\alpha,j}^*)}
(k'\otimes 1_{\rho_u}\circ k )}\otimes t_{\alpha,j}^*(\psi'\otimes\psi).$$
Since the isometries $t_{\alpha,j}:=\vartheta_{u,v}\circ s_{\alpha,j}\in
(u_\alpha, v\otimes u)$ give an orthogonal decomposition of $v\otimes u$ into
irreducibles, the last term above equals 
$(\overline{k'}\otimes\psi')(\overline{k}\otimes\psi)$, and
 the proof is complete.
\medskip

Recall that a {\it $q$--Hecke symmetry } for an object $\rho$ in a tensor
$C^*$--category ${\cal T}$
is given by representations 
$$\varepsilon_n: H_n(q)\to (\rho^{\otimes n},\rho^{\otimes n})$$
of the Hecke algebras $H_n(q)$
such that
$$\varepsilon_{n+1}(b)=\varepsilon_n(b)\otimes 1_\rho,\quad b\in H_n(q)\subset
H_{n+1}(q),$$
$$\varepsilon_{n+1}(\sigma(b))=1_\rho\otimes\varepsilon_n(b),\quad b\in
H_n(q),$$
where $\sigma: H_n(q)\to H_{n+1}(q)$ is the homomorphism taking each
generator $g_i$ of $H_n(q)$ to $g_{i+1}$.

Also recall that 
an object $\rho$ of ${\cal T}$ is called {\it $\mu$--special of dimension
$d$} if 
there is 
a $\mu^2$--Hecke symmetry for $\rho$  and an
intertwiner
$R\in(\iota, \rho^{\otimes d})$ for some $d\geq 2$, satisfying
$$R^*R=d!_q,\eqno(10.3)$$
$$R^*\otimes 1_\rho\circ 1_\rho\otimes R=(d-1)!_q(-\mu)^{d-1}
1_\rho,\eqno(10.4)$$
$$RR^*=\varepsilon_d(A_d),\eqno(10.5)$$
$$\varepsilon(g_1\dots g_d)R\otimes 1_\rho=-(-\mu)^{d-1} 1_\rho\otimes
R,\eqno(10.6)$$
where $q:=\mu^2$,
see \cite{P1} for notation. 
\medskip

If ${\cal T}$ admits a special object of dimension $d$ for some $\mu>0$,
Theorem 6.2 in \cite{P1} then assures then 
the existence of a tensor $^*$--functor $\text{Rep}(S_\mu U(d))\to {\cal
T}$
taking the fundamental representation $u$ of $S_\mu U(d)$ to $\rho$ and
the
canonical
intertwiner $S\in(\iota, u^{\otimes d})$ to $R$, where $u$ is the fundamental
representation of   $(S_\mu U(d))$.
 We are thus in a position to apply 
Theorem 10.1.\medskip

\noindent{\bf 10.3 Theorem} {\sl Let $\rho$ be a $\mu$--special object of
a tensor $C^*$--category ${\cal
T}$
with dimension $d\geq2$ and parameter $\mu>0$.
Then there exists an ergodic nondegenerate action of $S_\mu U(d)$
on a unital $C^*$--algebra ${\cal B}$ whose spectral functor can be 
identified on the objects with $u^{\otimes r}\to (\iota,\rho^r)$, with $u$
the fundamental representation of $S_\mu U(d)$.}\medskip

In the particular case where $d=2$ the notion of a $\mu$--special object of
dimension $2$ simplifies considerably.
\medskip

\noindent{\bf 10.4 Proposition} {\sl If
 an object $\rho$ of a tensor $C^*$--category
${\cal T}$
admits an intertwiner $R\in(\iota,\rho^{\otimes 2})$ satisfying relations 
$$R^*\circ R=(1+q) 1_\iota,\eqno(10.7)$$
$$R^*\otimes 1_\rho\circ 1_\rho\otimes R=-\mu 1_\rho,\eqno(10.8)$$
with $\mu>0$ and $q=\mu^2$, then $\rho$ can be made uniquely into a
$\mu$--special
object of dimension $2$ with intertwiner $R$}\medskip

\noindent{\bf Proof} If $\varepsilon$ is any $\mu^2$--Hecke symmetry for
$\rho$ making $\rho$ into a $\mu$--special object of dimension $2$ with
intertwiner $R$, then 
$(10.5)$ shows that $RR^*=\varepsilon_2(A_2)$. Since $A_2=1+g_1$,
$\varepsilon_2(g_1)=RR^*-1_{\rho^{\otimes 2}}$. Therefore for all $n$,
$i=1,\dots, n-1$,
$$\varepsilon_n(g_i)=\varepsilon_n(\sigma^{i-1}(g_1))=$$
$$1_{\rho^{\otimes
i-1}}\otimes\varepsilon_2(g_1)
=1_{\rho^{\otimes i-1}}\otimes(RR^*-1_{\rho^{\otimes 2}})$$
and the symmetry is uniquely determined. Let us then show the existence of 
a symmetry
using that formula.
 The orthogonal projection
$e=\frac{1}{1+q}RR^*\in(\rho^2,\rho^2)$
satisfies
the Temperley--Lieb relations 
$$e\otimes 1_\rho\circ 1_\rho\otimes e\circ e\otimes
1_\rho=(q+\frac{1}{q}+2)^{-1} e\otimes 1_\rho,$$
$$1_\rho\otimes e\circ e\otimes 
1_\rho\circ 1_\rho\otimes e=(q+\frac{1}{q}+2)^{-1} 1_\rho\otimes e$$
thanks to $(10.8)$.
Since the Temperley--Lieb algebra $TL_n((q+\frac{1}{q}+2)^{-1})$  is the
quotient of the Hecke algebra $H_n(q)$  
by the ideal generated by $A_3$ (see \cite{GHJ}), 
there does exist a  Hecke symmetry for $\rho$ such that
$\varepsilon_2(g_1)=(q+1)e-1_{\rho^{\otimes 2}}=RR^*-1_{\rho^{\otimes
2}}$.
Since $A_2=1+g_1$, $\varepsilon_2(A_2)=(q+1) e=R\circ R^*$,
so $(10.5)$ follows. 
We show $(10.6)$ for $d=2$:
$$\varepsilon_3(g_2)\circ R\otimes
1_\rho=(1_\rho\otimes(RR^*-1_{\rho^{\otimes 2}}))\circ R\otimes 1_\rho=$$
$$-\mu 1_\rho\otimes R-R\otimes 1_\rho,$$
hence
$$\varepsilon_3(g_1g_2)\circ R\otimes 1_\rho=((RR^*-1_{\rho^{\otimes
2}})\otimes 1_\rho)\circ
(-\mu 1_\rho\otimes R-R\otimes 1_\rho)=$$
$$\mu^2 R\otimes 1_\rho-(q+1)R\otimes 1_\rho
+\mu 1_\rho\otimes R+R\otimes 1_\rho=\mu 1_\rho\otimes R.$$
\medskip

\noindent{\it  Remark}
There is a canonical isomorphism from 
$\text{Rep}(S_\mu
U(2))$ to  $\text{Rep}(A_o(F))$ 
\medskip

Relations $(10.7)$ and
$(10.8)$ can be
implemented in Hilbert spaces with dimension $\geq2$ in the following way.
Let $j$ be any antilinear invertible map on a finite dimensional Hilbert
space $H$, and set $R=\sum_i je_i\otimes e_i\in H^{\otimes 2}$. 
Then, for this $R$, $(10.7)$ and $(10.8)$ become, 
$$\text{Trace}(j^*j)=1+q,$$
$$j^2=-\mu.$$
 Consider the involutive antiunitary map $c$ of $H$ acting trivially 
on the orthonormal basis $e_i$, and set $F=jc$ and $\overline{F}=cFc$.
Then the above conditions can be equivalently written
$$\text{Trace}(F^*F)=1+q,$$
$$F\overline{F}=-\mu.$$
The maximal compact quantum group with representation category generated 
by $R$ is the universal quantum group  $A_o(F)$ defined in
\cite{VW}. 
Therefore the fundamental representation of $A_o(F)$ is a $\mu$--special 
object of dimension $2$ in $\text{Rep}(A_o(F))$. Theorem 6.2 in \cite{P1}
then shows that
there is a unique isomorphism of
tensor $C^*$--categories $\text{Rep}(S_\mu U(2))\to\text{Rep}(A_o(F))$
taking the fundamental representation $u$ of $S_\mu U(2)$ to
the fundamental reprsentation of $A_o(F)$ and the quantum determinant
$S=\psi_1\otimes\psi_2-\mu\psi_2\otimes\psi_1\in(\iota, u^{\otimes 2})$
to $R=\sum j e_i\otimes e_i$.

For  related results,
see Cor.\ 5.4 in \cite{BRV}, where the authors find similar necessary and
sufficient conditions for  the existence of a
monoidal equivalence between a generic pair of  universal compact
quantum groups, and the result of Banica \cite{Banica}, where it is shown
that the fusion rules of $A_o(F)$ are the same as those of $SU(2)$.
\medskip

\noindent{\bf 10.5  Corollary} {\sl Let $\rho$ be an object of a tensor
$C^*$--category ${\cal T}$ with an intertwiner
$R\in(\iota,\rho^{\otimes 2})$ satisfying conditions $(10.7)$ and
$(10.8)$.
Then the unique tensor $^*$--functor $\text{Rep}(S_\mu U(2))\to {\cal
T}$ 
taking the fundamental representation $u$ to $\rho$  and the 
quantum determinant
$S=\psi_1\otimes\psi_2-\mu\psi_2\otimes\psi_1\in(\iota, u^{\otimes 2})$ 
to $R$
gives rise
to
an ergodic nondegenerate action of $S_\mu U(2)$ on a unital
$C^*$--algebra ${\cal B}$ with spectral subspaces $\ov{L}_{u^{\otimes
r}}=(\iota,\rho^r)$.}\medskip
\end{section}

\begin{section} {Actions embeddable into quantum
quotient spaces}

As a second application of the duality theorem 9.1, consider the
invariant vectors functor
$K$ associated with a compact quantum subgroup $K$ of $G$. 
We know that this is just a copy of the spectral functor of the 
quotient space $K\backslash G$.
\medskip

\noindent{\bf 11.1 Theorem} {\sl Let $K$ be a compact quantum subgroup of
a
maximal compact quantum group $G$. Then the nondegenerate ergodic
$G$--system associated with the invariant vectors functor $K$  is isomorphic
to the quotient $G$--space
$K\backslash G.$}\medskip

\noindent{\bf Proof} It is clear from the construction of the dense Hopf 
$^*$--algebra ${\cal C}_K$ in Proposition 2.6 that ${\cal C}_K$ is 
$^*$--isomorphic 
to the dense spectral subalgebra of ${\cal A}^\delta$, in such a way that
the constructed $G$--action  corresponds to the right $G$--action on the
right coset space. Therefore we are left to show that the maximal
$C^*$--seminorm $\|\ \dot{}\ \|_1$ on ${\cal C}_K={\cal A}^\delta_{sp}$
coincides with the
restriction of the maximal 
$C^*$--seminorm $\|\ \dot{}\ \|_2$ on ${\cal A}_\infty$. 
Any Hilbert space representation of ${\cal A}_\infty$ restricts to a
Hilbert space representation on ${\cal A}^\delta_{sp}$,
so $\|\ \dot{}\ \|_2\leq \|\ \dot{}\ \|_1$. Conversely, if 
$\pi$ is a Hilbert space representation for ${\cal A}^{\delta}_{sp}$, we 
can induce it up to a Hilbert space representation $\tilde{\pi}$ of ${\cal 
A}_\infty$ via the conditional expectation $m:{\cal A}_\infty\to{\cal
A}^\delta_{sp}$ obtained averaging over the action $\delta$ of the
subgroup. There is an isometry $V$ from the Hilbert space of $\pi$ to the
Hilbert space of $\tilde{\pi}$ such that $V^*\tilde{\pi}(a)V=\pi(m(a)),$
for $a\in{\cal A}_\infty$. In particular, if $a\in{\cal A}^\delta_{sp}$,
$\pi(a)=V^*\tilde{\pi}(a)V$, so
$\|\pi(a)\|\leq \|\tilde{\pi}(a)\|\leq\|a\|_2$, and $\|a\|_1\leq\|a\|_2$.
\medskip

\noindent{\bf 11.2 Definition}  An ergodic $G$--action $\eta:{\cal
B}\to{\cal
B}\otimes{\cal A}$ will be called {\it maximal} 
if $G$ is a maximal compact quantum group and if ${\cal B}$ is obtained
completing the  dense spectral $^*$--subagebra with respect to the maximal
$C^*$--seminorm.\medskip

Combining the previous result with the abstract characterization of the
invariant vectors functor $K$ given in Theorem 5.5 and with Prop. 9.4,
gives the following characterization of maximal ergodic systems isomorphic
to quotient
spaces.
\medskip

\noindent{\bf 11.3 Theorem} {\sl Let $({\cal B},\eta)$  be a maximal,
nondegenerate,
ergodic
$G$--action. If there exists, for each unitary representation $u$ of $G$, a 
subspace $K_u\subset H_u$ satisfying properties $(5.1)$, $(5.2)$, $(5.7)$,
$(5.8)$ and a quasitensor natural unitary transformation 
from the spectral functor $\ov{L}$ associated with $({\cal B},\eta)$ to
the functor $K$, then there exists 
a unique maximal compact quantum subgroup $K$ of $G$ such that 
$({\cal B},\eta)\simeq K\backslash G$.}\medskip

The last application concerns a functor $K$ satisfying 
conditions $(5.1)$--$(5.4)$, which are  weaker
than the conditions describing  the invariant vectors functor.
\medskip

\noindent{\bf 11.4 Theorem} {\sl Let $G=({\cal A},\Delta)$ be a compact
quantum group, and
$\zeta:{\cal C}\to {\cal C}\otimes{\cal A}$ a nondegenerate ergodic
$G$--action on a
unital $C^*$--algebra ${\cal C}$ with associated spectral functor
$\ov{L}$. Then the following properties are
equivalent:
\begin{description}
\item{\rm a)} ${\cal C}_{sp}$ has a $^*$--character,
\item{\rm b)} there is a subfunctor $K$ of the embedding functor $H$
satisfying properties $(5.1)$--$(5.4)$ and a quasitensor unitary natural
transformation from $\ov{L}$ to $K$,
\item{\rm c)} there is a compact quantum subgroup $K$ of $G$ and a 
faithful $^*$--homomorphism $\phi:{\cal C}_{sp}\to {\cal A}^\delta$
intertwining $\zeta$ with the $G$--action on the compact quantum quotient
space $K\backslash G$.
\end{description}}\medskip

\noindent{\bf Proof} We first show that a) implies b). Let $\chi$ be a
$^*$--character of ${\cal C}_{sp}$, and define, for $u\in\text{Rep}(G)$,
the map $\eta_u:\ov{L_u}\to H_u$ by
$$\eta_u(\ov{T})=\sum_i\chi(T(e_i)^*)e_i,$$
with $e_i$ an orthonormal basis of $H_u$. This map is an isometry, as
$$(\eta_u(\ov{T}),\eta_u(\ov{T'}))=\sum_i\ov{\chi(T(e_i)^*)}\chi(T'(e_i)^*)=
\chi(\sum_iT(e_i)T'(e_i)^*)=(\ov{T},\ov{T'}).$$ Actually $\eta$ is a
quasitensor natural transformation from $\ov{L}$ to $H$, as
for $A\in(u,v)$, $\ov{T}\in\ov{L_u}$,
$$\eta_v(\ov{L}_A(\ov{T}))=\eta_v(\ov{TA^*})=\sum_j\chi((TA^*(f_j))^*)f_j=$$
$$\sum_{i,j}\chi(T(e_i)^*)(f_j,
Ae_i)f_j=\sum_i\chi(T(e_i)^*)Ae_i=H_A(\eta_u(\ov{T})).$$
and, for $\overline{T}\in \ov{L_u}$, $\ov{T'}\in\ov{L_v}$,
$$\eta_{u\otimes v}(\ov{T}\otimes\ov{T'})=\eta_{u\otimes v}(\ov{T\otimes
T'})=$$
$$\sum_{i,j}\chi(T\otimes T'(e_i\otimes f_j)^*)e_i\otimes
f_j=\sum_{i,j}\chi((T(e_i)T'(f_j))^*)e_i\otimes
f_j=$$
$$\sum_{i}\chi(T(e_i)^*)e_i\otimes\sum_j\chi(T'(f_j)^*)f_j=
\eta_u(\ov{T})\otimes\eta_v(\ov{T'}).$$
It follows that the functor $K_u:=\eta_u(\ov{L_u})$,
$K_A:=A\upharpoonright_{K_u}$, for $A\in(u,v)$ and $u,v\in\text{Rep}(G)$
is a quasitensor $^*$--subfunctor of $H$. Therefore $K$ satisfies
properties $(5.1)$, $(5.2)$ and $(5.4)$. We are left to show that $(5.3)$ 
holds as well. Consider, for $S\in{L_v}$, the map
$\hat{S}: \ov{L_{u\otimes v}}\to \ov{L_u}$ defined in $(7.1)$.
A straightforward computation shows that 
$$\eta_u\circ\hat{S}=r_{\eta_v(\ov{S})}^*\circ\eta_{u\otimes v}.$$
Therefore 
$$r_{K_v}^*K_{u\otimes v}\subset\eta_u(\ov{L_u})=K_u,$$
and this is relation $(5.3)$.
We next show that b) implies c). Thanks to Lemma 4.1, $K$ is a quasitensor 
subfunctor of $H$. By Prop. 9.4 there is a faithful $^*$--homomorphism
$\phi:
{\cal C}^{\ov L}\to{\cal C}^{K}$ intertwining the $G$--actions.
Now, ${\cal C}^K$, regarded as  the $^*$--algebraic $G$--system defined by
$K$,
is $^*$--isomorphic to the $^*$--algebraic $G$--system
defined by the linear span of the coefficients $u_{k,\varphi}$, 
with $k\in K_u$, $\varphi\in H_u$, $u\in\text{Rep}(G)$, with the
restricted $G$--action, thanks to 
the $^*$--algebraic structure of ${\cal A}$ recalled at the end of
subsection 2.1. 
This is in turn a $^*$--subsystem of some quantum quotient space
$K\backslash G$, by Theorem 5.1. On the other hand the system $({\cal
C}^{\ov
L},\eta_{\ov{L}})$ is in turn $^*$--isomorphic to $({\cal
C}_{sp},\zeta)$, and the proof is now complete.
We are left to show that c) implies a).
The range of $\phi$ must be contained in the spectral $^*$--subalgebra
of ${\cal A}^\delta$ because of the intertwining relation between the
$G$--actions. Since ${\cal A}^\delta_{sp}$ is contained in the dense
$^*$--subalgebra of ${\cal A}$ generated by the matrix coefficients,
$u_{\varphi,\psi}$, we can define a $^*$--character $\chi$ on ${\cal
C}_{sp}$ simply by composing $\phi$ with the counit $e$ of $G$.
\medskip

\noindent{\it Remark} Under the assumptions of Theorem 11.4,
if one has an everywhere defined $^*$--character on ${\cal C}$
and if the action $\zeta$ and the Haar measure of $G$ are faithful,
one can construct a faithful embedding of the whole system $({\cal
C},\zeta)$ into a compact quantum quotient space $K\backslash G$, as shown
in Theorem 7.4 in \cite{P2}.
\medskip

\noindent{\bf Acknowledgements} The authors are very grateful to S.
Doplicher for many discussions during the preparation of the manuscript.
C.P. would also like to thank S. Vaes for drawing our attention to
\cite{BRV}. 

This research was supported by the European network
`Quantum Spaces - Noncommutative Geometry' HPRN-CT-2002-00280.

\end{section}

\end{document}